\documentclass[11pt]{article}

\usepackage{amsmath,amssymb,amsthm,mathtools}
\usepackage{geometry}
\usepackage{enumitem}
\usepackage{microtype}
\usepackage[hidelinks]{hyperref}

\numberwithin{equation}{section}

\newtheorem{theorem}{Theorem}[section]
\newtheorem{proposition}[theorem]{Proposition}
\newtheorem{lemma}[theorem]{Lemma}
\newtheorem{corollary}[theorem]{Corollary}
\theoremstyle{remark}
\newtheorem{remark}[theorem]{Remark}

\newcommand{\im}{{\rm Im}\,}
\newcommand{\re}{{\rm Re}\,}
\newcommand{\dist}{\hbox{ \rm dist}}

\newcommand {\C} {\mathbb C}
\newcommand {\R} {\mathbb R}
\newcommand {\Z} {\mathbb Z}

\newcommand{\N}{\mathbb N}

\title{Analytic Tools and Weighted Composition Operators on Paley--Wiener Spaces}

\author{
Pham Viet Hai\thanks{Faculty of Mathematics and Informatics,
Hanoi University of Science and Technology,
1 Dai Co Viet, Hanoi, Vietnam.
E-mail: \texttt{hai.phamviet@hust.edu.vn}}
\and
Vu Nhat Huy\thanks{Faculty of Mathematics, Mechanics and Informatics,
University of Science, Vietnam National University,
Hanoi, Vietnam.
E-mail: \texttt{nhat\_huy85@yahoo.com}}
\and
Pham Trong Tien\thanks{Faculty of Mathematics, Mechanics and Informatics,
University of Science, Vietnam National University,
Hanoi, Vietnam.
E-mail: \texttt{phamtien@vnu.edu.vn}}
}

\date{}

\begin{document}

\maketitle

\begin{abstract}

\smallskip
\noindent
We develop a scale-explicit analytic framework for operator-theoretic problems on the Paley--Wiener spaces \(B_\sigma^p\), \(\sigma>0\), throughout the full range \(0<p\leq\infty\). The main device is a family of localized Paley--Wiener packets obtained from powers of reproducing kernels. These packets yield sharp-order point-evaluation estimates, quantitative separation of distinct point-evaluation functionals, localized test functions, and direct scale-explicit sampling estimates and Carleson measure characterizations.

\smallskip
\noindent
As a principal application of this framework, we develop a comprehensive theory of weighted composition operators
\[
W_{\psi,\varphi}:B_\sigma^p\longrightarrow B_\tau^q,
\qquad
0<p,q\leq\infty,\quad \sigma,\tau>0.
\]
Boundedness forces every nonconstant symbol to be affine with real slope. A sharp dichotomy appears according to the relative size of \(p\) and \(q\): when \(p\leq q\), the weight is governed by the residual bandwidth \(\tau-\sigma|a|\), whereas when \(q<p\), the weight must also supply the missing integrability, with exponent \(\frac{pq}{p-q}\), interpreted as \(q\) when \(p=\infty\). In the latter range, boundedness and compactness coincide. We also obtain compactness criteria and essential-norm estimates in the range \(p\leq q\), together with quantitative results for composition and multiplication operators.

\smallskip
\noindent
Beyond individual operators, we establish a no-cancellation principle for differences of weighted composition operators: for distinct symbols, a bounded (respectively, compact) difference can occur only if both summands are bounded (respectively, compact).
For composition operators, the fixed-slope classes are precisely the components when \(p\leq q\), while no bounded composition operators exist otherwise. By contrast, for nonzero weighted composition operators, the whole space is path connected when \(q<p\), whereas in the complementary range it has exactly three connected components, all of which are path connected.
\end{abstract}



\clearpage
\tableofcontents
\clearpage

\section{Introduction}
\label{sec:introduction}

For \(\sigma>0\) and \(0<p\leq\infty\), let \(B_\sigma^p\) denote the
Paley--Wiener space of entire functions of exponential type at most \(\sigma\)
whose restrictions to the real line belong to \(L^p(\R)\). For
\(1\leq p\leq\infty\), the space \(B_\sigma^p\), equipped with the norm
\(\|\cdot\|_{L^p}\), is a Banach space, while for \(0<p<1\) it is a
quasi-Banach space with quasi-norm \(\|\cdot\|_{L^p}\). We denote by
\((B_\sigma^p)^*\) its continuous dual, endowed with the operator norm
\[
    \|\Lambda\|_{(B_\sigma^p)^*}
    :=
    \sup_{\|f\|_{L^p}\leq1}|\Lambda(f)|.
\]
For convenience, we also denote by \(B_0^\infty\) the space of all entire
functions of exponential type zero that are bounded on \(\R\), equipped
with the norm \(\|\cdot\|_{L^\infty}\). These spaces form a natural
meeting point of complex analysis, Fourier analysis, and operator theory.

\smallskip
\noindent
Operator theory on Paley--Wiener spaces presents several features that
distinguish it from the corresponding theories on Hardy, Bergman, and Fock
spaces. For \(0<p<\infty\), the spaces \(B_\sigma^p\) contain no nonzero
polynomials, while for \(p=\infty\) the only polynomial elements are
constants. Thus monomial- or Taylor-based arguments are unavailable
throughout the scale considered here. Moreover, in the Hilbert case the
classical Paley--Wiener reproducing kernel is a sinc kernel and, for each
fixed center, has infinitely many zeros. Consequently, arguments relying on
nonvanishing kernels or direct kernel-ratio methods do not transfer in a
straightforward way. At the same time, working in the full range
\(0<p\leq\infty\), including the quasi-Banach regime, requires a careful
combination of entire-function estimates, Fourier localization, and
real-variable discretization. These features make localization, sampling,
point-evaluation estimates, and entire-function rigidity particularly
natural tools in the Paley--Wiener setting.

\smallskip
\noindent
The first purpose of this paper is to develop an analytic framework adapted
to such operator-theoretic problems. A central feature of Paley--Wiener
analysis is the interaction between bandwidth, pointwise growth, and spatial
localization, and we organize these ingredients in a quantitative form
valid throughout \(0<p\leq\infty\). The main device is a family of
localized Paley--Wiener packets obtained from powers of reproducing
kernels. They yield sharp-order point-evaluation estimates, quantitative
separation of distinct point-evaluation functionals, localized test
functions and compactness sequences, as well as direct scale-explicit
sampling estimates and Carleson measure characterizations. While sampling, point-evaluation, and Carleson-measure questions have classical antecedents, the framework developed here organizes the tools needed for operator theory into a unified quantitative form across the full Banach and quasi-Banach Paley--Wiener scale. In particular, the quantitative separation of point evaluations developed below will be a key ingredient in the rigidity results for differences of operators.

\smallskip
\noindent
Point evaluations provide a natural link
with classical extremal questions in Paley--Wiener spaces. On the real
line, translation invariance implies that the norm of each point-evaluation
functional coincides with the norm of the embedding
\(B_\sigma^p\hookrightarrow B_\sigma^\infty\). Exact constants,
extremizers, and related off-axis estimates have been studied in
considerable detail; see, for instance,
\cite{BrevigChirreOrtegaCerdaSeip2024,Korevaar1949}.
For our operator-theoretic purposes, the packet
construction gives a simple scale-explicit two-sided estimate together with
the localization properties needed later. It also leads to a quantitative
estimate for the distance between two distinct point-evaluation
functionals modulo scalar multiples, a form of separation particularly
suited to the study of differences of weighted composition operators.

\smallskip
\noindent
Weighted composition operators provide a natural class on which to test
this framework. For entire functions \(\psi\) and \(\varphi\), one sets
\[
    W_{\psi,\varphi}f
    :=
    \psi\cdot(f\circ\varphi),
\]
which contains both the composition operator
\(C_\varphi\) and the multiplication operator \(M_\psi\) as the special
cases \(\psi\equiv1\) and \(\varphi=\operatorname{id}\), respectively.
Weighted composition operators have been studied extensively on classical
spaces of analytic functions, including Hardy and Bergman spaces
\cite{ContrerasHernandezDiaz2003,CuckovicZhao2004}, and on spaces of
entire functions such as Fock spaces
\cite{HaiKhoi2016,TienKhoi2019,Ueki2007}. In these settings one has
developed substantial theories of boundedness, compactness, essential
norms, differences of operators, and the topology of the corresponding
operator spaces.

\smallskip
\noindent
Differences of operators play a central role in such topological questions. In the Hardy-space setting, Berkson initiated the study of
isolation phenomena for composition operators, and Shapiro and Sundberg
subsequently related component questions to compact differences; see
\cite{Berkson1981,ShapiroSundberg1990}. Compact differences have since
been studied on several classical spaces of analytic functions, notably
on weighted Bergman spaces by Moorhouse \cite{Moorhouse2005} and, in the
weighted-composition setting, in \cite{ChoeChoiKooYang2020}. For Fock
spaces, compact differences and the topological structure of the spaces
of composition and weighted composition operators admit particularly
explicit descriptions; see, for example,
\cite{Dai2015,TienKhoiDifferences2019,TienKhoi2019}. Nevertheless, even
on classical spaces of analytic functions on the unit disc, the
topological structure problem can be delicate, and complete descriptions
of path-connected components remain unavailable in several settings; see,
for example,
\cite{ChoeChoiKooPark2022,GallardoGutierrezGonzalezNieminenSaksman2008,Tien2023}.

\smallskip
\noindent
By comparison, the literature on composition and weighted composition
operators on Paley--Wiener spaces is much more limited. Chac\'on,
Chac\'on and Gim\'enez \cite{ChaconChaconGimenez2007} initiated the study
of composition operators on the classical Paley--Wiener space
\(B_\pi^2\). They showed that a nonconstant entire symbol induces a
bounded composition operator precisely when
\[
    \varphi(z)=az+b,
    \qquad
    a\in\R,\quad 0<|a|\leq1,
\]
and proved that no such composition operator is compact. They also
considered composition operators on \(B_\pi^p\) and developed Carleson
measure criteria for \(1\leq p\leq q<\infty\), with the compactness
characterization in the range \(1<p\leq q<\infty\).

\smallskip
\noindent
Mukherjee, Jafari and McInroy \cite{MukherjeeJafariMcInroy2011} subsequently considered a weighted-composition problem in the Hilbert space \(B_a^2\). Assuming the symbol to be affine with nonzero real slope, they related the bandwidth of the image to Fourier-support conditions on the weight. In particular, the affine form of the symbol is assumed rather than derived from boundedness. Their result therefore provides an early connection between weighted composition operators and bandwidth control, but does not address the general boundedness and compactness theory across arbitrary exponents and bandwidths.

\smallskip
\noindent
More recently, Samanta and Sarkar \cite{SamantaSarkar2023} extended the boundedness classification for composition operators to several-variable Paley--Wiener spaces in the Hilbert setting. Other recent work has focused on operator-theoretic properties of composition operators on \(B_\sigma^2\). In particular, Hai, Noor and Severiano \cite{HaiNoorSeveriano2025} characterized cyclicity of bounded composition operators on these spaces, with complex symmetry playing an important role in their analysis.

\smallskip
\noindent
Taken together, these works show that the theory on Paley--Wiener spaces is still comparatively sparse: the existing results concern mainly composition operators, Hilbert-space settings, or particular questions of bandwidth control for weighted composition operators. To the best of our knowledge, there has been no systematic treatment of
\[
W_{\psi,\varphi}:B_\sigma^p\longrightarrow B_\tau^q
\]
for arbitrary bandwidths \(\sigma,\tau>0\) and throughout the full range \(0<p,q\leq\infty\), nor a corresponding theory of differences and of the topological structure of these operator spaces.

\smallskip
\noindent
The second purpose of the present paper, and the principal application of
the analytic framework developed above, is to provide such a theory. A
first structural result shows that boundedness forces every nonconstant
symbol to be affine with real slope,
\[
    \varphi(z)=az+b,\qquad a\in\R.
\]
Once this rigidity is established, two mechanisms govern the operator: the symbol determines the bandwidth available to the weight, while the weight may also have to compensate for a loss of real-line integrability. Their roles separate sharply according to the relative size of \(p\) and \(q\).

\smallskip
\noindent
If \(p\leq q\), Nikol'skii's inequality supplies the required gain of
integrability, and the weight is governed by the residual bandwidth
\[
    \tau-\sigma|a|.
\]
In contrast, when \(q<p\), the weight itself must supply the missing
integrability, leading to the exponent
\[
    \frac{pq}{p-q}.
\]
Thus the passage from the Hilbert-space setting to the full two-exponent
scale is not merely a formal extension: it reveals an additional
integrability mechanism that is invisible when \(p=q=2\). In the range
\(q<p\), this mechanism also produces the striking equivalence between
boundedness and compactness. In the complementary range \(p\leq q\), we
obtain a complete compactness criterion and two-sided essential-norm
estimates. We also derive sharp-order norm estimates for affine
composition operators and corresponding quantitative results for
multiplication operators.

\smallskip
\noindent
Beyond individual weighted composition operators, we study their
differences
\[
    W_{\psi,\varphi}-W_{\vartheta,\phi}.
\]
Here a new rigidity phenomenon appears. When \(\varphi\not\equiv\phi\), a bounded (respectively, compact) difference requires both summands to be bounded (respectively, compact); in other words, distinct symbols admit no cancellation. When the symbols coincide, the problem reduces to the single weighted composition operator \(W_{\psi-\vartheta,\varphi}\). In the range \(q<p\), boundedness and
compactness of the difference are again equivalent. The quantitative
separation of point evaluations developed in the analytic part of the
paper is one of the principal tools behind this rigidity.

\smallskip
\noindent
Finally, these results allow us to determine the topological structure of the corresponding operator spaces in the natural operator norm or quasi-norm topology. For composition operators, the fixed-slope classes are precisely the connected and path-connected components when \(p\leq q\), while no bounded composition operators exist otherwise. By contrast, for nonzero weighted composition operators, the whole space is path connected when \(q<p\), whereas in the complementary range it has exactly three connected components, each of which is path connected. These correspond to the strict bandwidth region
\(
    \sigma|a|<\tau
\)
and the two endpoint slopes
\[
    a=\frac{\tau}{\sigma}
    \qquad\text{and}\qquad
    a=-\frac{\tau}{\sigma}.
\]
In all cases the space of nonzero weighted composition operators has no isolated points.

\smallskip
\noindent
The remainder of the paper is organized as follows.
Section~\ref{sec:analytic-tools} develops the scale-explicit analytic
tools used throughout the paper, including vertical translations,
Paley--Wiener packets, sharp-order point-evaluation estimates,
quantitative separation of point evaluations, sampling estimates, and
Carleson measure characterizations.
Section~\ref{sec-wco} studies boundedness, compactness,
and essential norms of weighted composition operators, together with the
corresponding results for composition and multiplication operators.
Section~\ref{sec:differences} is devoted to differences of weighted
composition operators and establishes the no-cancellation phenomena for
distinct symbols in both exponent ranges.
Finally, Section~\ref{sec-topological-structure} determines the connected
and path-connected components of the spaces of composition operators and
nonzero weighted composition operators.

\bigskip
\noindent
\textbf{Notation.}
We use the convention \(\frac{1}{\infty}=0\). Whenever
\(0<q<p\leq\infty\), we also adopt the convention
\[
    \frac{pq}{p-q}:=q
    \qquad\text{when } p=\infty.
\]
We write
\(A\lesssim_{\alpha_1,\dots,\alpha_k}B\) if \(A\leq CB\) for a
constant \(C>0\) depending only on the indicated parameters
\(\alpha_1,\dots,\alpha_k\), and write
\(A\simeq_{\alpha_1,\dots,\alpha_k}B\) if
\(A\lesssim_{\alpha_1,\dots,\alpha_k}B\) and
\(B\lesssim_{\alpha_1,\dots,\alpha_k}A\). In statements of results, the
dependence of implicit constants is specified explicitly in the text. In
proofs, parameter dependence is indicated by subscripts whenever
necessary; an unsubscripted implicit constant is absolute.

\smallskip
\noindent
For an entire function \(f\) of finite exponential type, we write
\(\operatorname{type}(f)\) for its exponential type.

\section{Analytic tools in Paley--Wiener spaces}
\label{sec:analytic-tools}

The purpose of this section is to develop a scale-explicit collection
of analytic tools for Paley--Wiener spaces. We begin with scaling, vertical translations, and Nikol'skii inequalities.
We then introduce localized Paley--Wiener packets and use them to study
point evaluations, their quantitative separation, and spatial decay. A second family of localized packets leads to sampling
and discretization estimates on the real line, which in turn provide the
basis for the Carleson and vanishing Carleson measure characterizations
below.

\subsection{Scaling and basic estimates}
\label{subsec:scaling}

We first record the scaling normalization that will be used throughout this
section. For \(\sigma>0\) and \(0<p\leq\infty\), define
\begin{equation}\label{eq:scaling-formula}
    (\mathcal U_\sigma f)(z)
    :=
    \left(\frac{\pi}{\sigma}\right)^{\frac{1}{p}}
    f\left(\frac{\pi}{\sigma}z\right).
\end{equation}
Then
\(
    \mathcal U_\sigma:B_\sigma^p\longrightarrow B_\pi^p
\)
is an isometric isomorphism, with inverse
\[
    (\mathcal U_\sigma^{-1}g)(z)
    =
    \left(\frac{\sigma}{\pi}\right)^{\frac{1}{p}}
    g\left(\frac{\sigma}{\pi}z\right).
\]
Thus estimates at bandwidth \(\pi\) transfer to arbitrary
\(\sigma>0\) with the corresponding scale factors made explicit.

\begin{lemma}[Vertical translations]
\label{lem:vertical-translation}
Let \(\sigma>0\) and \(0<p\leq\infty\). For every \(y \in \R\),  the vertical translation operator
\[
    T_yf(z):=f(z+iy)
\]
maps \(B_\sigma^p\) boundedly into itself and satisfies
\[
    \|T_y\|_{B_\sigma^p\to B_\sigma^p}
    = e^{\sigma|y|}.
\]
\end{lemma}

\begin{proof}
The classical Plancherel--Pólya inequality for
\(0<p<\infty\) \cite[Theorem~4, p.~50]{Levin1996} and Bernstein's inequality for \(p=\infty\) \cite{Bernstein1923} give
\[
    \|f(\cdot+iy)\|_{L^p}
    \leq
    e^{\sigma|y|}\|f\|_{L^p}.
\]
Thus
\[
    \|T_y\|_{B_\sigma^p\to B_\sigma^p}
    \leq e^{\sigma|y|}.
\]

\smallskip
\noindent
For the reverse estimate, fix \(0<\varepsilon< \sigma\) and choose
\(k\in\mathbb N\) such that \(k> \frac{1}{p}\). Set
\[
    h_\varepsilon(z)
    :=
    \left(
        \frac{\sin(\varepsilon z/k)}
             {\varepsilon z/k}
    \right)^{k},
    \qquad h_\varepsilon(0):=1.
\]
This sinc-power is the prototype of the Paley--Wiener packets introduced
in Subsection~\ref{subsec:packets}.
The function \(h_\varepsilon\) has exponential type at most
\(\varepsilon\), and
\[
    |h_\varepsilon(x)|
    \lesssim_p
    (1+\varepsilon|x|)^{-k},
    \qquad x\in\R.
\]
Since \(k>\frac{1}{p}\), we have \(h_\varepsilon\in B_\varepsilon^p\).
Define
\[
    f_\varepsilon(z)
    :=
    \begin{cases}
        e^{-i(\sigma-\varepsilon)z}h_\varepsilon(z),
            & y\geq0,\\[2mm]
        e^{i(\sigma-\varepsilon)z}h_\varepsilon(z),
            & y<0.
    \end{cases}
\]
Thus \(f_\varepsilon\in B_\sigma^p\) and
\(
    \|f_\varepsilon\|_{L^p}
    =
    \|h_\varepsilon\|_{L^p}.
\)
Moreover,
\[
    \|T_yf_\varepsilon\|_{L^p}
    =
    e^{(\sigma-\varepsilon)|y|}
    \|h_\varepsilon(\cdot+iy)\|_{L^p}.
\]
Applying the upper estimate to \(h_\varepsilon(\cdot+iy)\in
B_\varepsilon^p\) with translation \(-iy\), we obtain
\[
    \frac{\|T_yf_\varepsilon\|_{L^p}}
         {\|f_\varepsilon\|_{L^p}}
    \geq
    e^{(\sigma-2\varepsilon)|y|}.
\]
Letting \(\varepsilon\to0^+\) gives the reverse inequality.
\end{proof}

\begin{lemma}[Nikol'skii inequality]
\label{lem:nikolskii}
Let \(\sigma>0\) and \(0<p\leq q\leq\infty\). Then
\[
    \|f\|_{L^q}
    \lesssim
    \sigma^{\frac{1}{p}-\frac{1}{q}}\|f\|_{L^p},
    \qquad f\in B_\sigma^p,
\]
where the implicit constant depends only on \(p\) and \(q\). In particular,
\(
    B_\sigma^p\subset B_\sigma^q
\)
continuously with
\[
    \|I\|_{B_\sigma^p\to B_\sigma^q}
    =
    \left(\frac{\sigma}{\pi}\right)^{\frac{1}{p}-\frac{1}{q}}
    \|I\|_{B_\pi^p\to B_\pi^q}.
\]
\end{lemma}

\begin{proof}
Nikol'skii's inequality at bandwidth \(\pi\) is classical; see
\cite[Corollary~1]{NesselWilmes1978}. By the scaling formula \eqref{eq:scaling-formula},
\[
    \frac{\|f\|_{L^q}}{\|f\|_{L^p}}
    =
    \left(\frac{\sigma}{\pi}\right)^{\frac{1}{p}-\frac{1}{q}}
    \frac{\|\mathcal U_\sigma f\|_{L^q}}
         {\|\mathcal U_\sigma f\|_{L^p}},
    \qquad 0\neq f\in B_\sigma^p.
\]
Since
\(
    \mathcal U_\sigma:B_\sigma^p\to B_\pi^p
\)
is an isometric isomorphism, taking the supremum over
\(0\neq f\in B_\sigma^p\) gives the desired equality.
\end{proof}

\medskip
\noindent
Combining Lemmas~\ref{lem:vertical-translation} and
\ref{lem:nikolskii}, we obtain
\begin{equation}\label{eq:vertical-translation}
\|f(\cdot+iy)\|_{L^q}
\lesssim
\sigma^{\frac{1}{p}-\frac{1}{q}}e^{\sigma|y|}
\|f\|_{L^p},
\qquad
f\in B_\sigma^p,\quad y\in\mathbb R,
\end{equation}
whenever \(0<p\leq q\leq\infty\) with the implicit constant depending only on \(p\) and \(q\).

\subsection{Paley--Wiener packets and point evaluations}
\label{subsec:packets}

We now turn to the localization mechanism underlying the subsequent
point-evaluation and sampling estimates.
For \(\sigma>0\) and \(w\in\C\), let
\[
    K_{\sigma,w}(z)
    :=
    \frac{\sin \sigma(z-\overline w)}
         {\pi(z-\overline w)},
    \qquad z\in\C,
\]
with the usual continuous interpretation at \(z=\overline w\).  This is the
reproducing kernel of \(B_\sigma^2\).
For \(m\in\N\) and \(w\in\C\), define
\[
    G_{\sigma,m,w}(z)
    :=
    \bigl(K_{\sigma/m,w}(z)\bigr)^m.
\]
The reduced bandwidth \(\frac{\sigma}{m}\) compensates for the \(m\)-th power,
so that \(G_{\sigma,m,w}\) has exponential type at most \(\sigma\).
These packets combine a large peak at \(w\) with quantitative spatial
decay.

\begin{lemma}[Paley--Wiener packet estimates]
\label{lem:packet-estimates}
Let \(\sigma > 0\), \(0<p\leq\infty\), \(w \in \C\), and let \(m\in\N\) satisfy
\(m>\frac1p\).
Then \(G_{\sigma,m,w}\in B_\sigma^p\) and
\[
    \|G_{\sigma,m,w}\|_{L^p}
    \simeq
    \sigma^{m-\frac{1}{p}}
    \frac{e^{\sigma |\im w|}}
         {(1+\sigma |\im w|)^{m-\frac{1}{p}}}
\qquad \text{and} \qquad
    |G_{\sigma,m,w}(w)|
    \simeq
    \sigma^m
    \frac{e^{2\sigma |\im w|}}
         {(1+\sigma |\im w|)^m},
\]
where the implicit constants depend only on \(m\) and \(p\).
\end{lemma}

\begin{proof}
We first record the corresponding kernel estimates. Let \(\rho>0\),
\(1<s\leq\infty\).
By translation invariance of the \(L^s(\mathbb R)\)-norm and the change of
variables \(x=\rho t\), for \(1<s<\infty\),
\[
    \|K_{\rho,w}\|_{L^s}^s
    =
    \int_{\mathbb R}
    \left|
        \frac{\sin\rho(t+i|\im w|)}
             {\pi(t+i|\im w|)}
    \right|^s\,dt =
    \frac{\rho^{s-1}}{\pi^s}
    \int_{\mathbb R}
    \left|
        \frac{\sin(x+i\rho|\im w| )}
             {x+i\rho |\im w|}
    \right|^s\,dx.
\]
We use the identity
\[
    |\sin(x+i \rho |\im w|)|^2
    =
    \sin^2 x+\sinh^2 (\rho |\im w|).
\]
There are two cases of $\rho |\im w|$.

\smallskip
\noindent
If \(0\leq \rho |\im w| \leq1\), then
\[
    \int_{\mathbb R}
    \left|
        \frac{\sin(x+i \rho |\im w|)}
             {x+i\rho |\im w|}
    \right|^s\,dx
    \simeq_s 1.
\]
Indeed, the entire extension of \(\frac{\sin\zeta}{\zeta}\) is bounded above
and bounded away from zero on the square
\[
    \{\zeta\in\C:|\operatorname{Re}\zeta|\leq1,\,
    |\operatorname{Im}\zeta|\leq1\},
\]
while, for \(|x|\geq1\), the above identity gives
\[
    \left|
        \frac{\sin(x+i\rho|\im w|)}
             {x+i\rho|\im w|}
    \right|
    \lesssim |x|^{-1}.
\]
Since \(s>1\), the assertion follows.

\smallskip
\noindent
If \(\rho |\im w|\geq1\), then
\[
    |\sin(x+i\rho |\im w|)| 
    \simeq e^{\rho |\im w|},
    \qquad x\in\mathbb R.
\]
Hence
\[
    \int_{\mathbb R}
    \left|
        \frac{\sin(x+i \rho |\im w|)}
             {x+i\rho |\im w|}
    \right|^s\,dx
    \simeq_s
    e^{s \rho |\im w|}
    \int_{\mathbb R}
    \frac{dx}{(x^2+\rho^2 |\im w|^2)^{\frac{s}{2}}} =
    e^{s\rho |\im w|} (\rho |\im w|)^{1-s}
    \int_{\mathbb R}
    \frac{dt}{(1+t^2)^{\frac{s}{2}}}.
\]
Combining the two cases gives
\begin{equation}\label{eq:kernel-norm-estimate}
    \|K_{\rho,w}\|_{L^s}
    \simeq_s
    \rho^{1-\frac1s}
    \frac{e^{\rho |\im w|}}
         {(1+\rho |\im w|)^{1-\frac1s}},
    \qquad 1<s<\infty.
\end{equation}
The same estimate for \(s=\infty\) follows directly from the preceding
two cases. Thus \eqref{eq:kernel-norm-estimate} holds for
\(1<s\leq\infty\). Moreover,
\[
    K_{\rho,w}(w)
    =
    \frac{\sin\rho(w-\overline w)}
         {\pi(w-\overline w)}
    =
    \frac{\sinh(2\rho\operatorname{Im}w)}
         {2\pi\operatorname{Im}w},
\]
with the limiting value \(\frac{\rho}{\pi}\) when
\(\operatorname{Im}w=0\). Therefore
\begin{equation}\label{eq:kernel-peak-estimate}
    |K_{\rho,w}(w)|
    \simeq
    \rho
    \frac{e^{2\rho |\im w|}}
         {1+\rho |\im w|}.
\end{equation}
Now set
\(
    \rho:=\frac{\sigma}{m}.
\)
Since \(K_{\rho,w}\) has exponential type at most \(\rho\), the function $G_{\sigma,m,w}
    =
    (K_{\rho,w})^m$ has exponential type at most \(\sigma\). If \(p<\infty\), then
\(mp>1\), and \eqref{eq:kernel-norm-estimate} with \(s=mp\) gives
\(K_{\rho,w}\in L^{mp}(\mathbb R)\). Hence
\(G_{\sigma,m,w}\in B_\sigma^p\) and
\[
    \|G_{\sigma,m,w}\|_{L^p}
    =
    \|K_{\rho,w}\|_{L^{mp}}^m \simeq_{p,m} 
    \rho^{m-\frac1p}
    \frac{e^{m\rho |\im w|}}
         {(1+\rho |\im w|)^{m-\frac1p}}.
\]
The same identity is immediate when \(p=\infty\). From this, it follows that
\[
    \|G_{\sigma,m,w}\|_{L^p}
    \simeq_{p,m}
    \sigma^{m-\frac1p}
    \frac{e^{\sigma |\im w|}}
         {(1+\sigma |\im w|)^{m-\frac1p}}.
\]
Finally, by \eqref{eq:kernel-peak-estimate},
\[
    |G_{\sigma,m,w}(w)|
    =
    |K_{\rho,w}(w)|^m \simeq
    \rho^m
    \frac{e^{2m\rho |\im w|}}
         {(1+\rho |\im w|)^m} \\
    \simeq_m
    \sigma^m
    \frac{e^{2\sigma |\im w|}}
         {(1+\sigma |\im w|)^m}.
\]
\end{proof}

\begin{lemma}[Packet decay on the real line]
\label{lem:packet-decay}
Under the assumptions of Lemma~\ref{lem:packet-estimates},
\[
    |G_{\sigma,m,w}(x)|
    \lesssim
    \sigma^m e^{\sigma|\im w|}
    \bigl(1+\sigma|x-\re w|\bigr)^{-m},
    \qquad x\in\R,
\]
where the implicit constant depends only on \(m\).
In particular, for fixed \(\sigma>0\), if \((w_j)_{j\in \N}\subset\C\) satisfies
\[
    |\re w_j|\longrightarrow\infty \quad \text{as} \quad j \to \infty,
    \qquad
    \sup_{j \in \N} |\im w_j|<\infty,
\]
then the normalized packets
\[
    \widetilde G_{\sigma,m,w_j}
    :=
    \frac{G_{\sigma,m,w_j}}
         {\|G_{\sigma,m,w_j}\|_{L^p}}
\]
converge to \(0\) uniformly on compact subsets of \(\C\).
\end{lemma}

\begin{proof}
Using the elementary bound
\begin{equation} \label{eq:sin-est}
\left|\frac{\sin\zeta}{\zeta}\right|
    \lesssim
    \frac{e^{|\im\zeta|}}{1+|\zeta|},
\end{equation}
we obtain, for \(x\in\R\),
\[
    |K_{\rho,w}(x)|
    =
    \frac{\rho}{\pi}
    \left|
        \frac{\sin\bigl(\rho(x-\overline w)\bigr)}
             {\rho(x-\overline w)}
    \right| \lesssim
    \frac{\rho e^{\rho|\im w|}}{1+\rho|x-\re w|}.
\]
Hence, with \(\rho = \frac{\sigma}{m}\),
\[
    |G_{\sigma,m,w}(x)|
    =
    |K_{\rho,w}(x)|^m \lesssim \frac{\rho^m e^{m\rho|\im w|}}{\bigl(1+\rho|x-\re w|\bigr)^m} \lesssim_m
    \frac{\sigma^m e^{\sigma|\im w|}}{
    \bigl(1+\sigma|x-\re w|\bigr)^{m}}.
\]

\smallskip
\noindent
For the final assertion, the same argument gives, for \(z\in\C\),
\[
    |G_{\sigma,m,w}(z)|
    \lesssim_m
    \frac{\sigma^m
    e^{\sigma(|\im w|+|\im z|)}}{
    \bigl(1+\sigma|\re z-\re w|\bigr)^{m}}.
\]
On the other hand, Lemma~\ref{lem:packet-estimates} yields
\[
    \|G_{\sigma,m,w}\|_{L^p}
    \simeq_{p,m}
    \sigma^{m-\frac1p}
    \frac{e^{\sigma|\im w|}}
         {(1+\sigma|\im w|)^{m-\frac1p}}.
\]
Thus, if \(K\subset\C\) is compact and
\(M :=\sup_j|\im w_j|<\infty\), then
\[
    \sup_{z\in K}
    |\widetilde G_{\sigma,m,w_j}(z)|
    \lesssim_{p,m,\sigma,K,M}
    \bigl(
        1+\sigma \dist(\re w_j,\re K)
    \bigr)^{-m}.
\]
Since \(|\re w_j|\to\infty\), the right-hand side tends to \(0\).
\end{proof}

\noindent
The preceding estimates show that the packets simultaneously control
the size of evaluation at their centers and decay away from those
centers. We first apply this to point evaluations. For later use, we also set
\[
    \Phi_{\sigma,p}(t)
    :=
    \frac{\sigma^{\frac{1}{p}}e^{\sigma t}}
         {(1+\sigma t)^{\frac{1}{p}}},
    \qquad t\geq0.
\]

\begin{theorem}[Sharp-order point evaluations]
\label{thm:point-evaluations}
Let \(\sigma>0\) and \(0<p\leq\infty\). Then, for every \(w\in\C\),
the point-evaluation functional
\[
    \delta_w(f):=f(w)
\]
is continuous on \(B_\sigma^p\), and
\[
    \|\delta_w\|_{(B_\sigma^p)^*}
    \simeq
    \Phi_{\sigma,p}(|\im w|),
\]
where the implicit constants depend only on \(p\).
For \(p=\infty\),
\[
    \|\delta_w\|_{(B_\sigma^\infty)^*}
    =
    e^{\sigma|\im w|}.
\]
\end{theorem}

\begin{proof}
We first consider \(0<p<\infty\).
Suppose first that \(\sigma|\im w|\leq1\). By
Lemmas~\ref{lem:vertical-translation} and~\ref{lem:nikolskii},
\[
    |f(w)|
    \leq
    \|f(\cdot+i \, \im w)\|_{L^\infty}  \lesssim_p
    \sigma^{\frac{1}{p}}\|f(\cdot+i \, \im w)\|_{L^p} \leq
    \sigma^{\frac{1}{p}}e^{\sigma|\im w|}\|f\|_{L^p}.
\]
Since \(1+\sigma|\im w|\simeq1\) in this case,
\[
    |f(w)|
    \lesssim_p
    \Phi_{\sigma,p}(|\im w|)\|f\|_{L^p}.
\]
Now suppose that \(\sigma|\im w|>1\). If \(\im w>0\), set $F(z):=e^{i\sigma z}f(z).$ By Lemma~\ref{lem:vertical-translation},
\[
    \|F(\cdot+it)\|_{L^p}
    =
    e^{-\sigma t}\|f(\cdot+it)\|_{L^p}
    \leq
    \|f\|_{L^p},
    \qquad t>0.
\]
Thus \(F\in H^p(\C_+)\) and $\|F\|_{H^p(\C_+)}
    \leq
    \|f\|_{L^p},$ where \(\C_+\) is the upper half-plane.
The standard Hardy-space point estimate
\cite[Lemma~1, p.~188]{Duren-Hp} gives
\[
    |F(w)|
    \lesssim_p
    \frac{1}{(\im w)^{\frac{1}{p}}}\|f\|_{L^p}.
\]
Consequently,
\[
    |f(w)|
    \lesssim_p
    \frac{e^{\sigma \im w}}{(\im w)^{\frac{1}{p}}}\|f\|_{L^p}
    \lesssim_p
    \Phi_{\sigma,p}(\im w)\|f\|_{L^p}.
\]
The case \(\im w<0\) is identical, using $F(z):=e^{-i\sigma z}f(z)$ in the lower half-plane. Hence
\[
    \|\delta_w\|_{(B_\sigma^p)^*}
    \lesssim_p
    \Phi_{\sigma,p}(|\im w|).
\]

\smallskip
\noindent
For the reverse estimate, choose \(m\in\N\), depending only on \(p\),
such that \(m>\frac{1}{p}\). By Lemma~\ref{lem:packet-estimates},
\[
    \|\delta_w\|_{(B_\sigma^p)^*}
    \geq
    \frac{|G_{\sigma,m,w}(w)|}
         {\|G_{\sigma,m,w}\|_{L^p}} \gtrsim_p
    \sigma^{\frac{1}{p}}
    \frac{e^{\sigma|\im w|}}
         {(1+\sigma|\im w|)^{\frac{1}{p}}} =
    \Phi_{\sigma,p}(|\im w|),
\]
which proves the assertion for \(0<p<\infty\). Finally, let \(p=\infty\). By
Lemma~\ref{lem:vertical-translation},
\[
    |f(w)|
    \leq
    \|f(\cdot+i\im w)\|_{L^\infty}
    \leq
    e^{\sigma|\im w|}\|f\|_{L^\infty},
\]
and hence
\[
    \|\delta_w\|_{(B_\sigma^\infty)^*}
    \leq
    e^{\sigma|\im w|}.
\]
Take
\[
    f(z):= \begin{cases}
        e^{-i\sigma z}, \qquad & \im w\geq0, \\
        e^{i\sigma z}, \qquad & \im w<0. 
    \end{cases}
\]
Then
\(  \|f\|_{L^\infty}=1\) and \(    |f(w)|=e^{\sigma|\im w|}\).
Therefore
\[
    \|\delta_w\|_{(B_\sigma^\infty)^*}
    =
    e^{\sigma|\im w|}.
\]
\end{proof}

\begin{remark}
\label{rem:packet-philosophy}
The packets \(G_{\sigma,m,w}\) will be used in two
complementary ways: their peak at \(w\) yields sharp lower bounds for
point evaluations, while their off-center decay provides localized test
functions and sequences converging to zero uniformly on compact subsets.
\end{remark}

\begin{remark}
For \(w\in\R\), translation invariance shows that the norm of
\(\delta_w\) coincides with the norm of the embedding
\(B_\sigma^p\hookrightarrow B_\sigma^\infty\). The corresponding
extremal problem, including the structure of extremizers, has been
studied in considerable detail; see \cite{BrevigChirreOrtegaCerdaSeip2024}.
The study of point evaluation off the real axis goes back at least to
Korevaar~\cite{Korevaar1949}; see also
\cite[Section~9.2]{BrevigChirreOrtegaCerdaSeip2024}.
Theorem~\ref{thm:point-evaluations} provides the scale-explicit
two-sided estimate needed in the present paper.
\end{remark}

\noindent
We next record a quantitative separation estimate for point-evaluation
functionals. Besides the size of a single evaluation functional, we shall
need to measure how far two evaluations are from being scalar multiples of
one another.

\begin{lemma}[Quantitative separation of point evaluations]\label{lem-two-point-separation}
Let \(\sigma>0\) and \(0<p\leq\infty\). For \(w,v\in\C\), put
\[
    d_{\sigma,p}(w,v)
    :=
    \operatorname{dist}\bigl(\delta_w,\C\delta_v\bigr)
    =
    \inf_{\lambda\in\C}
    \|\delta_w-\lambda\delta_v\|_{(B_\sigma^p)^*}.
\]
Then
\[
    d_{\sigma,p}(w,v)
    \gtrsim
    \Phi_{\sigma,p}(|\im w|)
    \frac{\sigma|w-v|}
    {1+\sigma|\im w|+\sigma|w-v|},
    \qquad w,v\in\C,
\]
where the implicit constant depends only on \(p\).
\end{lemma}

\begin{proof}
The conclusion is trivial when \(w=v\), so suppose that \(w\neq v\).
Choose \(m\in\N\), depending only on \(p\), such that $m>1+\frac1p$. Let
\[
    G_{\sigma,m,w}(z)
    =
    \bigl(K_{\sigma/m,w}(z)\bigr)^m
\]
be the packet from Lemma~\ref{lem:packet-estimates}, and define
\[
    H_{w,v}(z):=(z-v)G_{\sigma,m,w}(z).
\]
Multiplication by a polynomial does not increase exponential type. We shall
see below that \(H_{w,v}|_{\R}\in L^p(\R)\), and hence
\(H_{w,v}\in B_\sigma^p\). Indeed, for every \(\lambda\in\C\), we have
\[
    \|\delta_w-\lambda\delta_v\|_{(B_\sigma^p)^*}
    \geq\frac{|(\delta_w-\lambda\delta_v)H_{w,v}|}{\|H_{w,v}\|_{L^p}}
    =
    \frac{|w-v|\,|G_{\sigma,m,w}(w)|}
    {\|H_{w,v}\|_{L^p}}.
\]
We next estimate the denominator. Since $z-v=(z-\overline{w})+(\overline{w}-v)$, the triangle inequality, or the \(p\)-triangle inequality when \(0<p<1\),
gives
\[
    \|H_{w,v}\|_{L^p}
    \lesssim_p
    \|(\cdot-\overline{w})G_{\sigma,m,w}\|_{L^p}
    +
    |\overline{w}-v|\,\|G_{\sigma,m,w}\|_{L^p}.
\]
Set
\[
    \alpha:=\frac{\sigma}{m},
    \qquad
    \rho:=(m-1)\alpha=\frac{m-1}{m}\sigma.
\]
Then
\[
\begin{aligned}
    (z-\overline{w})G_{\sigma,m,w}(z)
    &=
    \frac{\sin\alpha(z-\overline{w})}{\pi}
    \bigl(K_{\alpha,w}(z)\bigr)^{m-1} \\
    &=
    \frac{\sin\alpha(z-\overline{w})}{\pi}
    G_{\rho,m-1,w}(z).
\end{aligned}
\]
For \(x\in\R\),
\[
    |\sin\alpha(x-\overline{w})|
    \leq
    \cosh\bigl(\alpha |\im w|\bigr)
    \leq
    e^{\alpha |\im w|}.
\]
Since \(m-1>\frac{1}{p}\), Lemma~\ref{lem:packet-estimates} gives
\[
\begin{aligned}
    \|(\cdot-\overline{w})G_{\sigma,m,w}\|_{L^p}
    &\lesssim_p
    e^{\alpha |\im w|}
    \rho^{\,m-1-\frac1p}
    \frac{e^{\rho |\im w|}}
    {\bigl(1+\rho |\im w|\bigr)^{\,m-1-\frac1p}} \\
    &\lesssim_p
    \sigma^{\,m-1-\frac1p}
    \frac{e^{\sigma |\im w|}}
    {\bigl(1+\sigma |\im w|\bigr)^{\,m-1-\frac1p}}.
\end{aligned}
\]
Also, by Lemma~\ref{lem:packet-estimates},
\[
    \|G_{\sigma,m,w}\|_{L^p}
    \lesssim_p
    \sigma^{\,m-\frac1p}
    \frac{e^{\sigma |\im w|}}
    {\bigl(1+\sigma |\im w|\bigr)^{\,m-\frac1p}}.
\]
Consequently,
\[
    \|H_{w,v}\|_{L^p}
    \lesssim_p
    \sigma^{\,m-1-\frac1p}
    \frac{e^{\sigma |\im w|}}
    {\bigl(1+\sigma |\im w|\bigr)^{\,m-\frac1p}}
    \bigl(1+\sigma |\im w|+\sigma|w-v|\bigr).
\]
In particular, \(H_{w,v}\in B_\sigma^p\).

\smallskip
\noindent
On the other hand, Lemma~\ref{lem:packet-estimates} yields
\[
    |G_{\sigma,m,w}(w)|
    \gtrsim_p
    \sigma^m
    \frac{e^{2\sigma |\im w|}}
    {\bigl(1+\sigma |\im w|\bigr)^m}.
\]
Therefore, for every \(\lambda\in\C\),
\[
\begin{aligned}
    \|\delta_w-\lambda\delta_v\|_{(B_\sigma^p)^*}
    &\gtrsim_p
    \sigma^{1+\frac1p}
    \frac{e^{\sigma |\im w|}}
    {\bigl(1+\sigma |\im w|\bigr)^{\frac1p}}\cdot
    \frac{|w-v|}
    {1+\sigma |\im w|+\sigma|w-v|} \\
    &=
    \Phi_{\sigma,p}(|\im w|)
    \frac{\sigma|w-v|}
    {1+\sigma |\im w|+\sigma|w-v|}.
\end{aligned}
\]
Taking the infimum over \(\lambda\in\C\) proves the result.
\end{proof}

\noindent
We also record the following compactness principle, which will be used
repeatedly below.

\begin{lemma}[Compactness criterion]\label{lem-compactness-criterion}
Let \(\sigma>0\), \(0<p\leq\infty\), and let \(Y\) be a quasi-Banach
space. Suppose that \(T:B_\sigma^p\to Y\) is bounded and has the following
uniqueness property: whenever \((f_j)_{j \in \N}\) is bounded in \(B_\sigma^p\),
\(f_j\to0\) uniformly on compact subsets of \(\C\), and
\(    Tf_{j_k}\longrightarrow g\)
    in \(Y\)
along a subsequence, then \(g=0\).

\smallskip
\noindent
Then \(T\) is compact if and only if $\|Tf_j\|_Y\longrightarrow0$ for every bounded sequence \((f_j)_{j \in \N}\subset B_\sigma^p\) converging to
\(0\) uniformly on compact subsets of \(\C\).
\end{lemma}

\begin{proof}
Suppose first that \(T\) is compact, and let \((f_j)_{j \in \N}\) be bounded in
\(B_\sigma^p\) with \(f_j\to0\) uniformly on compact subsets of \(\C\).
If \(\|Tf_j\|_Y\not\to0\), then some subsequence satisfies
\(\|Tf_{j_k}\|_Y\geq\varepsilon>0\). By compactness, a further
subsequence converges in \(Y\), and the uniqueness property forces its
limit to be \(0\), a contradiction.

\smallskip
\noindent
Conversely, let \((f_j)_{j \in \N}\) be bounded in \(B_\sigma^p\). By
Theorem~\ref{thm:point-evaluations}, it is locally bounded on
\(\C\); hence, after passing to a subsequence, $f_{j_k}\longrightarrow f$ uniformly on compact subsets of \(\C\). The same point-evaluation estimate
shows that \(f\) has exponential type at most \(\sigma\). Moreover, Fatou's
lemma gives \(f\in L^p(\R)\) when \(p<\infty\), while for \(p=\infty\)
the conclusion follows directly from local uniform convergence on
\(\R\). Thus \(f\in B_\sigma^p\).

\smallskip
\noindent
The sequence \(f_{j_k}-f\) is bounded in \(B_\sigma^p\) and converges
locally uniformly to \(0\). By hypothesis,
\[
    \|T(f_{j_k}-f)\|_Y\longrightarrow0,
\]
so \(Tf_{j_k}\to Tf\) in \(Y\). Hence \(T\) is compact.
\end{proof}

\subsection{Localization and sampling}
\label{subsec:localization-sampling}

We next pass from localization around a single complex point to a
real-line localization adapted to intervals of the natural spatial
scale.
For \(\rho > 0\) and \(n\in\Z\), set
\[
    I_{\rho,n}
    :=
    \left[\frac{n}{\rho},\frac{n+1}{\rho}\right).
\]
We now construct packets adapted to the cells \(I_{\rho,n}\), which
provide a convenient discretization of real-line \(L^p\) data.

\smallskip
\noindent
For each \(0<p\leq\infty\), fix once and for all an even integer
\(N_p\) such that
\[
    N_p>\max\left\{1,\frac1p\right\}.
\]
\begin{lemma}[Localized Paley--Wiener packets]
\label{lem:localized-packets}
Let  \(\sigma>0\), \(0<p\leq\infty\), and define
\[
    h_\sigma(z)
    :=
    \left(
        \frac{\sin(\sigma z/N_p)}
             {\sigma z/N_p}
    \right)^{N_p},
    \qquad h_\sigma(0):=1,
\qquad \text{and} \qquad
    h_{\sigma,n}(z)
    :=
    h_\sigma\left(
        z-\frac{n+1/2}{\sigma}
    \right),
    \qquad n\in\Z.
\]
Then the following properties hold.

\begin{enumerate}[label=\textnormal{(\alph*)}]
\item
\(h_{\sigma,n}\in B_\sigma^p\) and
\(   \|h_{\sigma,n}\|_{L^p}
    \simeq
    \sigma^{-\frac{1}{p}}\) for every \( n\in\Z\).

\item
There exists \(\eta>0\) depending only on \(p\) such that
\[
    h_{\sigma,n}(x)\geq0,
    \qquad x\in\R,
\qquad \text{and} \qquad
    h_{\sigma,n}(x)\geq\eta,
    \qquad x\in I_{\sigma,n}.
\]

\item
For fixed \(\sigma>0\),
\(
    h_{\sigma,n}\longrightarrow0
\)
uniformly on compact subsets of \(\C\) as \(|n|\to\infty\).

\item
For every finite set \(E\subset\Z\) and every complex sequence
\((c_n)_{n\in E}\),
\[
    \sum_{n\in E}c_nh_{\sigma,n}\in B_\sigma^p
\qquad \text{and} \qquad
    \left\|
        \sum_{n\in E}c_nh_{\sigma,n}
    \right\|_{L^p}
    \lesssim
    \sigma^{-\frac{1}{p}}
    \|(c_n)_{n\in E}\|_{\ell^p}.
\]
\end{enumerate}

Here the implicit constants in
\textup{(a)} and \textup{(d)} depend only on \(p\).
\end{lemma}

\begin{proof}
Put
\[
    H(z)
    :=
    \left(
        \frac{\sin(z/N_p)}
             {z/N_p}
    \right)^{N_p},
    \qquad H(0):=1.
\]
Then
\[
    h_{\sigma,n}(z)
    =
    H\left(\sigma z-n-\frac12\right).
\]
The function \(H\) has exponential type at most \(1\), and
\begin{equation}\label{eq:H}
        |H(x)|
    \lesssim_p
    (1+|x|)^{-N_p},
    \qquad x\in\R.
\end{equation}
Since \(pN_p>1\) when \(p<\infty\), it follows that
\(H\in L^p(\R)\); for \(p=\infty\), \(H\) is bounded. Thus
\(h_{\sigma,n}\) has exponential type at most \(\sigma\), and hence
belongs to \(B_\sigma^p\). Moreover, translation and scaling give
\[
    \|h_{\sigma,n}\|_{L^p}
    =
    \sigma^{-\frac{1}{p}}\|H\|_{L^p},
\]
with the usual interpretation when \(p=\infty\). This proves
\textup{(a)}.

\smallskip
\noindent
For \(x\in I_{\sigma,n}\), we have
\[
    \left|\sigma x-n-\frac12\right|\leq\frac12.
\]
Since \(N_p\) is even, \(H\geq0\) on \(\R\), and \(H\) has no zero on
\([-\frac{1}{2},\frac{1}{2}]\). Hence
\[
    h_{\sigma,n}(x)
    \geq
    \eta
    :=
    \min_{|t|\leq1/2}H(t)>0,
\]
which proves \textup{(b)}.

\smallskip
\noindent
For \textup{(c)}, let \(K\subset\C\) be compact. The elementary estimate \eqref{eq:sin-est}
gives
\[
    \sup_{z\in K}|h_{\sigma,n}(z)|
    \lesssim_{p,\sigma,K}
    (1+|n|)^{-N_p}.
\]
Thus \(h_{\sigma,n}\to0\) uniformly on \(K\) as \(|n|\to\infty\).

\smallskip
\noindent
It remains to prove \textup{(d)}. Let
\[
    f(z):=\sum_{n\in E}c_nh_{\sigma,n}(z).
\]
Since the sum is finite and each \(h_{\sigma,n}\) has exponential type
at most \(\sigma\), \(f\) has exponential type at most \(\sigma\).

\smallskip
\noindent
Suppose first that \(0<p<1\). By the \(p\)-triangle inequality,
\[
    \|f\|_{L^p}^p
    \leq
    \sum_{n\in E}|c_n|^p
\|h_{\sigma,n}\|_{L^p}^p 
    =
\sigma^{-1}\|H\|_{L^p}^p
    \sum_{n\in E}|c_n|^p.
\]
Hence
\[
    \|f\|_{L^p}
    \lesssim_{p}
    \sigma^{-\frac{1}{p}}
    \|(c_n)_{n\in E}\|_{\ell^p}.
\]
Now let \(1\leq p\leq\infty\). For \(j\in\Z\), set
\[
    A_j
    :=
    \sup_{|t|\leq1/2}|H(j+t)|.
\]
Since
\(
    |H(x)|\lesssim_p (1+|x|)^{-N_p}
\)
and \(N_p>1\), we have
\[
    A:=(A_j)_{j\in\Z}\in\ell^1,
    \qquad
    \|A\|_{\ell^1}\lesssim_p 1.
\]
If \(k\in\Z\) and \(x\in I_{\sigma,k}\), then
\(
    t:=\sigma x-k-\frac12
    \in\left[-\frac12,\frac12\right),
\)
and hence
\[
    |h_{\sigma,n}(x)|
    =
    |H(k-n+t)|
    \leq
    A_{k-n}.
\]
Therefore
\[
    |f(x)|
    \leq
    \sum_{n\in\Z}|c_n|A_{k-n},
    \qquad x\in I_{\sigma,k},
\]
where the sequence \((c_n)_{n\in E}\) is extended by zero outside
\(E\).

\smallskip
\noindent
If \(1\leq p<\infty\), then
\[
    \|f\|_{L^p}^p
    \leq
    \frac1\sigma
    \sum_{k\in\Z}
    \left(
        \sum_{n\in\Z}|c_n|A_{k-n}
    \right)^p =
    \frac1\sigma
    \|A*|c|\|_{\ell^p}^p.
\]
Now Young's convolution inequality gives
\[
    \|A*|c|\|_{\ell^p}
    \leq
    \|A\|_{\ell^1}\|c\|_{\ell^p},
\]
and therefore
\[
    \|f\|_{L^p}
    \lesssim_p
    \sigma^{-\frac{1}{p}}\|c\|_{\ell^p}.
\]
For \(p=\infty\), the same pointwise estimate yields
\[
    \|f\|_{L^\infty}
    \leq
    \|A\|_{\ell^1}\|c\|_{\ell^\infty}
    \lesssim_p
    \|c\|_{\ell^\infty},
\]
which proves \textup{(d)}.
\end{proof}

\noindent
The localized packets will be used below for reverse and lower estimates.
For the upper sampling bound, however, separation alone is sufficient.

\begin{lemma}[Separated sampling]
\label{lem:separated-sampling}
Let \(0<\sigma,p<\infty\), and let
\(\Lambda=(\lambda_k)_{k\in \Z}\subset\C\) satisfy
\[
    \inf_{j\neq k}|\lambda_j-\lambda_k|
    \geq\delta>0
    \qquad\text{and}\qquad
    \sup_{k \in \Z}|\im\lambda_k|
    \leq M<\infty.
\]
Then
\[
    \|(f(\lambda_k))_{k \in \Z}\|_{\ell^p}
    \lesssim
    \sigma^{\frac1p}\|f\|_{L^p},
    \qquad f\in B_\sigma^p,
\]
where the implicit constant depends only on
\(p\), \(\sigma\delta\), and \(\sigma M\).
\end{lemma}

\begin{proof}
The disks
\[
    D_k
    :=
    \left\{
        z\in\C:
        |z-\lambda_k|<\frac{\delta}{2}
    \right\}
\]
are pairwise disjoint and contained in the strip
\[
    \left\{
        z\in\C:
        |\im z|<M+\frac{\delta}{2}
    \right\}.
\]
Since \(|f|^p\) is subharmonic,
\[
    |f(\lambda_k)|^p
    \leq
    \frac{4}{\pi\delta^2}
    \int_{D_k}|f(z)|^p\,dA(z).
\]
Summing over \(k\) and using
Lemma~\ref{lem:vertical-translation}, we obtain
\[
\sum_{k \in \Z}|f(\lambda_k)|^p
    \leq
    \frac{4}{\pi\delta^2}
    \int_{-(M+\delta/2)}^{M+\delta/2}
    \|f(\cdot+iy)\|_{L^p}^p\,dy \leq
    \frac{4}{\pi\delta^2}
    \left(
        \int_{-(M+\delta/2)}^{M+\delta/2}
        e^{p\sigma|y|}\,dy
    \right)
    \|f\|_{L^p}^p.
\]
With the change of variables \(t=\sigma y\),
\[
    \sum_{k \in \Z}|f(\lambda_k)|^p
    \lesssim_{p,\sigma\delta,\sigma M}
    \sigma\|f\|_{L^p}^p.
\]
Taking \(p\)-th roots proves the assertion.
\end{proof}

\begin{remark}
Lemma~\ref{lem:separated-sampling} is a Plancherel--P\'olya upper
sampling estimate. Such estimates are classical; see, for instance,
\cite[Theorem~6.7.15]{Boas1954} and
\cite[Section~20.1]{Levin1996}. Sampling and interpolation in the
quasi-Banach Paley--Wiener spaces \(0<p\leq1\) have also been studied;
see \cite{Flornes1998}. We include the short proof above because it
applies directly to separated sequences in a bounded horizontal strip
and makes the dependence on the dimensionless parameters
\(\sigma\delta\) and \(\sigma M\) explicit. For the broader sampling
and interpolation theory in Paley--Wiener spaces, see, for instance,
Landau~\cite{Landau1967}, Lyubarskii and Seip~\cite{LyubarskiiSeip1997},
and Ortega-Cerd\`a and Seip~\cite{OrtegaCerdaSeip2002}.
\end{remark}

\begin{corollary}
\label{cor:local-supremum-sampling}
Let \(0<\sigma, p<\infty\). Then
\[
    \sum_{n\in\Z}
    \left(
        \sup_{x\in I_{\sigma,n}}|f(x)|
    \right)^p
    \lesssim
    \sigma\|f\|_{L^p}^p,
    \qquad f\in B_\sigma^p,
\]
where the implicit constant depends only on \(p\).
\end{corollary}

\begin{proof}
For each \(n\in\Z\), choose
\(x_n\in \left[\frac{n}{\sigma},\frac{n+1}{\sigma}\right]\) such that
\[
    |f(x_n)|
    =
    \sup_{x\in I_{\sigma,n}}|f(x)|.
\]
The sequences \((x_{2n})_{n\in\Z}\) and \((x_{2n+1})_{n\in\Z}\) are
\(\frac{1}{\sigma}\)-separated and lie on the real axis. Hence
Lemma~\ref{lem:separated-sampling}, applied to each sequence, gives
\[
    \sum_{n\in\Z}
    \left(
        \sup_{x\in I_{\sigma,n}}|f(x)|
    \right)^p
    =
    \sum_{n\in\Z}|f(x_{2n})|^p
    +
    \sum_{n\in\Z}|f(x_{2n+1})|^p \lesssim_p
    \sigma\|f\|_{L^p}^p.
\]
\end{proof}

\noindent
The preceding estimate discretizes local suprema. The next lemma gives
a complementary estimate in terms of local \(L^q\)-masses.

\begin{lemma}[Local \(L^q\)-to-\(L^p\) estimate]
\label{lem:local-q-to-p}
Let \(\tau,\rho>0\) and \(0<q<p<\infty\). If \(f\) is an entire
function of exponential type at most \(\tau\), then
\[
    \|f\|_{L^p}^p
    \lesssim
    (\tau+\rho)^{\frac pq-1}
    \sum_{n\in\Z}
    \left(
        \int_{I_{\rho,n}}|f(x)|^q\,dx
    \right)^{\frac pq},
\]
where the implicit constant depends only on \(p\) and \(q\).
\end{lemma}

\begin{proof}
Set
\[
    s:=\frac{p}{q}>1,
    \qquad
    A_n:=\int_{I_{\rho,n}}|f(x)|^q\,dx.
\]
There is nothing to prove if \((A_n)_{n\in\Z}\notin\ell^s\), so assume
\(A:=(A_n)_{n\in\Z}\in\ell^s\).

\smallskip
\noindent
Let \(h_{\rho,n}\) be the localized packets of
Lemma~\ref{lem:localized-packets} corresponding to the exponent \(q\). If \(x\in I_{\rho,k}\), then the definition of \(h_{\rho,n}\) and
the elementary estimate \eqref{eq:H}
give
\[
    |h_{\rho,n}(x)|
    \lesssim_q
    (1+|k-n|)^{-N_q}.
\] 
Hence, with
\[
    \omega_j:=(1+|j|)^{-qN_q},
    \qquad j\in\Z,
\]
we obtain
\[
\|fh_{\rho,n}\|_{L^q}^q
    =
    \sum_{k\in\Z}
    \int_{I_{\rho,k}}
    |f(x)|^q|h_{\rho,n}(x)|^q\,dx \lesssim_q
    \sum_{k\in\Z}
    \omega_{n-k}A_k
    =
    (\omega*A)_n.
\]
Since \(qN_q>1\), we have \(\omega\in\ell^1\). Moreover,
\(A\in\ell^s\subset\ell^\infty\), and hence
\[
    (\omega*A)_n
    \leq
    \|\omega\|_{\ell^1}\|A\|_{\ell^\infty}
    <\infty.
\]
Thus \(fh_{\rho,n}\in L^q\).
Since \(fh_{\rho,n}\) has exponential type at most \(\tau+\rho\), we have $fh_{\rho,n}\in B_{\tau+\rho}^q$. Lemma \ref{lem:nikolskii} therefore yields
\[
    \|fh_{\rho,n}\|_{L^\infty}^q
    \lesssim_q
    (\tau+\rho)\|fh_{\rho,n}\|_{L^q}^q
    \lesssim_q
    (\tau+\rho)(\omega*A)_n.
\]
By Lemma~\ref{lem:localized-packets},
\(h_{\rho,n}\geq\eta>0\) on \(I_{\rho,n}\), where \(\eta\) depends
only on \(q\). Hence
\[
    \sup_{x\in I_{\rho,n}}|f(x)|^q
    \lesssim_q
    (\tau+\rho)(\omega*A)_n.
\]
Consequently,
\[
\int_{I_{\rho,n}}|f(x)|^p\,dx
    \leq
    \left(
        \sup_{x\in I_{\rho,n}}|f(x)|^q
    \right)^{s-1}A_n  \lesssim_{p,q}
    (\tau+\rho)^{s-1}
    (\omega*A)_n^{s-1}A_n.
\]
Since \(\omega_0=1\), we have \(A_n\leq(\omega*A)_n\), and therefore
\[
    \int_{I_{\rho,n}}|f(x)|^p\,dx
    \lesssim_{p,q}
    (\tau+\rho)^{s-1}
    (\omega*A)_n^s.
\]
Summing over \(n\) and applying Young's inequality on \(\ell^s\),
\[
\begin{aligned}
    \|f\|_{L^p}^p
    &\lesssim_{p,q}
    (\tau+\rho)^{s-1}
    \|\omega*A\|_{\ell^s}^s \leq
    (\tau+\rho)^{s-1}
    \|\omega\|_{\ell^1}^s
    \|A\|_{\ell^s}^s \\
    &\lesssim_{p,q}
    (\tau+\rho)^{\frac pq-1}
    \sum_{n\in\Z}
    \left(
        \int_{I_{\rho,n}}|f(x)|^q\,dx
    \right)^{\frac pq}.
\end{aligned}
\]
\end{proof}


\subsection{Carleson measures}
\label{subsec:carleson}

The Carleson measure problem for Paley--Wiener spaces was considered
in connection with composition operators in
\cite[Theorems~2.9 and~2.10]{ChaconChaconGimenez2007}, while a
substantially more general full-range theory for Bernstein spaces is
developed in \cite{CalziPeloso2023}. For the applications below, we use
a discrete formulation adapted to the natural scale \(\frac{1}{\sigma}\).
Besides giving direct proofs from the localization and
sampling estimates above, this formulation yields quantitative norm
estimates with the dependence on the bandwidth \(\sigma\) made explicit.
These estimates will be important in the operator-theoretic arguments
later in the paper.

\smallskip
\noindent
Fix \(\sigma>0\), and let \(\mu\) be a positive Borel measure on \(\R\). For
\(0<p\leq\infty\) and \(0<q<\infty\), consider the natural embedding
\[
    J_\mu:B_\sigma^p\longrightarrow L^q(\mu),
    \qquad
    J_\mu f=f.
\]
We call \(\mu\) a \textit{\((p,q)\)-Carleson measure} for
\(B_\sigma^p\) if \(J_\mu\) is bounded, and a
\textit{\((p,q)\)-vanishing Carleson measure} for \(B_\sigma^p\)
if \(J_\mu\) is compact. We write
\[
    \|\mu\|_{\mathcal C_{p,q}(\sigma)}
    :=
    \|J_\mu\|_{B_\sigma^p\to L^q(\mu)}
    =
    \sup_{\|f\|_{L^p}\leq1}
    \left(
        \int_{\R}|f(x)|^q\,d\mu(x)
    \right)^{\frac{1}{q}},
\]
with the convention that
\(\|\mu\|_{\mathcal C_{p,q}(\sigma)}=\infty\) if the embedding is
not bounded. We put
\[
    \mu_{\sigma,n}:=\mu(I_{\sigma,n}),
    \quad n\in\Z.
\]

\begin{theorem}\label{thm-carleson-p-le-q}
Let \(\sigma>0\), \(0<p\leq q<\infty\), and let \(\mu\) be a positive
Borel measure on \(\R\). Then the following assertions hold.
\begin{enumerate}[label=\textnormal{(\alph*)}]
\item \(\mu\) is a \((p,q)\)-Carleson measure for \(B_\sigma^p\) if
and only if
\( (\mu_{\sigma,n})_{n\in \Z} \in \ell^{\infty}\).
Moreover,
\[
    \|\mu\|_{\mathcal C_{p,q}(\sigma)}
    \simeq
    \sigma^{\frac{1}{p}}
    \|(\mu_{\sigma,n})_{n \in \Z}\|_{\ell^{\infty}}^{\frac{1}{q}},
\]
where the implicit constants depend only on \(p\) and \(q\).

\item \(\mu\) is a  \((p,q)\)-vanishing Carleson measure for
\(B_\sigma^p\) if and only if
\(\mu_{\sigma,n}\longrightarrow0\) as \(|n|\to\infty\).
\end{enumerate}
\end{theorem}

\begin{proof}
Suppose first that \((\mu_{\sigma,n})_{n\in\Z} \in \ell^{\infty}\). For \(f\in B_\sigma^p\), since \(p\leq q\), Corollary~\ref{cor:local-supremum-sampling} gives
\[
\begin{aligned}
    \|J_\mu f\|_{L^q(\mu)}^q
    & \leq
    \|(\mu_{\sigma,n})_{n\in\Z}\|_{\ell^\infty}\sum_{n\in\Z}\left( \sup_{x\in I_{\sigma,n}}|f(x)| \right)^q \\
   & \leq
    \|(\mu_{\sigma,n})_{n\in\Z}\|_{\ell^\infty}\left(\sum_{n\in\Z}\left( \sup_{x\in I_{\sigma,n}}|f(x)| \right)^p\right)^{\frac{q}{p}} \lesssim_p
    \sigma^{\frac{q}{p}} \|(\mu_{\sigma,n})_{n\in\Z}\|_{\ell^\infty} \|f\|_{L^p}^q.
    \end{aligned}
\]
Hence
\[
    \|\mu\|_{\mathcal C_{p,q}(\sigma)}
    \lesssim_{p,q}
    \sigma^{\frac{1}{p}}\|(\mu_{\sigma,n})_{n\in\Z}\|_{\ell^\infty}^{\frac{1}{q}}.
\]
Conversely, let
\(  \|\mu\|_{\mathcal C_{p,q}(\sigma)}<\infty.\)
Choose the localized packets \(h_{\sigma,n}\) from
Lemma~\ref{lem:localized-packets} corresponding to the exponent \(p\). Since \(h_{\sigma,n}\geq\eta>0\) on \(I_{\sigma,n}\),
\[
    \eta^q\mu_{\sigma,n}
    \leq
    \|J_\mu h_{\sigma,n}\|_{L^q(\mu)}^q
    \leq
    \|\mu\|_{\mathcal C_{p,q}(\sigma)}^q\|h_{\sigma,n}\|_{L^p}^q
    \lesssim_p
    \|\mu\|_{\mathcal C_{p,q}(\sigma)}^q\sigma^{-\frac{q}{p}}.
\]
Taking the supremum over \(n\) yields
\[
    \sigma^{\frac{1}{p}}
    \|(\mu_{\sigma,n})_{n\in\Z}\|_{\ell^\infty}^{\frac{1}{q}}
    \lesssim_{p,q}
    \|\mu\|_{\mathcal C_{p,q}(\sigma)},
\]
which proves (a).

\smallskip
\noindent
We turn to compactness. First observe that \(J_\mu\) has the uniqueness
property in Lemma~\ref{lem-compactness-criterion}. Indeed, if
\(f_j\to0\) locally uniformly and \(J_\mu f_{j_k}\to g\) in
\(L^q(\mu)\), then a further subsequence converges to \(g\)
\(\mu\)-almost everywhere. Since \(f_{j_k}(x)\to0\) for every
\(x\in\R\), we have \(g=0\) \(\mu\)-almost everywhere.

\smallskip
\noindent
Suppose that \(J_\mu\) is compact. By
Lemma~\ref{lem:localized-packets}, the sequence
\((h_{\sigma,n})_{n \in \Z}\) is bounded in \(B_\sigma^p\) and
\(
    h_{\sigma,n}\longrightarrow0
\)
uniformly on compact subsets of \(\C\) as \(|n|\to\infty\).
Lemma~\ref{lem-compactness-criterion} therefore gives
\[
    \|J_\mu h_{\sigma,n}\|_{L^q(\mu)}\longrightarrow0.
\]
Since \(h_{\sigma,n}\geq\eta\) on \(I_{\sigma,n}\),
\(
    \eta^q\mu_{\sigma,n}
    \leq
    \|J_\mu h_{\sigma,n}\|_{L^q(\mu)}^q,
\)
and hence $\mu_{\sigma,n}\longrightarrow0$ as $|n|\to\infty$.

\smallskip
\noindent
Conversely, suppose that \(\mu_{\sigma,n}\to0\). Then
\((\mu_{\sigma,n})_n\in\ell^\infty\), so \(J_\mu\) is bounded.
Let \((f_j)_{j \in \N}\) be bounded in \(B_\sigma^p\) and converge to \(0\)
uniformly on compact subsets of \(\C\).
Since \(p \leq q\), by Corollary~\ref{cor:local-supremum-sampling}, for \(N\in \N\),
\[
\begin{aligned}
    \|J_\mu f_j\|_{L^q(\mu)}^q
    &\leq
    \sum_{|n|\leq N}\mu_{\sigma,n}\left(\sup_{x\in I_{\sigma,n}}|f_j(x)| \right)^q
    +
    \left(\sup_{|n|>N}\mu_{\sigma,n}\right)
    \sum_{|n|>N}\left(\sup_{x\in I_{\sigma,n}}|f_j(x)| \right)^q \\
    &\leq
    \sum_{|n|\leq N}\mu_{\sigma,n}\left(\sup_{x\in I_{\sigma,n}}|f_j(x)| \right)^q
    +
    \left(\sup_{|n|>N}\mu_{\sigma,n}\right)
    \left(\sum_{n\in\Z}\left(\sup_{x\in I_{\sigma,n}}|f_j(x)| \right)^p\right)^{\frac{q}{p}} \\
    &\lesssim_p
    \sum_{|n|\leq N}\mu_{\sigma,n}\left(\sup_{x\in I_{\sigma,n}}|f_j(x)| \right)^q
    +
    \sigma^{\frac{q}{p}} \left(\sup_{j \in \N}\|f_j\|_{L^p}\right)^q
    \left(\sup_{|n|>N}\mu_{\sigma,n}\right).
\end{aligned}
\]
For fixed \(N\), the first term tends to \(0\) as \(j\to\infty\);
the second tends to \(0\) as \(N\to\infty\), uniformly in \(j\).
Thus
\[
    \|J_\mu f_j\|_{L^q(\mu)}\longrightarrow0.
\]
Lemma~\ref{lem-compactness-criterion} now implies that \(J_\mu\) is
compact.
\end{proof}

\begin{theorem}\label{thm-carleson-q-lt-p}
Let \(\sigma>0\), \(0<q<p\leq\infty\), and let \(\mu\) be a positive
Borel measure on \(\R\). Then the following are equivalent:
\begin{enumerate}[label=\textnormal{(\alph*)}]
\item \(\mu\) is a \((p,q)\)-Carleson measure for \(B_\sigma^p\);
\item \(\mu\) is a  \((p,q)\)-vanishing Carleson measure for
\(B_\sigma^p\);
\item
\(
    (\mu_{\sigma,n})_{n\in\Z}
    \in
    \ell^{\left(\frac{p}{q}\right)'},
\)
\end{enumerate}
where \(\left(\frac{p}{q}\right)'\) denotes the H\"older conjugate of \(\frac{p}{q}\);
in particular, \(\left(\frac{\infty}{q}\right)'=1\). Moreover,
\[
    \|\mu\|_{\mathcal C_{p,q}(\sigma)}
    \simeq
    \sigma^{\frac{1}{p}}
    \|(\mu_{\sigma,n})_{n\in\Z}\|_{\ell^{\left(\frac{p}{q}\right)'}}^{\frac{1}{q}},
\]
where the implicit constants depend only on \(p\) and \(q\). For
\(p=\infty\), the last estimate is in fact the identity
\[
    \|\mu\|_{\mathcal C_{\infty,q}(\sigma)}
    =
    \mu(\R)^{\frac{1}{q}}.
\]
\end{theorem}

\begin{proof}
Suppose first that \(p<\infty\). Then
\( \left( \frac{p}{q}\right)'= \frac{p}{p-q}\).

\smallskip
\noindent
Assume that
\((\mu_{\sigma,n})_{n \in \Z}\in\ell^{\frac{p}{p-q}}\). For \(f\in B_\sigma^p\),
H\"older's inequality and
Corollary~\ref{cor:local-supremum-sampling} give
\[
\begin{aligned}
    \|J_\mu f\|_{L^q(\mu)}^q
    &\leq
    \sum_{n\in\Z}\mu_{\sigma,n}\left(\sup_{x\in I_{\sigma,n}}|f(x)|\right)^q \\
    &\leq
    \|(\mu_{\sigma,n})\|_{\ell^\frac{p}{p-q}}
    \left(\sum_{n\in\Z}\left(\sup_{x\in I_{\sigma,n}}|f(x)|\right)^p\right)^{\frac{q}{p}} \\
    &\lesssim_p
    \sigma^{\frac{q}{p}}
    \|(\mu_{\sigma,n})\|_{\ell^\frac{p}{p-q}}
    \|f\|_{L^p}^q.
\end{aligned}
\]
Hence
\begin{equation}\label{eq-carleson-q-lt-p-upper}
    \|\mu\|_{\mathcal C_{p,q}(\sigma)}
    \lesssim_{p,q}
    \sigma^{\frac{1}{p}}
    \|(\mu_{\sigma,n})\|_{\ell^\frac{p}{p-q}}^{\frac{1}{q}}.
\end{equation}
Conversely, let
\(
    \|\mu\|_{\mathcal C_{p,q}(\sigma)}<\infty.
\)
Choose the localized packets \(h_{\sigma,n}\) from
Lemma~\ref{lem:localized-packets} corresponding to the exponent \(p\). For a finite set \(E\subset\Z\)
and numbers \(c_n\geq0\), put
\[
    f:=\sum_{n\in E}c_nh_{\sigma,n}.
\]
Since the packets are nonnegative on \(\R\) and
\(h_{\sigma,n}\geq\eta\) on \(I_{\sigma,n}\),
\[
    \eta^q\sum_{n\in E}c_n^q\mu_{\sigma,n}
    \leq
    \|J_\mu f\|_{L^q(\mu)}^q.
\]
By the synthesis estimate in Lemma~\ref{lem:localized-packets},
\begin{equation}\label{eq-carleson-packet-testing}
    \sum_{n\in E}c_n^q\mu_{\sigma,n}
    \lesssim_{p,q}
    \left( \|\mu\|_{\mathcal C_{p,q}(\sigma)}\right)^q\sigma^{-\frac{q}{p}}
    \left(\sum_{n\in E}c_n^p\right)^{q/p}.
\end{equation}
Taking \(E=\{n\}\) shows that every \(\mu_{\sigma,n}\) is finite. Setting $c_n=(\mu_{\sigma,n})^{\frac{1}{p-q}}$ in \eqref{eq-carleson-packet-testing}, we find
\[
    \left( \sum_{n\in E}\mu_{\sigma,n}^{\frac{p}{p-q}}\right)^{\frac{p-q}{p}}
    \lesssim_{p,q}
    \left( \|\mu\|_{\mathcal C_{p,q}(\sigma)}\right)^q\sigma^{-\frac{q}{p}} .
\]
Since \(E\) is arbitrary,
\(
    (\mu_{\sigma,n})_{n\in\Z}\in\ell^\frac{p}{p-q}
\)
and
\begin{equation*}\label{eq-carleson-q-lt-p-lower}
    \sigma^{\frac{1}{p}}
    \|(\mu_{\sigma,n})\|_{\ell^\frac{p}{p-q}}^{\frac{1}{q}}
    \lesssim_{p,q}
    \|\mu\|_{\mathcal C_{p,q}(\sigma)}.
\end{equation*}
Together with \eqref{eq-carleson-q-lt-p-upper}, this proves the
boundedness characterization and the norm estimate.

\smallskip
\noindent
It remains to prove compactness. As in the proof of
Theorem~\ref{thm-carleson-p-le-q}, \(J_\mu\) satisfies the uniqueness
hypothesis of Lemma~\ref{lem-compactness-criterion}. Let therefore
\((f_j)_{j \in \N}\) be bounded in \(B_\sigma^p\) and converge to \(0\) uniformly
on compact subsets of \(\C\). 
For \(N\in \N\), H\"older's inequality and
Corollary~\ref{cor:local-supremum-sampling} give
\[
\begin{aligned}
    \|J_\mu f_j\|_{L^q(\mu)}^q
    &\leq
    \sum_{|n|\leq N}\mu_{\sigma,n}\left( \sup_{x\in I_{\sigma,n}}|f_j(x)|\right)^q 
    +
    \|(\mu_{\sigma,n})_{|n|>N}\|_{\ell^\frac{p}{p-q}}
    \left(\sum_{n\in\Z}\left( \sup_{x\in I_{\sigma,n}}|f_j(x)|\right)^p\right)^{\frac{q}{p}} \\
    &\lesssim_p
    \sum_{|n|\leq N}\mu_{\sigma,n}\left( \sup_{x\in I_{\sigma,n}}|f_j(x)|\right)^q
    +
    \sigma^{\frac{q}{p}}\left( \sup_{j \in \N}\|f_j\|_{L^p}\right)^q
    \|(\mu_{\sigma,n})_{|n|>N}\|_{\ell^\frac{p}{p-q}}.
\end{aligned}
\]
For fixed \(N\), the first term tends to \(0\) as \(j\to\infty\), while
the second tends to \(0\) as \(N\to\infty\), uniformly in \(j\).
Hence
\[
    \|J_\mu f_j\|_{L^q(\mu)}\longrightarrow0.
\]
Lemma~\ref{lem-compactness-criterion} shows that \(J_\mu\) is compact.
Thus (c) implies (b), while (b) trivially implies (a).

\smallskip
\noindent
Finally, let \(p=\infty\). Since the cells \(I_{\sigma,n}\) partition
\(\R\),
\[
    \|(\mu_{\sigma,n})_{n \in \Z}\|_{\ell^1}
    =
    \sum_{n\in\Z}\mu_{\sigma,n}
    =
    \mu(\R).
\]
If \(\mu(\R)<\infty\), then
\[
    \|J_\mu f\|_{L^q(\mu)}
    \leq
    \mu(\R)^{\frac{1}{q}}\|f\|_{L^\infty},
\]
while testing on the constant function \(1\) gives the reverse
inequality. Therefore
\[
    \|\mu\|_{\mathcal C_{\infty,q}(\sigma)}
    =
    \mu(\R)^{\frac{1}{q}}.
\]
Moreover, if \((f_j)_{j \in \N}\) is bounded in \(B_\sigma^\infty\) and
\(f_j\to0\) locally uniformly, then \(f_j(x)\to0\) for every
\(x\in\R\), and the dominated convergence theorem gives
\[
    \|J_\mu f_j\|_{L^q(\mu)}\longrightarrow0.
\]
Hence \(J_\mu\) is compact by
Lemma~\ref{lem-compactness-criterion}. This completes the proof.
\end{proof}

\noindent
To relate the preceding discrete criteria to the usual continuous
local-mass formulation, for \(r>0\) set
\[
    \widetilde{\mu}_{\sigma,r}(x)
    :=
    \mu\left(
        x-\frac r\sigma,
        x+\frac r\sigma
    \right),
    \qquad x\in\R.
\]

\begin{lemma}\label{lem-local-mass-equivalences}
Let \(\sigma, r > 0\) and \(\mu\) be a positive Borel measure on \(\mathbb R\). Then the following
assertions hold.

\begin{enumerate}[label=\textnormal{(\alph*)}]
\item
\( \displaystyle
\sup_{n\in\Z}\mu_{\sigma,n} \simeq
   \sup_{x\in\R}\widetilde{\mu}_{\sigma,r}(x)\).

\item
\(   \mu_{\sigma,n}\longrightarrow0\) as \(|n|\to\infty\) 
if and only if
\(  \widetilde{\mu}_{\sigma,r}(x)\longrightarrow0\) as \(|x|\to\infty\).

\item If \(1\leq s<\infty\), then
\[
    \|\widetilde{\mu}_{\sigma,r}\|_{L^s}
    \simeq
    \sigma^{-\frac{1}{s}}
    \|(\mu_{\sigma,n})_{n\in\Z}\|_{\ell^s}.
\]
\end{enumerate}
The implicit constants in \textnormal{(a)} and \textnormal{(c)}
depend only on \(r\) and on \(r,s\), respectively.
\end{lemma}

\begin{proof}
If \(x\in I_{\sigma,n}\), then
\[
    \left(
        x-\frac r\sigma,
        x+\frac r\sigma
    \right)
    \subset
    \bigcup_{|k-n|\leq \lceil r\rceil+1}I_{\sigma,k},
\]
so that
\begin{equation}\label{eq-local-mass-upper}
    \widetilde{\mu}_{\sigma,r}(x)
    \leq
    \sum_{|k-n|\leq \lceil r\rceil+1}\mu_{\sigma,k},
    \qquad x\in I_{\sigma,n}.
\end{equation}
Conversely, with
\(
    I_{\sigma,n}^{(r)}
    :=
    \left(
        \frac{n-r}{\sigma},
        \frac{n+1+r}{\sigma}
    \right),
\)
Tonelli's theorem gives
\begin{equation}\label{eq-local-mass-lower}
    \int_{I_{\sigma,n}^{(r)}}
        \widetilde{\mu}_{\sigma,r}(x)\,dx
    \geq
    \frac{2r}{\sigma}\mu_{\sigma,n}.
\end{equation}
Since
\(
    |I_{\sigma,n}^{(r)}|
    =
    \frac{1+2r}{\sigma},
\)
we obtain
\[
    \mu_{\sigma,n}
    \leq
    \frac{1+2r}{2r}
    \sup_{x\in I_{\sigma,n}^{(r)}}
    \widetilde{\mu}_{\sigma,r}(x).
\]
Together with \eqref{eq-local-mass-upper}, this proves (a) and (b).

\smallskip
\noindent
Now let \(1\leq s<\infty\). From \eqref{eq-local-mass-upper},
\[
    \|\widetilde{\mu}_{\sigma,r}\|_{L^s}^s
    \lesssim_{r,s}
    \frac1\sigma
    \sum_{n\in\Z}\mu_{\sigma,n}^s.
\]
On the other hand, H\"older's inequality and
\eqref{eq-local-mass-lower} yield
\[
    \frac1\sigma\mu_{\sigma,n}^s
    \lesssim_{r,s}
    \int_{I_{\sigma,n}^{(r)}}
        \widetilde{\mu}_{\sigma,r}(x)^s\,dx.
\]
Since the intervals \(I_{\sigma,n}^{(r)}\) have uniformly bounded
overlap, summing over \(n\) gives
\[
    \frac1\sigma
    \sum_{n\in\Z}\mu_{\sigma,n}^s
    \lesssim_{r,s}
    \|\widetilde{\mu}_{\sigma,r}\|_{L^s}^s.
\]
\end{proof}

\begin{remark}\label{rem-carleson-continuous-form}
Lemma~\ref{lem-local-mass-equivalences} converts the preceding discrete
criteria into the usual continuous local-mass formulation. For every
\(r>0\), if \(0<p\leq q<\infty\), then \(\mu\) is a
\((p,q)\)-Carleson measure for \(B_\sigma^p\) if and only if
\[
    \widetilde{\mu}_{\sigma,r}\in L^\infty(\R),
\qquad \text{and} \qquad
    \|\mu\|_{\mathcal C_{p,q}(\sigma)}
    \simeq_{p,q,r}
    \sigma^{\frac{1}{p}}
    \|\widetilde{\mu}_{\sigma,r}\|_{L^\infty}^{\frac{1}{q}}.
\]
Moreover,
\(\mu\) is a  \((p,q)\)-vanishing Carleson measure if and only if
\[
    \widetilde{\mu}_{\sigma,r}(x)\longrightarrow0
    \quad \text{as} \quad |x|\to\infty.
\]
If \(0<q<p\leq\infty\), then \(\mu\) is a
\((p,q)\)-Carleson measure if and only if it is a 
\((p,q)\)-vanishing Carleson measure, and these conditions are equivalent to
\[
    \widetilde{\mu}_{\sigma,r}
    \in L^{\left(\frac{p}{q}\right)'}(\R)
\qquad \text{and} \qquad
    \|\mu\|_{\mathcal C_{p,q}(\sigma)}
    \simeq_{p,q,r}
    \sigma^{\frac{1}{q}}
    \|\widetilde{\mu}_{\sigma,r}\|_{L^{\left(\frac{p}{q}\right)'}}^{\frac{1}{q}}.
\]
\end{remark}

\section{Weighted composition operators}
\label{sec-wco}

We now apply the analytic tools developed in Section~\ref{sec:analytic-tools} to weighted
composition operators between Paley--Wiener spaces of possibly different
bandwidths. Let \(\sigma,\tau>0\) and \(0<p,q\leq\infty\). For entire
functions \(\psi\) and \(\varphi\), define
\[
    W_{\psi,\varphi}f
    :=
    \psi\cdot(f\circ\varphi).
\]
We study bounded and compact operators
\[
    W_{\psi,\varphi}:
    B_\sigma^p\longrightarrow B_\tau^q,
\]
together with their essential norms and two basic specializations:
composition and multiplication operators. The case \(\psi\equiv0\) is
trivial and will be excluded unless otherwise stated.

\smallskip
\noindent
The first reduction is by scaling. Under the isometries from
Subsection~\ref{subsec:scaling}, set
\begin{equation}\label{eq-symbols-scaling}
    \widetilde\psi(z)
    :=
    \left(\frac{\pi}{\tau}\right)^{\frac{1}{q}}
    \left(\frac{\sigma}{\pi}\right)^{\frac{1}{p}}
    \psi\left(\frac{\pi}{\tau}z\right)
\qquad \text{and} \qquad
    \widetilde\varphi(z)
    :=
    \frac{\sigma}{\pi}
    \varphi\left(\frac{\pi}{\tau}z\right).
\end{equation}
Then
\begin{equation}\label{eq-wco-scaling}
    \mathcal U_\tau
    W_{\psi,\varphi}
    \mathcal U_\sigma^{-1}
    =
    W_{\widetilde\psi,\widetilde\varphi}
    :
    B_\pi^p\longrightarrow B_\pi^q.
\end{equation}
Thus boundedness, compactness, and operator norms may be studied in the
normalized setting
\(
B_\pi^p\longrightarrow B_\pi^q.
\)
Whenever this simplifies the argument, we shall work first in this
setting and then return to arbitrary \(\sigma,\tau>0\) by
\eqref{eq-symbols-scaling}-\eqref{eq-wco-scaling}.

\subsection{Structural reduction}
\label{subsec-wco-structural}

We begin with the structural consequences of boundedness. These arguments
do not depend on the relative size of \(p\) and \(q\), and will therefore
be used in both exponent ranges considered below.

\smallskip
\noindent
The zero-bandwidth space introduced above satisfies
\[
    B_0^\infty=\mathbb C.
\]
Indeed, if
\(f\) is an entire function of exponential type zero and
\(f|_{\R}\in L^\infty(\R)\), then \(f\in B_\varepsilon^\infty\) for
every \(\varepsilon>0\). By Lemma~\ref{lem:vertical-translation},
\[
    |f(x+iy)|
    \leq
    e^{\varepsilon|y|}\|f\|_{L^\infty},
    \qquad x,y\in\R.
\]
Letting \(\varepsilon\to0^+\), we obtain
\[
    |f(z)|\leq\|f\|_{L^\infty},
    \qquad z\in\C,
\]
so \(f\) is constant by Liouville's theorem.

\medskip
\noindent
We first dispose of constant symbols.

\begin{proposition}[Constant symbols]
\label{prop-wco-constant-symbol}
Let \(\sigma,\tau>0\), \(0<p,q\leq\infty\), and let
\(\varphi\equiv b\in\C\). Let \(\psi\not\equiv0\) be entire. Then the
following conditions are equivalent:
\begin{enumerate}
\item[\textnormal{(a)}]
\(
    W_{\psi,\varphi}:
    B_\sigma^p\longrightarrow B_\tau^q
\)
is bounded;

\item[\textnormal{(b)}]
\(
    W_{\psi,\varphi}:
    B_\sigma^p\longrightarrow B_\tau^q
\)
is compact;

\item[\textnormal{(c)}]
\(
    \psi\in B_\tau^q.
\)
\end{enumerate}
Moreover,
\[
    \|W_{\psi,\varphi}\|_{B_\sigma^p\to B_\tau^q}
    =
    \|\psi\|_{L^q}
    \|\delta_b\|_{(B_\sigma^p)^*}
    \simeq
    \|\psi\|_{L^q}
    \Phi_{\sigma,p}(|\operatorname{Im}b|),
\]
where the implicit constants depend only on \(p\).
\end{proposition}

\begin{proof}
Since $W_{\psi,\varphi}f=\psi f(b),$ the assertion follows from the point-evaluation estimate in
Theorem~\ref{thm:point-evaluations}. If \(\psi\in B_\tau^q\), then the range of
\(W_{\psi,\varphi}\) is contained in the one-dimensional space
spanned by \(\psi\).
\end{proof}

\noindent
Having disposed of constant symbols, we henceforth assume, unless
otherwise stated, that
\[
    \psi\not\equiv0
    \qquad\text{and}\qquad
    \varphi\ \text{is nonconstant}.
\]
The next two rigidity lemmas will be used to show that boundedness forces
\(\varphi\) to be affine with real slope. We include their proofs for
completeness.

\begin{lemma}[Product rigidity]
\label{lem-product-rigidity}
Let \(f\not\equiv0\) and \(g\) be entire functions. If both \(f\) and
\(fg\) are of finite exponential type, then \(g\) is of finite
exponential type.
\end{lemma}

\begin{proof}
Suppose, to the contrary, that \(g\) is not of finite exponential type. In Levin's terminology, \(g\) belongs to a strictly higher growth
category than \(f\): since \(f\) is of finite exponential type whereas
\(g\) is not, either \(g\) has order greater than one or it has order
one and maximal type. By the classical product theorem
\cite[Chapter~I, Section~9, Theorem~12]{Levin1980}, the product \(fg\)
then belongs to the same growth category as \(g\). This contradicts the
assumption that \(fg\) is of finite exponential type.
\end{proof}

\begin{lemma}[Exponential-decay rigidity]
\label{lem-exponential-decay-rigidity}
Let \(f\) be an entire function of finite exponential type. If there
exist constants \(C,\gamma>0\) such that
\(|f(x)|
    \leq
    Ce^{-\gamma|x|}\) on \( x\in\R\),
then \(f\equiv0\).
\end{lemma}

\begin{proof}
Choose \(\rho>0\) such that \(f\) has exponential type at most
\(\rho\). The assumed decay implies that
\[
    f|_{\R}\in L^1(\R)\cap L^2(\R).
\]
Hence \(f\in B_\rho^2\), and the Paley--Wiener theorem \cite[Theorem~6.8.1]{Boas1954} gives
\[
    \operatorname{supp}\widehat f\subset[-\rho,\rho].
\]
On the other hand, the assumed decay also shows that
the integral
\[
    \widehat f(z)
    =
    \frac1{\sqrt{2\pi}}
    \int_{\R}f(x)e^{-ixz}\,dx
\]
converges absolutely and locally uniformly  for
\(|\operatorname{Im}z|<\gamma\). 
Thus \(\widehat f\) extends holomorphically to the strip
\(\{z\in\C:|\operatorname{Im}z|<\gamma\}.\) Since \(f\in L^1(\R)\), its Fourier transform is continuous on \(\R\).
The support condition therefore implies
\[
    \widehat f(x)=0,
    \qquad |x|>\rho.
\]
By the identity theorem, the holomorphic extension of \(\widehat f\)
vanishes identically on the strip. Hence \(\widehat f=0\) on \(\R\),
and Fourier uniqueness gives \(f=0\) on \(\R\). The identity theorem
then yields \(f\equiv0\).
\end{proof}

\medskip
\noindent
The point-evaluation estimates from Theorem~\ref{thm:point-evaluations} provide the quantitative
quantity that governs the structural reduction. Define
\[
    M_{\sigma,\tau;p,q}(\psi,\varphi)
    :=
    \sup_{z\in\C}
    |\psi(z)|
    \frac{\Phi_{\sigma,p}(|\operatorname{Im}\varphi(z)|)}
         {\Phi_{\tau,q}(|\operatorname{Im}z|)}.
\]
A direct calculation from \eqref{eq-symbols-scaling} and the definition of
\(\Phi_{\sigma,p}\) gives
\[
    M_{\sigma,\tau;p,q}(\psi,\varphi)
    =
    M_{\pi,\pi;p,q}
    \left(\widetilde\psi,\widetilde\varphi\right).
\]
The following proposition combines this pointwise estimate with the
preceding rigidity lemmas.

\begin{proposition}[Structural reduction]
\label{prop-wco-structural}
Let \(\sigma,\tau>0\) and \(0<p,q\leq\infty\). Let \(\psi,\varphi\)
be entire functions such that \(\psi\not\equiv0\) and \(\varphi\) is
nonconstant. Suppose that
\(
    W_{\psi,\varphi}:
    B_\sigma^p\longrightarrow B_\tau^q
\)
is bounded. Then
\[
    M_{\sigma,\tau;p,q}(\psi,\varphi)<\infty
\qquad \text{and} \qquad
    \varphi(z)=az+b,
\]
where
\(a\in\R\) with \(0<\sigma|a|\leq\tau\) and \(b\in\C\).
Moreover,
\[
    M_{\sigma,\tau;p,q}(\psi,\varphi)
    \lesssim
    \|W_{\psi,\varphi}\|_{B_\sigma^p\to B_\tau^q},
\]
where the implicit constant depends only on \(p\) and \(q\).
\end{proposition}

\begin{proof}
We first prove the assertion in the normalized case
\(W_{\psi,\varphi}:B_\pi^p\longrightarrow B_\pi^q.\)

\smallskip
\noindent
\emph{Step 1. Point evaluations.}
For every \(z\in\C\),
\(
    \delta_z\circ W_{\psi,\varphi}
    =
    \psi(z)\delta_{\varphi(z)}.
\)
Hence
\[
    |\psi(z)|
    \|\delta_{\varphi(z)}\|_{(B_\pi^p)^*}
    \leq
    \|\delta_z\|_{(B_\pi^q)^*}
    \|W_{\psi,\varphi}\|.
\]
By Theorem~\ref{thm:point-evaluations},
\[
    |\psi(z)|
    \frac{\Phi_{\pi,p}(|\operatorname{Im}\varphi(z)|)}
         {\Phi_{\pi,q}(|\operatorname{Im}z|)}
    \lesssim_{p,q}
    \|W_{\psi,\varphi}\|,
    \qquad z\in\C.
\]
Thus
\[
    M_{\pi,\pi;p,q}(\psi,\varphi)
    \lesssim_{p,q}
    \|W_{\psi,\varphi}\|.
\]

\smallskip
\noindent
\emph{Step 2. The symbol is affine.}
Choose \(m\in\N\) such that \(m>\frac{1}{p}\), and set
\[
    h:=G_{\pi/2,m,0}.
\]
By Lemma~\ref{lem:packet-estimates},
\(
    0\not\equiv h\in B_{\pi/2}^p.
\)
Choose \(0<\varepsilon<\frac{\pi}{2}\) and put
\[
    h_\varepsilon(z):=e^{-i\varepsilon z}h(z).
\]
Since
\(|h_\varepsilon(x)|=|h(x)|\) on \(\R\),
and \(h_\varepsilon\) has exponential type at most
\(\frac{\pi}{2}+\varepsilon<\pi\), we have
\(
    h,h_\varepsilon\in B_\pi^p.
\)
Therefore
\[
    H:=W_{\psi,\varphi}h
    =
    \psi\,(h\circ\varphi)
\qquad \text{and} \qquad
    H_\varepsilon
    :=
    W_{\psi,\varphi}h_\varepsilon
    =
    e^{-i\varepsilon\varphi}H
\]
belong to \(B_\pi^q\), and hence are of finite exponential type.

\smallskip
\noindent
Since \(\varphi\) is nonconstant and \(h\not\equiv0\), we have
\(h\circ\varphi\not\equiv0\). Thus \(H\not\equiv0\), because
\(\psi\not\equiv0\). 
Lemma~\ref{lem-product-rigidity}, applied to
\(
    H_\varepsilon
    =
    e^{-i\varepsilon\varphi}H,
\)
shows that \(e^{-i\varepsilon\varphi}\) is of finite exponential type.
This function is zero-free, so the Hadamard factorization theorem gives $e^{-i\varepsilon\varphi(z)}=e^{P(z)}$, where \(P\) is a polynomial of degree at most one. Hence
\[
    -i\varepsilon\varphi(z)-P(z)
    \in2\pi i\Z,
    \qquad z\in\C.
\]
The left-hand side is entire and takes values in a discrete set; it is
therefore constant. Thus \(\varphi\) is affine. Since \(\varphi\) is
nonconstant,
\[
    \varphi(z)=az+b,
    \qquad
    a,b\in\C,\quad a\neq0.
\]

\smallskip
\noindent
\emph{Step 3. The slope is real.}
We first note that \(\psi\) is of finite exponential type. Indeed,
\(h\circ\varphi\) is a nonzero entire function of finite exponential
type, while
\(
    \psi(h\circ\varphi)=H
\)
is also of finite exponential type. 
Lemma~\ref{lem-product-rigidity} therefore shows that \(\psi\) has
finite exponential type.

\smallskip
\noindent
Write
\[
    a=a_1+ia_2,
    \qquad
    b=b_1+ib_2,
\]
with \(a_j,b_j\in\R\). Suppose that \(a_2\neq0\). Taking \(z=x\in\R\)
in the estimate from Step~1 gives
\[
    |\psi(x)|
    \leq
    M_{\pi,\pi;p,q}(\psi,\varphi)
    \frac{\Phi_{\pi,q}(0)}
         {\Phi_{\pi,p}(|a_2x+b_2|)}.
\]
By the definition of \(\Phi_{\pi,p}\),
\[
    |\psi(x)|
    \lesssim_{p,q}
    M_{\pi,\pi;p,q}(\psi,\varphi)
    (1+|a_2x+b_2|)^{\frac{1}{p}}
    e^{-\pi|a_2x+b_2|}.
\]
Since \(a_2\neq0\), the polynomial factor can be absorbed into the
exponential. Thus there exist \(C_0,\gamma_0>0\) such that
\[
    |\psi(x)|
    \leq
    C_0e^{-\gamma_0|x|},
    \qquad x\in\R.
\]
Lemma~\ref{lem-exponential-decay-rigidity} then gives
\(\psi\equiv0\), a contradiction. Hence \(a_2=0\), and therefore
\(
    a\in\R.
\)

\smallskip
\noindent
\emph{Step 4. The slope satisfies \(|a|\leq1\).}
Suppose, to the contrary, that \(|a|>1\). Since \(a\in\R\), taking
\(z=iy\) in the estimate from Step~1 gives
\[
    |\psi(iy)|
    \leq
    M_{\pi,\pi;p,q}(\psi,\varphi)
    \frac{\Phi_{\pi,q}(|y|)}
         {\Phi_{\pi,p}(|ay+b_2|)}.
\]
Hence
\[
    |\psi(iy)|
    \lesssim_{p,q}
    M_{\pi,\pi;p,q}(\psi,\varphi)
    \frac{(1+|ay+b_2|)^{\frac{1}{p}}}
         {(1+|y|)^{\frac{1}{q}}}
    e^{\pi(|y|-|ay+b_2|)}.
\]
Since $|ay+b_2|-|y|
    \geq
    (|a|-1)|y|-|b_2|,$ the polynomial factor is again absorbed by the exponential. Thus, there
exist \(C_1,\gamma_1>0\) such that
\[
    |\psi(iy)|
    \leq
    C_1e^{-\gamma_1|y|},
    \qquad y\in\R.
\]
The entire function
\(
    g(z):=\psi(iz)
\)
has finite exponential type, so
Lemma~\ref{lem-exponential-decay-rigidity} yields \(g\equiv0\).
Hence \(\psi\equiv0\), again a contradiction. Therefore
\(
    0<|a|\leq1.
\)

\medskip
\noindent
\emph{Scaling back.}
We have proved the normalized assertion. For general \(\sigma,\tau\),
apply the normalized result to
\[
    W_{\widetilde\psi,\widetilde\varphi}:
    B_\pi^p\longrightarrow B_\pi^q,
\]
where
\[
    \widetilde\varphi(z)
    =
    \widetilde a z+\widetilde b,
    \qquad
    \widetilde a\in\R,
    \qquad
    0<|\widetilde a|\leq1.
\]
Hence
\[
    \varphi(z)=az+b,
    \qquad
    a\in\R,
    \qquad
    0<\sigma|a|\leq\tau.
\]
Finally, the scaling invariance of \(M_{\sigma,\tau;p,q}\) and the
norm identity give
\[
    M_{\sigma,\tau;p,q}(\psi,\varphi)
    \lesssim_{p,q}
    \|W_{\psi,\varphi}\|_{B_\sigma^p\to B_\tau^q}.
\]

\end{proof}

\subsection{The case \(0<p\leq q\leq\infty\)}
\label{subsec-wco-p-le-q}

In this range, Nikol'skii's inequality supplies the necessary gain in
integrability. Consequently, once the symbol has been reduced to a real
affine map, boundedness is governed by the residual bandwidth available
to the weight.

\begin{theorem}[Boundedness]
\label{thm-wco-p-le-q-bounded}
Let \(\sigma,\tau>0\) and \(0<p\leq q\leq\infty\), and let
$\psi,\varphi$ be entire functions such that
$\psi\not\equiv0$ and $\varphi$ is nonconstant. Then
\(
    W_{\psi,\varphi}:
    B_\sigma^p\longrightarrow B_\tau^q
\)
is bounded if and only if 
the following conditions hold:
\begin{itemize}
    \item[\textup{(a)}] \(\varphi(z)=az+b\),
where \(a\in\mathbb R\), \(0<\sigma|a|\leq \tau\), and \(b\in\mathbb C\).
\item[\textup{(b)}] \(\psi\in B_{\tau - \sigma |a|}^\infty\).
\end{itemize}
Moreover,
\[
    M_{\sigma,\tau;p,q}(\psi,\varphi)
    \lesssim \|W_{\psi,\varphi}\|_{B_\sigma^p\to B_\tau^q}
    \lesssim
    |a|^{-\frac{1}{q}}
    \sigma^{\frac{1}{p}-\frac{1}{q}}
    e^{\sigma|\operatorname{Im}b|}
    \|\psi\|_{L^\infty},
\]
where the implicit constants depend only on \(p\) and \(q\).
In particular, when \(p=q=\infty\),
\[
    \|W_{\psi,\varphi}\|_{B_\sigma^\infty\to B_\tau^\infty}
    =
    e^{\sigma|\operatorname{Im}b|}
    \|\psi\|_{L^\infty}.
\]
\end{theorem}

\begin{proof}
We first prove the assertion in the normalized setting
\[
    W_{\psi,\varphi}:B_\pi^p\longrightarrow B_\pi^q.
\]

\smallskip
\noindent
\emph{Necessity.}
Suppose that \(W_{\psi,\varphi}\) is bounded. By
Proposition~\ref{prop-wco-structural},
\[
    \varphi(z)=az+b,
    \qquad
    a\in\R,\quad 0<|a|\leq1,
\]
and
\[
    M_{\pi,\pi;p,q}(\psi,\varphi)
    \lesssim_{p,q}
    \|W_{\psi,\varphi}\|.
\]
In particular,
\[
    M:=M_{\pi,\pi;p,q}(\psi,\varphi)<\infty.
\]
Write \(b=b_1+ib_2\). For \(z=x+iy\), we have $\operatorname{Im}\varphi(z)=ay+b_2,$ and hence
\[
\begin{aligned}
    |\psi(x+iy)|
    &\leq
    M
    \frac{\Phi_{\pi,q}(|y|)}
         {\Phi_{\pi,p}(|ay+b_2|)}
\\
    &=
    M\pi^{\frac1q-\frac1p}
    \frac{(1+\pi|ay+b_2|)^{\frac1p}}
         {(1+\pi|y|)^{\frac1q}}
    e^{\pi(|y|-|ay+b_2|)}\\
    &\lesssim_{p,q,b}
    M(1+|y|)^{\frac1p}
    e^{\pi(1-|a|)|y|}.
\end{aligned}
\]
Therefore, for every \(\varepsilon>0\),
\[
    |\psi(z)|
    \lesssim_{\varepsilon,p,q,b}
    M e^{(\pi(1-|a|)+\varepsilon)|z|},
    \qquad z\in\C.
\]
Thus \(\psi\) has exponential type at most
\(\pi(1-|a|)\).

\smallskip
\noindent
Taking \(y=0\) in the defining estimate for \(M\), we also obtain
\[
    |\psi(x)|
    \leq
    M
    \frac{\Phi_{\pi,q}(0)}
         {\Phi_{\pi,p}(|b_2|)},
    \qquad x\in\R.
\]
Hence \(\psi|_{\R}\in L^\infty(\R)\), and consequently
\[
    \psi\in B_{\pi(1-|a|)}^\infty.
\]

\smallskip
\noindent
\emph{Sufficiency.}
Conversely, suppose that
\[
    \varphi(z)=az+b,
    \qquad
    a\in\R,\quad 0<|a|\leq1,
\]
and
\[
    \psi\in B_{\pi(1-|a|)}^\infty.
\]
If \(f\in B_\pi^p\), then \(f\circ\varphi\) has exponential type at
most \(\pi|a|\). Hence $W_{\psi,\varphi}f
    =
    \psi(f\circ\varphi)$ has exponential type at most
\[
    \pi(1-|a|)+\pi|a|=\pi.
\]
Moreover,
\[
    \|W_{\psi,\varphi}f\|_{L^q}
    \leq
    \|\psi\|_{L^\infty}
    \|f(a\cdot+b)\|_{L^q}.
\]
A change of variables gives
\[
    \|f(a\cdot+b)\|_{L^q}
    =
    |a|^{-\frac1q}
    \|f(\cdot+ib_2)\|_{L^q},
\]
with the usual interpretation when \(q=\infty\). By \eqref{eq:vertical-translation},
\[
    \|f(\cdot+ib_2)\|_{L^q}
    \lesssim_{p,q}
    \pi^{\frac1p-\frac1q}
    e^{\pi|b_2|}
    \|f\|_{L^p}.
\]
Therefore
\[
    \|W_{\psi,\varphi}f\|_{L^q}
    \lesssim_{p,q}
    |a|^{-\frac1q}
    \pi^{\frac1p-\frac1q}
    e^{\pi|\operatorname{Im}b|}
    \|\psi\|_{L^\infty}
    \|f\|_{L^p}.
\]
Thus \(W_{\psi,\varphi}:B_\pi^p\to B_\pi^q\) is bounded and
\[
    \|W_{\psi,\varphi}\|
    \lesssim_{p,q}
    |a|^{-\frac1q}
    \pi^{\frac1p-\frac1q}
    e^{\pi|\operatorname{Im}b|}
    \|\psi\|_{L^\infty}.
\]
We now return to arbitrary \(\sigma,\tau>0\). Under the scaling
\eqref{eq-symbols-scaling}, 
\[
    0<|\widetilde a|\leq1
    \quad\Longleftrightarrow\quad
    0<\sigma|a|\leq\tau.
\]
Moreover, from the definition of \(\widetilde\psi\),
\[
    \widetilde\psi
    \in
    B_{\pi(1-|\widetilde a|)}^\infty
    \quad\Longleftrightarrow\quad
    \psi
    \in
    B_{\tau-\sigma|a|}^\infty.
\]
Hence the normalized characterization gives precisely conditions
\textup{(a)}--\textup{(b)} in the general case.

\smallskip
\noindent
The lower estimate follows from
Proposition~\ref{prop-wco-structural}. For the upper estimate, by \eqref{eq-symbols-scaling}-\eqref{eq-wco-scaling}, we obtain
\[
\begin{aligned}
    \|W_{\psi,\varphi}\|
    =
    \|W_{\widetilde\psi,\widetilde\varphi}\|
    &\lesssim_{p,q}
    |\widetilde a|^{-\frac1q}
    \pi^{\frac1p-\frac1q}
    e^{\pi|\operatorname{Im}\widetilde b|}
    \|\widetilde\psi\|_{L^\infty}
\\
    &=
    |a|^{-\frac1q}
    \sigma^{\frac1p-\frac1q}
    e^{\sigma|\operatorname{Im}b|}
    \|\psi\|_{L^\infty}.
\end{aligned}
\]
Finally, let \(p=q=\infty\). Lemma~2.1 gives
\[
    \|W_{\psi,\varphi}f\|_{L^\infty}
    \leq
    e^{\sigma|\operatorname{Im}b|}
    \|\psi\|_{L^\infty}
    \|f\|_{L^\infty}.
\]
For the reverse inequality, take
\[
    f_0(z)
    :=
    \begin{cases}
        e^{-i\sigma z}, & \operatorname{Im}b\geq0,\\
        e^{i\sigma z},  & \operatorname{Im}b<0.
    \end{cases}
\]
Then \(f_0\in B_\sigma^\infty\), \(\|f_0\|_{L^\infty}=1\), and, since
\(a\in\R\),
\[
    |f_0(ax+b)|
    =
    e^{\sigma|\operatorname{Im}b|},
    \qquad x\in\R.
\]
Consequently,
\[
    \|W_{\psi,\varphi}f_0\|_{L^\infty}
    =
    e^{\sigma|\operatorname{Im}b|}
    \|\psi\|_{L^\infty},
\]
which proves the asserted exact formula.
\end{proof}

\noindent
Compactness requires one additional condition: the weight must vanish
along the real axis at infinity. This also forces strict inequality in
the bandwidth condition.

\begin{theorem}[Compactness]
\label{thm-wco-p-le-q-compact}
Let \(\sigma,\tau>0\) and \(0<p\leq q\leq\infty\), and let
$\psi,\varphi$ be entire functions such that
$\psi\not\equiv0$ and $\varphi$ is nonconstant. Then
\(
    W_{\psi,\varphi}:
    B_\sigma^p\longrightarrow B_\tau^q
\)
is compact if and only if the following conditions hold:
\begin{itemize}
\item[\textup{(a)}]
\(\varphi(z)=az+b\), 
where \(a\in\mathbb R\), \(0< \sigma |a|< \tau\), and \(b\in\mathbb C\);

\item[\textup{(b)}]
\(\psi\in B_{\tau - \sigma |a|}^\infty\) and 
\( \psi(x) \longrightarrow 0\) as \(|x|\to\infty\).
\end{itemize}
\end{theorem}

\begin{proof}
We first prove the assertion in the normalized setting
\[
    W_{\psi,\varphi}:B_\pi^p\longrightarrow B_\pi^q.
\]

\smallskip
\noindent
\emph{Necessity.}
Suppose that \(W_{\psi,\varphi}\) is compact. Then it is bounded, and
Theorem~\ref{thm-wco-p-le-q-bounded} gives
\[
    \varphi(z)=az+b,
    \qquad
    a\in\R,\quad 0<|a|\leq1,
\qquad \text{and} \qquad
    \psi\in B_{\pi(1-|a|)}^\infty.
\]
We first show that
\[
    \psi(x)\longrightarrow0
    \qquad (|x|\to\infty).
\]
Suppose not. Then there exist \(\varepsilon>0\) and a sequence
\((x_j)_{j\in\N}\subset\R\) such that
\[
    |x_j|\longrightarrow\infty,
    \qquad
    |\psi(x_j)|\geq\varepsilon.
\]
Put
\[
    w_j:=\varphi(x_j)=ax_j+b.
\]
Since \(a\in\R\setminus\{0\}\),
\[
    |\operatorname{Re}w_j|\longrightarrow\infty,
    \qquad
    \operatorname{Im}w_j=\operatorname{Im}b.
\]
Choose \(m\in\N\) such that \(m>\frac{1}{p}\), and define
\[
    F_j
    :=
    \frac{G_{\pi,m,w_j}}
         {\|G_{\pi,m,w_j}\|_{L^p}}.
\]
Then $\|F_j\|_{L^p}=1$ and Lemma~\ref{lem:packet-decay} gives
\[
    F_j\longrightarrow0
\]
uniformly on compact subsets of \(\C\). Since
\(W_{\psi,\varphi}\) is compact, Lemma~\ref{lem-compactness-criterion} yields
\[
    \|W_{\psi,\varphi}F_j\|_{L^q}
    \longrightarrow0.
\]
On the other hand, by Theorem~\ref{thm:point-evaluations}, applied at the real point \(x_j\),
\[
\begin{aligned}
    \|W_{\psi,\varphi}F_j\|_{L^q}
    &\gtrsim_q
    |(W_{\psi,\varphi}F_j)(x_j)|
=
    |\psi(x_j)|\,|F_j(w_j)|.
\end{aligned}
\]
By Lemma~\ref{lem:packet-estimates},
\[
    |F_j(w_j)|
    =
    \frac{|G_{\pi,m,w_j}(w_j)|}
         {\|G_{\pi,m,w_j}\|_{L^p}}
    \simeq_p
    \Phi_{\pi,p}(|\operatorname{Im}b|).
\]
Therefore
\[
    \|W_{\psi,\varphi}F_j\|_{L^q}
    \gtrsim_{p,q}
    \varepsilon\,
    \Phi_{\pi,p}(|\operatorname{Im}b|)
    >0,
\]
which is a contradiction. Thus
\[
    \psi(x)\longrightarrow0
    \qquad (|x|\to\infty).
\]
It remains to show that \(|a|<1\). If \(|a|=1\), then $\psi\in B_0^\infty=\C.$ Hence \(\psi\) is constant, and the preceding limit forces
\(\psi\equiv0\), contrary to our standing assumption. Therefore
\[
    0<|a|<1.
\]

\smallskip
\noindent
\emph{Sufficiency.}
Conversely, suppose that
\[
    \varphi(z)=az+b,
    \qquad
    a\in\R,\quad 0<|a|<1,
\qquad \text{and} \qquad
    \psi\in B_{\pi(1-|a|)}^\infty,
    \qquad
    \psi(x)\longrightarrow0
    \quad (|x|\to\infty).
\]
By Theorem~\ref{thm-wco-p-le-q-bounded},
\(W_{\psi,\varphi}:B_\pi^p\to B_\pi^q\) is bounded.

\smallskip
\noindent
Let \((f_j)_{j\in\N}\) be a bounded sequence in \(B_\pi^p\) such that $f_j\longrightarrow0$ uniformly on compact subsets of \(\C\). Write $b=b_1+ib_2$ and put
\[
    A:=\sup_{j\in\N}\|f_j\|_{L^p}<\infty.
\]
Suppose first that \(q<\infty\). Fix \(R>0\). Since the set $\{ax+b:|x|\leq R\}$ is compact in \(\C\), we observe
\[
    \sup_{|x|\leq R}|f_j(ax+b)|
    \longrightarrow0.
\]
As \(\psi\) is bounded on \(\R\),
\[
    \int_{|x|\leq R}
    |\psi(x)|^q|f_j(ax+b)|^q\,dx
    \longrightarrow0.
\]
For the tail,
\[
\begin{aligned}
    \int_{|x|>R}
    |\psi(x)|^q|f_j(ax+b)|^q\,dx
\leq
    \left(
        \sup_{|x|>R}|\psi(x)|
    \right)^q
    \|f_j(a\cdot+b)\|_{L^q}^q.
\end{aligned}
\]
A change of variables and \eqref{eq:vertical-translation} give
\[
\begin{aligned}
    \|f_j(a\cdot+b)\|_{L^q}
    &=
    |a|^{-\frac{1}{q}}
    \|f_j(\cdot+ib_2)\|_{L^q}
\lesssim_{p,q}
    |a|^{-\frac{1}{q}}
    \pi^{\frac{1}{p}-\frac{1}{q}}
    e^{\pi|b_2|}
    \|f_j\|_{L^p}.
\end{aligned}
\]
Consequently,
\[
  \int_{|x|>R}
    |\psi(x)|^q|f_j(ax+b)|^q\,dx
\lesssim_{p,q}
    |a|^{-1}
    \pi^{\frac{q}{p}-1}
    e^{q\pi|b_2|}
    A^q
    \left(
        \sup_{|x|>R}|\psi(x)|
    \right)^q.
\]
The right-hand side tends to \(0\) as \(R\to\infty\), uniformly in
\(j\). Hence
\[
    \|W_{\psi,\varphi}f_j\|_{L^q}
    \longrightarrow0.
\]
If \(q=\infty\), then for fixed \(R>0\),
\[
    \sup_{|x|\leq R}
    |\psi(x)|\,|f_j(ax+b)|
    \longrightarrow0.
\]
Moreover, by \eqref{eq:vertical-translation},
\[
    \sup_{|x|>R}
    |\psi(x)|\,|f_j(ax+b)|
    \leq
    \left(
        \sup_{|x|>R}|\psi(x)|
    \right)
    \|f_j(a\cdot+b)\|_{L^\infty}
\lesssim_p
    \pi^{\frac{1}{p}}
    e^{\pi|b_2|}
    A
    \sup_{|x|>R}|\psi(x)|.
\]
Again the last quantity tends to \(0\) as \(R\to\infty\), uniformly in
\(j\). Therefore
\[
    \|W_{\psi,\varphi}f_j\|_{L^\infty}
    \longrightarrow0.
\]
The uniqueness hypothesis in Lemma~\ref{lem-compactness-criterion} is automatic for
\(W_{\psi,\varphi}\). Indeed, if \(f_j\to0\) locally uniformly, then
\[
    W_{\psi,\varphi}f_j
    =
    \psi(f_j\circ\varphi)
    \longrightarrow0
\]
locally uniformly, while convergence in \(B_\pi^q\) implies local
uniform convergence by Theorem~\ref{thm:point-evaluations}. Hence Lemma~\ref{lem-compactness-criterion} shows that
\(W_{\psi,\varphi}:B_\pi^p\to B_\pi^q\) is compact.

\smallskip
\noindent
We now return to arbitrary \(\sigma,\tau>0\). By the scaling reduction
\eqref{eq-symbols-scaling}--\eqref{eq-wco-scaling}, the normalized conditions
\[
    0<|\widetilde a|<1,
    \qquad
    \widetilde\psi
    \in B_{\pi(1-|\widetilde a|)}^\infty,
    \qquad
    \widetilde\psi(x)\longrightarrow0
\]
are equivalent to
\[
    0<\sigma|a|<\tau,
    \qquad
    \psi\in B_{\tau-\sigma|a|}^\infty,
    \qquad
    \psi(x)\longrightarrow0
    \quad (|x|\to\infty).
\]
Since compactness is preserved under \eqref{eq-wco-scaling}, the general
statement follows.
\end{proof}

\noindent
The preceding compactness criterion suggests that the failure of
compactness should be measured by the asymptotic size of the weight on
the real axis. We now make this quantitative. For quasi-Banach spaces
\(X,Y\), let \(\mathcal K(X,Y)\) denote the space of compact operators
from \(X\) into \(Y\), and define
\[
    \|T\|_e
    :=
    \inf_{K\in\mathcal K(X,Y)}
    \|T-K\|_{X\to Y}.
\]
We continue to call this quantity the \emph{essential norm}, with the usual
quasi-Banach interpretation when necessary.

\begin{theorem}[Essential norm]
\label{thm-wco-essential-norm}
Let \(\sigma,\tau>0\) and \(0<p\leq q\leq\infty\). Suppose that
\(
    W_{\psi,\varphi}:
    B_\sigma^p\longrightarrow B_\tau^q
\)
is bounded, where \(\varphi(z)=az+b\) with \(a\in\R\) such that \(0<\sigma|a|\leq\tau\),
and put
\[
    L_\psi
    :=
    \limsup_{|x|\to\infty}|\psi(x)|.
\]
Then
\[
    \tau^{-\frac{1}{q}}
    \Phi_{\sigma,p}(|\operatorname{Im}b|)
    L_\psi
    \lesssim
    \|W_{\psi,\varphi}\|_e
    \lesssim
    |a|^{-\frac{1}{q}}
    \sigma^{\frac{1}{p}-\frac{1}{q}}
    e^{\sigma|\operatorname{Im}b|}
    L_\psi,
\]
where the implicit constants depend only on \(p\) and \(q\).
\end{theorem}

\begin{proof}
We first prove the estimates in the normalized setting
\[
    W_{\psi,\varphi}:B_\pi^p\longrightarrow B_\pi^q,
\]
where
\[
    \varphi(z)=az+b,
    \qquad
    a\in\R,\quad 0<|a|\leq1.
\]
By Theorem~\ref{thm-wco-p-le-q-bounded}, the function $\psi\in B_{\pi(1-|a|)}^\infty$. Write $b=b_1+ib_2.$

\smallskip
\noindent
\emph{Lower estimate.}
If \(L_\psi=0\), there is nothing to prove. Assume that \(L_\psi>0\).
Choose a sequence \((x_j)_{j\in\N}\subset\R\) such that
\[
    |x_j|\longrightarrow\infty,
    \qquad
    |\psi(x_j)|\longrightarrow L_\psi.
\]
Passing to a subsequence, we may assume that
\[
    |\operatorname{Re}w_j-\operatorname{Re}w_k|
    \longrightarrow\infty
    \qquad (j+k\to\infty,\ j\neq k),
\]
where
\[
    w_j:=\varphi(x_j)=ax_j+b.
\]
Notice that
\[
    \operatorname{Im}w_j=b_2,
    \qquad j\in\N.
\]
Choose \(m\in\N\), depending only on \(p\), such that $m>\frac1p$ and define the normalized packets
\[
    F_j
    :=
    \frac{G_{\pi,m,w_j}}
         {\|G_{\pi,m,w_j}\|_{L^p}}.
\]
Then $\|F_j\|_{L^p}=1$ and Lemma~\ref{lem:packet-decay} gives
\[
    F_j\longrightarrow0
\]
uniformly on compact subsets of \(\C\). Moreover, by Lemma~\ref{lem:packet-estimates},
\[
    |F_j(w_j)|
    \simeq_p
    \Phi_{\pi,p}(|b_2|).
\]
Let $K:B_\pi^p\longrightarrow B_\pi^q$ be a compact operator. Since \((F_j)_{j \in \N}\) is bounded, after passing to a further
subsequence we may assume that \((KF_j)_{j \in \N}\) converges in \(B_\pi^q\).
Relabeling the subsequence, put
\[
    H_j:=F_j-F_{j+1}.
\]
Then $\|H_j\|_{L^p}\lesssim_p1$ and $\|KH_j\|_{L^q}\longrightarrow0.$ The off-center decay of the packets also gives
\[
    F_{j+1}(w_j)\longrightarrow0.
\]
Indeed, the real parts of \(w_j\) and \(w_{j+1}\) separate while their
imaginary parts remain fixed. Consequently,
\[
    |H_j(w_j)|
    =
    |F_j(w_j)-F_{j+1}(w_j)|
    \gtrsim_p
    \Phi_{\pi,p}(|b_2|)
\]
for all sufficiently large \(j\).

\smallskip
\noindent
By Theorem~\ref{thm:point-evaluations}, point evaluation on the real axis in \(B_\pi^q\) satisfies
\[
    \|(W_{\psi,\varphi}-K)H_j\|_{L^q}
    \gtrsim_q
    \pi^{-\frac{1}{q}}
    |((W_{\psi,\varphi}-K)H_j)(x_j)|
    \geq
    \pi^{-\frac{1}{q}}
    \Bigl(
        |\psi(x_j)|\,|H_j(w_j)|
        -
        |KH_j(x_j)|
    \Bigr).
\]
Again by Theorem~\ref{thm:point-evaluations},
\[
    |KH_j(x_j)|
    \lesssim_q
    \pi^{\frac{1}{q}}\|KH_j\|_{L^q}
    \longrightarrow0.
\]
It follows that
\[
    \limsup_{j\to\infty}
    \|(W_{\psi,\varphi}-K)H_j\|_{L^q}
    \gtrsim_{p,q}
    \pi^{-\frac{1}{q}}
    \Phi_{\pi,p}(|b_2|)
    L_\psi.
\]
Since \(\|H_j\|_{L^p}\lesssim_p1\),
\[
    \|W_{\psi,\varphi}-K\|
    \gtrsim_{p,q}
    \pi^{-\frac{1}{q}}
    \Phi_{\pi,p}(|b_2|)
    L_\psi.
\]
Taking the infimum over all compact
\(K:B_\pi^p\to B_\pi^q\), we obtain
\[
    \|W_{\psi,\varphi}\|_e
    \gtrsim_{p,q}
    \pi^{-\frac{1}{q}}
    \Phi_{\pi,p}(|b_2|)
    L_\psi.
\]

\smallskip
\noindent
\emph{Upper estimate.}
For \(0<\lambda<1\), set
\[
    \varphi_\lambda(z):=\lambda az+b
\qquad \text{and} \qquad
    \varepsilon_\lambda
    :=
    \pi|a|(1-\lambda).
\]
Define
\[
    h_\lambda(z)
    :=
    \frac{\sin(\varepsilon_\lambda z)}
         {\varepsilon_\lambda z},
    \qquad
    h_\lambda(0):=1,
\qquad \text{and} \qquad
    \psi_\lambda:=\psi h_\lambda.
\]
Since $\operatorname{type}(\psi)
    \leq
    \pi(1-|a|)$ and $\operatorname{type}(h_\lambda)
    \leq
    \varepsilon_\lambda$, we have
\[
\begin{aligned}
    \operatorname{type}(\psi_\lambda)
    &\leq
    \pi(1-|a|)
    +
    \pi|a|(1-\lambda)
=
    \pi(1-\lambda|a|).
\end{aligned}
\]
Moreover,
\[
    |h_\lambda(x)|\leq1,
    \qquad x\in\R,
\]
and, for each fixed \(0<\lambda<1\),
\[
    h_\lambda(x)\longrightarrow0
    \qquad (|x|\to\infty).
\]
Since \(\psi\in L^\infty(\R)\),
\[
    \psi_\lambda(x)\longrightarrow0
    \qquad (|x|\to\infty).
\]
Thus Theorem~\ref{thm-wco-p-le-q-compact} shows that
\[
    W_{\psi_\lambda,\varphi_\lambda}
    :
    B_\pi^p\longrightarrow B_\pi^q
\]
is compact.

\smallskip
\noindent
We estimate $W_{\psi,\varphi}-W_{\psi_\lambda,\varphi_\lambda}$ as follows. For \(f\in B_\pi^p\) and \(x\in\R\),
\[
\begin{aligned}
    (W_{\psi,\varphi}-W_{\psi_\lambda,\varphi_\lambda})f(x)
    &=
    \psi(x)
    \bigl(
        f(ax+b)-h_\lambda(x)f(\lambda ax+b)
    \bigr)
\\
    &=
    \psi(x)
    \bigl(
        f(ax+b)-f(\lambda ax+b)
    \bigr)
+
    \psi(x)(1-h_\lambda(x))
    f(\lambda ax+b).
\end{aligned}
\]
Fix \(R>0\) and restrict for the moment to $\frac12\leq\lambda<1.$ Set
\[
    \Omega_R
    :=
    \{
        tax+b:
        |x|\leq R,\ \tfrac12\leq t\leq1
    \},
\]
which is a compact subset of \(\C\). Theorem~\ref{thm:point-evaluations} and Cauchy's estimate
give
\[
    \sup_{\|f\|_{L^p}\leq1}
    \sup_{z\in\Omega_R}
    \bigl(|f(z)|+|f'(z)|\bigr)
    <\infty.
\]
Consequently,
\[
    \sup_{\|f\|_{L^p}\leq1}
    \sup_{|x|\leq R}
    |f(ax+b)-f(\lambda ax+b)|
    \longrightarrow0
    \qquad (\lambda\to1^-).
\]
Also,
\[
    \sup_{|x|\leq R}|1-h_\lambda(x)|
    \longrightarrow0
    \qquad (\lambda\to1^-).
\]
Since \(\psi\) is bounded on \([-R,R]\), it follows that
\[
    \sup_{\|f\|_{L^p}\leq1}
    \|(W_{\psi,\varphi}-W_{\psi_\lambda,\varphi_\lambda})f\|_
        {L^q(|x|\leq R)}
    \longrightarrow0
    \qquad (\lambda\to1^-).
\]
For the tail, put
\[
    L_R
    :=
    \sup_{|x|>R}|\psi(x)|.
\]
Since \(|h_\lambda(x)|\leq1\) on \(\R\),
\[
\begin{aligned}
    |(W_{\psi,\varphi}-W_{\psi_\lambda,\varphi_\lambda})f(x)|
    &\leq
    L_R
    \bigl(
        |f(ax+b)|
        +
        |f(\lambda ax+b)|
    \bigr),
    \qquad |x|>R.
\end{aligned}
\]
Let $s:=\min\{1,q\},$ with \(s=1\) when \(q=\infty\). Using the \(s\)-triangle inequality
and \eqref{eq:vertical-translation}, we obtain
\[
\begin{aligned}
    \|(W_{\psi,\varphi}-W_{\psi_\lambda,\varphi_\lambda})f\|_
        {L^q(|x|>R)}^s
    & \lesssim_q
    L_R^s
    \left(
        \|f(a\cdot+b)\|_{L^q}^s
        +
        \|f(\lambda a\cdot+b)\|_{L^q}^s
    \right)
\\
    &\lesssim_{p,q}
    L_R^s
    |a|^{-s/q}
    \pi^{s(\frac{1}{p}-\frac{1}{q})}
    e^{s\pi|b_2|}
    \|f\|_{L^p}^s,
\end{aligned}
\]
where in the last step we used \(\lambda\geq \frac{1}{2}\). Hence
\[
    \sup_{\|f\|_{L^p}\leq1}
    \|(W_{\psi,\varphi}-W_{\psi_\lambda,\varphi_\lambda})f\|_
        {L^q(|x|>R)}
    \lesssim_{p,q}
    |a|^{-\frac{1}{q}}
    \pi^{\frac{1}{p}-\frac{1}{q}}
    e^{\pi|b_2|}
    L_R.
\]
Combining the compact and tail parts, and then letting
\(\lambda\to1^-\), gives
\[
    \|W_{\psi,\varphi}\|_e
    \lesssim_{p,q}
    |a|^{-\frac{1}{q}}
    \pi^{\frac{1}{p}-\frac{1}{q}}
    e^{\pi|b_2|}
    L_R.
\]
Since \(\psi\in L^\infty(\R)\),
\[
    L_R\longrightarrow L_\psi
    \qquad (R\to\infty).
\]
Therefore
\[
    \|W_{\psi,\varphi}\|_e
    \lesssim_{p,q}
    |a|^{-\frac{1}{q}}
    \pi^{\frac{1}{p}-\frac{1}{q}}
    e^{\pi|b_2|}
    L_\psi.
\]
We have thus proved, in the normalized setting,
\[
    \pi^{-\frac{1}{q}}
    \Phi_{\pi,p}(|\operatorname{Im}b|)
    L_\psi
    \lesssim
    \|W_{\psi,\varphi}\|_e
    \lesssim
    |a|^{-\frac{1}{q}}
    \pi^{\frac{1}{p}-\frac{1}{q}}
    e^{\pi|\operatorname{Im}b|}
    L_\psi.
\]

\smallskip
\noindent
\emph{Scaling back.}
Let
\[
    W_{\widetilde\psi,\widetilde\varphi}
    =
    \mathcal U_\tau
    W_{\psi,\varphi}
    \mathcal U_\sigma^{-1}
    :
    B_\pi^p\longrightarrow B_\pi^q.
\]
Since the scaling maps are isometric isomorphisms and isometric
conjugation preserves compact operators, we obtain
\[
\begin{aligned}
    \|W_{\psi,\varphi}\|_e=\|W_{\widetilde\psi,\widetilde\varphi}\|_e
    &\gtrsim
    \pi^{-\frac{1}{q}}
    \Phi_{\pi,p}(|\operatorname{Im}\widetilde b|)
    L_{\widetilde\psi}
=
    \tau^{-\frac{1}{q}}
    \Phi_{\sigma,p}(|\operatorname{Im}b|)
    L_\psi.
\end{aligned}
\]
Similarly, the normalized upper estimate gives
\[
\begin{aligned}
    \|W_{\psi,\varphi}\|_e
    &\lesssim
    |\widetilde a|^{-\frac{1}{q}}
    \pi^{\frac{1}{p}-\frac{1}{q}}
    e^{\pi|\operatorname{Im}\widetilde b|}
    L_{\widetilde\psi}
=
    |a|^{-\frac{1}{q}}
    \sigma^{\frac{1}{p}-\frac{1}{q}}
    e^{\sigma|\operatorname{Im}b|}
    L_\psi.
\end{aligned}
\]
This proves the assertion.
\end{proof}

\begin{remark}
Under the assumptions of Theorem~\ref{thm-wco-essential-norm},
\[
    \|W_{\psi,\varphi}\|_e=0
    \quad\Longleftrightarrow\quad
    \psi(x)\longrightarrow0
    \quad (|x|\to\infty),
\]
in agreement with Theorem~\ref{thm-wco-p-le-q-compact}.
\end{remark}

\subsection{The case \(0<q<p\leq\infty\)}
\label{subsec-wco-q-lt-p}

The situation changes when \(q<p\). Nikol'skii's inequality can no
longer provide the required improvement of integrability; instead, the
weight itself must supply the missing integrability. This leads naturally
to the exponent
\[
    \frac{pq}{p-q}.
\]
In this range, boundedness already implies compactness.

\begin{theorem}
\label{thm-wco-q-lt-p}
Let \(\sigma,\tau>0\) and \(0<q<p\leq\infty\), and let
$\psi,\varphi$ be entire functions such that
$\psi\not\equiv0$ and $\varphi$ is nonconstant.
Then the following conditions are equivalent:
\begin{enumerate}
\item[\textnormal{(a)}]
\(
    W_{\psi,\varphi}:
    B_\sigma^p\longrightarrow B_\tau^q
\)
is bounded;

\item[\textnormal{(b)}]
\(
    W_{\psi,\varphi}:
    B_\sigma^p\longrightarrow B_\tau^q
\)
is compact;

\item[\textnormal{(c)}]
\(
    \varphi(z)=az+b,
\)
where
\(a\in\R\) such that \(
    0<\sigma|a|<\tau\), and \(b\in\C\),
and
\(\psi\in B_{\tau-\sigma|a|}^{\frac{pq}{p-q}}.
\)
\end{enumerate}
Moreover,
\[
    \left(\frac{\sigma}{\tau}\right)^{\frac{1}{p}}
    e^{-\sigma|\operatorname{Im}b|}
    \|\psi\|_{L^{\frac{pq}{p-q}}}
    \lesssim
    \|W_{\psi,\varphi}\|_{B_\sigma^p\to B_\tau^q}
    \le
    |a|^{-\frac{1}{p}}
    e^{\sigma|\operatorname{Im}b|}
    \|\psi\|_{L^{\frac{pq}{p-q}}},
\]
where the implicit constants depend only on \(p\) and \(q\).
\end{theorem}

\begin{proof}
We first prove the assertion in the normalized setting
\[
    W_{\psi,\varphi}:B_\pi^p\longrightarrow B_\pi^q.
\]

\smallskip
\noindent
(a)$\Longrightarrow$(c). Suppose that \(W_{\psi,\varphi}\) is bounded. By
Proposition~\ref{prop-wco-structural},
\[
    \varphi(z)=az+b,
    \qquad
    a\in\R,\quad 0<|a|\leq1,
\]
and
\[
    M:=M_{\pi,\pi;p,q}(\psi,\varphi)<\infty.
\]
Write \(b=b_1+ib_2\). For \(z=x+iy\), the imaginary part $\operatorname{Im}\varphi(z)=ay+b_2$ and hence
\[
    |\psi(x+iy)|
    \leq
    M
    \frac{\Phi_{\pi,q}(|y|)}
         {\Phi_{\pi,p}(|ay+b_2|)}
=
    M\pi^{\frac1q-\frac1p}
    \frac{(1+\pi|ay+b_2|)^{\frac{1}{p}}}
         {(1+\pi|y|)^{\frac{1}{q}}}
    e^{\pi(|y|-|ay+b_2|)}.
\]
Since $|ay+b_2|\geq |a||y|-|b_2|$, the polynomial factor can be absorbed into an arbitrarily small
exponential. It follows that
\begin{align}\label{202608182159}
\operatorname{type}(\psi)\leq\pi(1-|a|).
\end{align}
Next, we show that \(|a|<1\). If \(|a|=1\), then the preceding estimate
gives
\[
    |\psi(x+iy)|
    \lesssim_{p,q,b}
    M(1+|y|)^{\frac1p-\frac1q},
    \qquad x,y\in\R.
\]
Since \(q<p\),
\[
    \frac1p-\frac1q<0.
\]
Thus \(\psi\) is bounded on \(\C\), and hence constant by Liouville's
theorem. Letting \(|y|\to\infty\) in the preceding estimate shows that
this constant is zero, contrary to \(\psi\not\equiv0\). Therefore
\[
    0<|a|<1.
\]
Assume first that \(p<\infty\). Define the positive measure
\[
    d\mu(t)
    :=
    \frac1{|a|}
    \left|
        \psi\left(\frac{t-b_1}{a}\right)
    \right|^qdt.
\]
For \(g\in B_\pi^p\), put \(f=T_{-b_2}g\). Then $f(t+ib_2)=g(t),$ and consequently
\[
    \|g\|_{L^q(\mu)}
    =
    \|W_{\psi,\varphi}f\|_{L^q}
    \leq
    \|W_{\psi,\varphi}\|\cdot
    \|T_{-b_2}g\|_{L^p}.
\]
By Lemma~2.1,
\[
    \|\mu\|_{\mathcal C_{p,q}(\pi)}
    \leq
    e^{\pi|b_2|}
    \|W_{\psi,\varphi}\|.
\]
Theorem~\ref{thm-carleson-q-lt-p} therefore gives
\[
    \pi^{\frac{1}{p}}
    \|(\mu_{\pi,n})_{n\in\Z}\|_
        {\ell^{(\frac pq)'}}^{\frac{1}{q}}
    \lesssim_{p,q}
    e^{\pi|b_2|}
    \|W_{\psi,\varphi}\|.
\]
For each \(n\in\Z\),
\[
    \mu_{\pi,n}
    =
    \int_{J_n}|\psi(u)|^q\,du,
\]
where \((J_n)_{n\in\Z}\) is, up to a reversal of the indexing, a real
translate of the partition
\[
    \bigl(I_{\pi|a|,n}\bigr)_{n\in\Z}.
\]
Since \(\psi\) has exponential type at most
\(\pi(1-|a|)\), Lemma~\ref{lem:local-q-to-p}, applied after the corresponding real
translation with the exponents
\[
    q
    \quad\text{and}\quad
    \frac{pq}{p-q},
\]
gives
\[
    \|\psi\|_{L^{\frac{pq}{p-q}}}
    \lesssim_{p,q}
    \pi^{\frac{1}{p}}
    \|(\mu_{\pi,n})_{n\in\Z}\|_
        {\ell^{(\frac pq)'}}^{\frac{1}{q}}.
\]
Consequently,
\[
    e^{-\pi|b_2|}
    \|\psi\|_{L^{\frac{pq}{p-q}}}
    \lesssim_{p,q}
    \|W_{\psi,\varphi}\|.
\]
Together with the type estimate,
\[
    \psi\in
    B_{\pi(1-|a|)}^{\frac{pq}{p-q}}.
\]
Finally, suppose that \(p=\infty\). Since
\(1\in B_\pi^\infty\) and \(W_{\psi,\varphi}1=\psi\), we have $\|\psi\|_{L^q}
    \leq
    \|W_{\psi,\varphi}\|.$ Together with \eqref{202608182159}, this gives $\psi\in B_{\pi(1-|a|)}^q,$ which is precisely the asserted condition.

\smallskip
\noindent
(c)$\Longrightarrow$(a). Suppose that
\[
    \varphi(z)=az+b,
    \qquad
    a\in\R,\quad0<|a|<1,
\qquad \text{and} \qquad
    \psi\in
    B_{\pi(1-|a|)}^{\frac{pq}{p-q}}.
\]
For \(f\in B_\pi^p\), the function
\(\psi(f\circ\varphi)\) has exponential type at most \(\pi\).
Moreover, H\"older's inequality gives
\[
    \|W_{\psi,\varphi}f\|_{L^q}
    \leq
    \|\psi\|_{L^{\frac{pq}{p-q}}}
    \|f(a\cdot+b)\|_{L^p}.
\]
By a change of variables and Lemma~2.1,
\[
    \|f(a\cdot+b)\|_{L^p}
    =
    |a|^{-\frac{1}{p}}
    \|f(\cdot+ib_2)\|_{L^p}
    \leq
    |a|^{-\frac{1}{p}}e^{\pi|b_2|}
    \|f\|_{L^p}.
\]
Therefore
\[
    \|W_{\psi,\varphi}\|
    \lesssim
    |a|^{-\frac{1}{p}}e^{\pi|b_2|}
    \|\psi\|_{L^{\frac{pq}{p-q}}}.
\]
It remains to prove (c)$\Longrightarrow$(b). Let \((f_j)_{j\in\N}\) be bounded in
\(B_\pi^p\) and satisfy \(f_j\to0\) uniformly on compact subsets of
\(\C\). For \(R>0\), H\"older's inequality gives
\[
\|W_{\psi,\varphi}f_j\|_{L^q}^q
    \leq
    \|\psi\|_{L^{\frac{pq}{p-q}}([-R,R])}^q
    \|f_j(a\cdot+b)\|_{L^p([-R,R])}^q
+
    \|\psi\|_{L^{\frac{pq}{p-q}}(\{|x|>R\})}^q
    \|f_j(a\cdot+b)\|_{L^p}^q.
\]
For fixed \(R\), the first term tends to zero because
\(f_j\to0\) locally uniformly. By Lemma~\ref{lem:vertical-translation},
\[
    \|f_j(a\cdot+b)\|_{L^p}
    \leq
    |a|^{-\frac{1}{p}}e^{\pi|\operatorname{Im}b|}
    \|f_j\|_{L^p},
\]
with the usual interpretation when \(p=\infty\). Hence the second
term tends to zero uniformly in \(j\) as \(R\to\infty\), since
\(\psi\in L^{\frac{pq}{p-q}}(\R)\). Thus
\[
    \|W_{\psi,\varphi}f_j\|_{L^q}\longrightarrow0.
\]
As in the proof of Theorem~\ref{thm-wco-p-le-q-compact}, the uniqueness
hypothesis in Lemma~\ref{lem-compactness-criterion} is automatic for \(W_{\psi,\varphi}\).
Therefore \(W_{\psi,\varphi}\) is compact.

\medskip
\noindent
The general case follows from the scaling reduction
\eqref{eq-symbols-scaling}--\eqref{eq-wco-scaling}. Indeed,
\[
    \|\widetilde\psi\|_
        {L^{\frac{pq}{p-q}}}
    =
    \left(\frac{\sigma}{\tau}\right)^{\frac{1}{p}}
    \|\psi\|_{L^{\frac{pq}{p-q}}},\quad \pi|\operatorname{Im}\widetilde b|
    =
    \sigma|\operatorname{Im}b|.
\]
The normalized conditions and estimates therefore give
\[
    0<\sigma|a|<\tau,
    \qquad
    \psi\in
    B_{\tau-\sigma|a|}^{\frac{pq}{p-q}},
\]
together with the asserted two-sided norm estimate.
\end{proof}

\begin{remark}
In the range \(0<q<p\leq\infty\),
\[
    W_{\psi,\varphi}:B_\sigma^p\to B_\tau^q
    \text{ is bounded}
    \quad\Longleftrightarrow\quad
    W_{\psi,\varphi}:B_\sigma^p\to B_\tau^q
    \text{ is compact}.
\]
Consequently, there is no nontrivial essential-norm problem in this
range.
\end{remark}

\subsection{Composition operators}
\label{subsec-composition-operators}

We first specialize the preceding results to the unweighted case
\(\psi\equiv1\). Besides giving the boundedness and compactness
classification, this case provides a sharp quantitative model for the
general theory. For affine symbols, the operator norm reduces exactly
to the norm of a vertical translation, yielding a sharp scale-explicit
estimate for \(0<p\leq q\leq\infty\) and an exact formula when \(p=q\).

\begin{corollary}
\label{cor-composition-classification}
Let \(\sigma,\tau>0\), \(0<p,q\leq\infty\), and let \(\varphi\) be
a nonconstant entire function.

\begin{enumerate}
\item[\textnormal{(a)}]
Let \(0<p\leq q\leq\infty\). Then
\(
    C_\varphi:
    B_\sigma^p\longrightarrow B_\tau^q
\)
is bounded if and only if
\(    \varphi(z)=az+b,
\)
where
\(a\in\R\) with \(
    0<\sigma|a|\leq\tau\) and \(
    b\in\C\).
    
No such nonconstant composition operator is compact.

\item[\textnormal{(b)}]
Let \(0<q<p\leq\infty\). Then there is no bounded nonconstant
composition operator
\(
    C_\varphi:
    B_\sigma^p\longrightarrow B_\tau^q.
\)
\end{enumerate}
\end{corollary}

\begin{proof}
Apply Theorems~\ref{thm-wco-p-le-q-bounded},
\ref{thm-wco-p-le-q-compact}, and
\ref{thm-wco-q-lt-p} with \(\psi\equiv1\).

\smallskip
\noindent
If \(p\leq q\), the condition $1\in B_{\tau-\sigma|a|}^\infty$ is automatic whenever \(\sigma|a|\leq\tau\), while the compactness
condition \(1\to0\) at infinity never holds.

\smallskip
\noindent
If \(q<p\), Theorem~\ref{thm-wco-q-lt-p} would require $1\in L^{\frac{pq}{p-q}}(\R)$, which is
impossible.
\end{proof}

\begin{remark}\label{rem-constant-composition}
If \(\varphi\equiv b\), then
\(
    C_\varphi f=f(b)\mathbf 1.
\)
Thus
\(
    C_\varphi:
    B_\sigma^p\longrightarrow B_\tau^q
\)
is bounded if and only if \(q=\infty\). In this case it has rank one and
\[
    \|C_\varphi\|_{B_\sigma^p\to B_\tau^\infty}
    =
    \|\delta_b\|_{(B_\sigma^p)^*}
    \simeq
    \Phi_{\sigma,p}(|\operatorname{Im}b|),
\]
where the implicit constants depend only on \(p\).
\end{remark}

\smallskip
\noindent
We now determine the norm of the bounded affine composition operators.

\begin{theorem}[Affine composition operators]
\label{thm-composition-norm}
Let \(\sigma,\tau>0\), \(0<p\leq q\leq\infty\), and
\(\varphi(z)=az+b\) with \(
    a\in\R\) such that \(
    0<\sigma|a|\leq\tau\).
Then
\[
    \|C_\varphi\|_{B_\sigma^p\to B_\tau^q}
    \simeq
    |a|^{-\frac{1}{q}}
    \sigma^{\frac{1}{p}-\frac{1}{q}}
    \frac{e^{\sigma|\operatorname{Im}b|}}
         {(1+\sigma|\operatorname{Im}b|)^{\frac{1}{p}-\frac{1}{q}}},
\]
where the implicit constants depend only on \(p\) and \(q\).
In particular, when \(p=q\),
\[
    \|C_\varphi\|_{B_\sigma^p\to B_\tau^p}
    =
    |a|^{-\frac{1}{p}}
    e^{\sigma|\operatorname{Im}b|}.
\]
\end{theorem}

\begin{proof}
Write $b=b_1+ib_2$ and set
\[
    \alpha:=\frac1p-\frac1q.
\]
For \(f\in B_\sigma^p\), the function \(f\circ\varphi\) has exponential
type at most $\sigma|a|\leq\tau.$ Moreover, if \(q<\infty\), a change of variables gives
\[
\begin{aligned}
    \|C_\varphi f\|_{L^q}^q
    &=
    \int_{\R}|f(ax+b_1+ib_2)|^q\,dx
=
    \frac1{|a|}
    \int_{\R}|f(u+ib_2)|^q\,du.
\end{aligned}
\]
Hence
\[
    \|C_\varphi f\|_{L^q}
    =
    |a|^{-\frac{1}{q}}
    \|T_{b_2}f\|_{L^q}.
\]
For \(q=\infty\), since \(a\R+b_1=\R\),
\[
    \|C_\varphi f\|_{L^\infty}
    =
    \|T_{b_2}f\|_{L^\infty},
\]
so the same identity holds.
Taking the supremum over nonzero \(f\in B_\sigma^p\), we obtain
\begin{equation}\label{eq:composition-translation-reduction}
    \|C_\varphi\|_{B_\sigma^p\to B_\tau^q}
    =
    |a|^{-\frac{1}{q}}
    \|T_{b_2}\|_{B_\sigma^p\to B_\sigma^q}.
\end{equation}
It remains to estimate the vertical translation. Put $y:=|b_2|.$

\smallskip
\noindent
Suppose first that \(p<q\). By
Lemma~\ref{lem:vertical-translation},
\[
    \|T_{b_2}f\|_{L^p}
    \leq
    e^{\sigma y}\|f\|_{L^p}.
\]
On the other hand, Theorem~\ref{thm:point-evaluations} gives
\[
\begin{aligned}
    \|T_{b_2}f\|_{L^\infty}
    &=
    \sup_{x\in\R}|f(x+ib_2)|
\lesssim_p
    \Phi_{\sigma,p}(y)\|f\|_{L^p}
=
    \sigma^{\frac{1}{p}}
    \frac{e^{\sigma y}}
         {(1+\sigma y)^{\frac{1}{p}}}
    \|f\|_{L^p}.
\end{aligned}
\]
If \(q<\infty\), then
\[
    \|g\|_{L^q}
    \leq
    \|g\|_{L^p}^{\frac{p}{q}}
    \|g\|_{L^\infty}^{1-\frac{p}{q}},
\]
and therefore
\[
\begin{aligned}
    \|T_{b_2}f\|_{L^q}
    &\lesssim_{p,q}
    \bigl(e^{\sigma y}\|f\|_{L^p}\bigr)^{\frac{p}{q}}
    \left(
        \sigma^{\frac{1}{p}}
        \frac{e^{\sigma y}}
             {(1+\sigma y)^{\frac{1}{p}}}
        \|f\|_{L^p}
    \right)^{1-\frac{p}{q}}
=
    \sigma^\alpha
    \frac{e^{\sigma y}}
         {(1+\sigma y)^\alpha}
    \|f\|_{L^p}.
\end{aligned}
\]
For \(q=\infty\), the same estimate follows directly from the
point-evaluation bound. Thus
\begin{equation}\label{eq:vertical-pq-upper}
    \|T_{b_2}\|_{B_\sigma^p\to B_\sigma^q}
    \lesssim_{p,q}
    \sigma^\alpha
    \frac{e^{\sigma y}}
         {(1+\sigma y)^\alpha}.
\end{equation}
We prove the reverse estimate by adapting the test functions used in
the proof of the corresponding result in the equal-bandwidth case.
Choose \(N\in\N\) such that $N>\frac1p$ and set
\[
    H_N(z)
    :=
    \left(
        \frac{\sin(z/N)}
             {z/N}
    \right)^N,
    \qquad H_N(0):=1.
\]
For \(0<\varepsilon<\sigma\), put $g_\varepsilon(z):=H_N(\varepsilon z).$ Then \(g_\varepsilon\) has exponential type at most \(\varepsilon\),
and
\[
    \|g_\varepsilon\|_{L^s}
    =
    \varepsilon^{-\frac{1}{s}}\|H_N\|_{L^s},
    \qquad 0<s<\infty,
\]
with the usual interpretation when \(s=\infty\).

\smallskip
\noindent
Define
\[
    f_\varepsilon(z)
    :=
    \begin{cases}
        e^{-i(\sigma-\varepsilon)z}g_\varepsilon(z),
            & b_2\geq0,\\[2mm]
        e^{i(\sigma-\varepsilon)z}g_\varepsilon(z),
            & b_2<0.
    \end{cases}
\]
Then \(f_\varepsilon\in B_\sigma^p\) and
\[
    \|f_\varepsilon\|_{L^p}
    =
    \varepsilon^{-\frac{1}{p}}\|H_N\|_{L^p}.
\]
Moreover,
\[
    \|T_{b_2}f_\varepsilon\|_{L^q}
    =
    e^{(\sigma-\varepsilon)y}
    \varepsilon^{-\frac{1}{q}}
    \|H_N(\cdot+i\varepsilon b_2)\|_{L^q}.
\]
Consequently,
\begin{equation}\label{eq:vertical-test-ratio}
    \frac{\|T_{b_2}f_\varepsilon\|_{L^q}}
         {\|f_\varepsilon\|_{L^p}}
    =
    e^{(\sigma-\varepsilon)y}
    \varepsilon^\alpha
    \frac{\|H_N(\cdot+i\varepsilon b_2)\|_{L^q}}
         {\|H_N\|_{L^p}}.
\end{equation}
Choose
\[
    \varepsilon
    :=
    \frac{\sigma}{2(1+\sigma y)}.
\]
Then
\[
    0<\varepsilon<\sigma,
    \qquad
    \varepsilon y\leq\frac12,
    \qquad
    e^{-\varepsilon y}\geq e^{-1/2}.
\]
If \(q<\infty\), the polynomial decay of \(H_N\) on horizontal strips
and dominated convergence show that
\[
    s\longmapsto
    \|H_N(\cdot+is)\|_{L^q}
\]
is continuous on \([-1/2,1/2]\). Since none of these functions is
identically zero,
\[
    \inf_{|s|\leq1/2}
    \|H_N(\cdot+is)\|_{L^q}>0.
\]
If \(q=\infty\), then
\[
    \|H_N(\cdot+is)\|_{L^\infty}
    \geq
    |H_N(is)|
    \geq1,
    \qquad |s|\leq\frac12.
\]
Thus \eqref{eq:vertical-test-ratio} yields
\[
\begin{aligned}
    \|T_{b_2}\|_{B_\sigma^p\to B_\sigma^q}
    &\gtrsim_{p,q}
    e^{\sigma y}\varepsilon^\alpha
\simeq_{p,q}
    \sigma^\alpha
    \frac{e^{\sigma y}}
         {(1+\sigma y)^\alpha}.
\end{aligned}
\]
Together with \eqref{eq:vertical-pq-upper}, this gives
\[
    \|T_{b_2}\|_{B_\sigma^p\to B_\sigma^q}
    \simeq
    \sigma^{\frac1p-\frac1q}
    \frac{e^{\sigma|b_2|}}
         {(1+\sigma|b_2|)^{\frac1p-\frac1q}},
\]
whenever \(p<q\).

\smallskip
\noindent
If \(p=q\), Lemma~\ref{lem:vertical-translation} gives directly $\|T_{b_2}\|_{B_\sigma^p\to B_\sigma^p}
    =
    e^{\sigma|b_2|}.$ Combining these estimates with
\eqref{eq:composition-translation-reduction}, we obtain
\[
    \|C_\varphi\|_{B_\sigma^p\to B_\tau^q}
    \simeq
    |a|^{-\frac{1}{q}}
    \sigma^{\frac1p-\frac1q}
    \frac{e^{\sigma|\operatorname{Im}b|}}
         {(1+\sigma|\operatorname{Im}b|)^{\frac{1}{p}-\frac{1}{q}}},
\]
with constants depending only on \(p\) and \(q\).

\smallskip
\noindent
Finally, when \(p=q\),
\[
\begin{aligned}
    \|C_\varphi\|_{B_\sigma^p\to B_\tau^p}
    &=
    |a|^{-\frac{1}{p}}
    \|T_{\operatorname{Im}b}\|_{B_\sigma^p\to B_\sigma^p}
=
    |a|^{-\frac{1}{p}}
    e^{\sigma|\operatorname{Im}b|}.
\end{aligned}
\]
This proves the theorem.
\end{proof}

\begin{remark}
The preceding composition operators also show that the lower estimate
in Theorem~\ref{thm-wco-p-le-q-bounded} is sharp in general. Indeed,
fix \(\sigma>0\) and \(0<p\leq\infty\), and consider the vertical
translations
\[
    \varphi_\beta(z):=z+i\beta,
    \qquad \beta>0.
\]
By Theorem~\ref{thm-composition-norm},
\[
    \|C_{\varphi_\beta}\|_{B_\sigma^p\to B_\sigma^p}
    =
    e^{\sigma\beta}.
\]
On the other hand,
\[
    M_{\sigma,\sigma;p,p}(1,\varphi_\beta)
    =
    \sup_{y\in\R}
    \frac{\Phi_{\sigma,p}(|y+\beta|)}
         {\Phi_{\sigma,p}(|y|)}  \geq \lim_{y\to+\infty}
\frac{\Phi_{\sigma,p}(y+\beta)}
     {\Phi_{\sigma,p}(y)}
=
e^{\sigma\beta},
\]
while Proposition~\ref{prop-wco-structural} gives the reverse estimate,
up to a constant depending only on \(p\). Hence
\[
    M_{\sigma,\sigma;p,p}(1,\varphi_\beta)
    \simeq
    e^{\sigma\beta}
    =
    \|C_{\varphi_\beta}\|_{B_\sigma^p\to B_\sigma^p}.
\]
Thus the lower estimate in terms of
\(M_{\sigma,\tau;p,q}(\psi,\varphi)\) is sharp already for unweighted
composition operators.

\smallskip
\noindent
Moreover, on the diagonal \(p=q\), for
\[
    \varphi(z)=az+b,
    \qquad
    a\in\R,\quad 0<|a|\leq1,\quad b\in\C,
\]
the exact formula
\[
    \|C_\varphi\|_{B_\sigma^p\to B_\sigma^p}
    =
    |a|^{-\frac{1}{p}}e^{\sigma|\operatorname{Im}b|}
\]
shows that the dependence on \(|a|\) and
\(\operatorname{Im}b\) in the corresponding upper estimate cannot be
improved in general.
\end{remark}

\subsection{Multiplication operators}
\label{subsec-multiplication-operators}

For an entire function \(\psi\), let
\[
    M_\psi f:=\psi f.
\]
Since
\[
    M_\psi=W_{\psi,\operatorname{id}},
\]
the preceding results immediately yield a complete description of
multipliers
\[
    M_\psi:B_\sigma^p\longrightarrow B_\tau^q.
\]
Here the two parts of the weighted-composition theory separate cleanly:
the difference \(\tau-\sigma\) is the bandwidth available to the
multiplier, while the relation between \(p\) and \(q\) determines the
integrability required of its restriction to \(\R\).

\begin{theorem}[Multiplication operators]
\label{thm-multiplication-operators}
Let \(\sigma,\tau>0\), \(0<p,q\leq\infty\), and let
\(\psi\not\equiv0\) be entire.

\begin{enumerate}
\item[\textnormal{(a)}]
Let \(0<p\leq q\leq\infty\). Then
\(
    M_\psi:
    B_\sigma^p\longrightarrow B_\tau^q
\)
is bounded if and only if
\[
    \sigma\leq\tau
    \qquad\text{and}\qquad
    \psi\in B_{\tau-\sigma}^\infty.
\]
In this case,
\[
    \sigma^{\frac{1}{p}}\tau^{-\frac{1}{q}}\|\psi\|_{L^\infty}
    \lesssim
    \|M_\psi\|_{B_\sigma^p\to B_\tau^q}
    \lesssim
    \sigma^{\frac{1}{p}-\frac{1}{q}}\|\psi\|_{L^\infty},
\]
where the implicit constants depend only on \(p\) and \(q\).

Moreover,
\(
    M_\psi:
    B_\sigma^p\longrightarrow B_\tau^q
\)
is compact if and only if
\[
    \sigma<\tau,
    \qquad
    \psi\in B_{\tau-\sigma}^\infty,
\qquad \text{and} \qquad
    \psi(x)\longrightarrow0
    \quad \text{as} \quad |x|\to\infty.
\]

\item[\textnormal{(b)}]
Let \(0<q<p\leq\infty\).
Then boundedness and compactness of
\(
    M_\psi:
    B_\sigma^p\longrightarrow B_\tau^q
\)
are equivalent. These conditions hold if and only if
\[
    \sigma<\tau
    \qquad\text{and}\qquad
    \psi\in B_{\tau-\sigma}^{\frac{pq}{p-q}}.
\]
Moreover,
\[
    \left(\frac{\sigma}{\tau}\right)^{\frac1p}
    \|\psi\|_{L^{\frac{pq}{p-q}}}
    \lesssim
    \|M_\psi\|_{B_\sigma^p\to B_\tau^q}
    \leq
    \|\psi\|_{L^{\frac{pq}{p-q}}},
\]
where the implicit constant depends only on \(p\) and \(q\).

If \(p=\infty\), then
\[
    \|M_\psi\|_{B_\sigma^\infty\to B_\tau^q}
    =
    \|\psi\|_{L^q}.
\]
\end{enumerate}
\end{theorem}

\begin{proof}
Since $M_\psi=W_{\psi,\operatorname{id}}$, the boundedness and compactness characterizations follow immediately
from Theorems~\ref{thm-wco-p-le-q-bounded},
\ref{thm-wco-p-le-q-compact}, and
\ref{thm-wco-q-lt-p}, applied with 
\[
    \varphi(z)=z,
    \qquad a=1,
    \qquad b=0.
\]
It remains to prove the norm estimates.

\smallskip
\noindent
\emph{The case \(0<p\leq q\leq\infty\).}
Suppose that \(M_\psi:B_\sigma^p\to B_\tau^q\) is bounded.
For every \(x\in\R\),
\[
    \delta_x\circ M_\psi
    =
    \psi(x)\delta_x,
\]
where the evaluation on the left is taken on \(B_\tau^q\), while the
one on the right is taken on \(B_\sigma^p\). Hence
\[
    |\psi(x)|
    \|\delta_x\|_{(B_\sigma^p)^*}
    \leq
    \|\delta_x\|_{(B_\tau^q)^*}
    \|M_\psi\|_{B_\sigma^p\to B_\tau^q}.
\]
By Theorem~\ref{thm:point-evaluations},
\[
    \|\delta_x\|_{(B_\sigma^p)^*}
    \simeq_p
    \sigma^{\frac{1}{p}},
    \qquad
    \|\delta_x\|_{(B_\tau^q)^*}
    \simeq_q
    \tau^{\frac{1}{q}},
    \qquad x\in\R.
\]
Therefore
\[
    \sigma^{\frac{1}{p}}\tau^{-\frac{1}{q}}|\psi(x)|
    \lesssim_{p,q}
    \|M_\psi\|_{B_\sigma^p\to B_\tau^q}.
\]
Taking the supremum over \(x\in\R\), we obtain
\[
    \sigma^{\frac{1}{p}}\tau^{-\frac{1}{q}}
    \|\psi\|_{L^\infty}
    \lesssim_{p,q}
    \|M_\psi\|_{B_\sigma^p\to B_\tau^q}.
\]
Conversely, if
\[
    \sigma\leq\tau,
    \qquad
    \psi\in B_{\tau-\sigma}^\infty,
\]
then, for \(f\in B_\sigma^p\), the product \(\psi f\) has exponential
type at most \(\tau\). Moreover, Nikol'skii's inequality gives
\[
\begin{aligned}
    \|M_\psi f\|_{L^q}
    &=
    \|\psi f\|_{L^q}
\leq
    \|\psi\|_{L^\infty}\|f\|_{L^q}
\lesssim_{p,q}
    \sigma^{\frac{1}{p}-\frac{1}{q}}
    \|\psi\|_{L^\infty}\|f\|_{L^p}.
\end{aligned}
\]
Thus
\[
    \|M_\psi\|_{B_\sigma^p\to B_\tau^q}
    \lesssim_{p,q}
    \sigma^{\frac{1}{p}-\frac{1}{q}}\|\psi\|_{L^\infty}.
\]
This proves the two-sided estimate in part \textnormal{(a)}.

\smallskip
\noindent
\emph{The case \(0<q<p\leq\infty\).}
Assume first that \(p<\infty\). Since
\[
    \frac1q
    =
    \frac1p+\frac{p-q}{pq},
\]
H\"older's inequality gives, for \(f\in B_\sigma^p\),
\[
\begin{aligned}
    \|M_\psi f\|_{L^q}
    &=
    \|\psi f\|_{L^q}
\leq
    \|\psi\|_{L^{\frac{pq}{p-q}}}
    \|f\|_{L^p}.
\end{aligned}
\]
Hence
\[
    \|M_\psi\|_{B_\sigma^p\to B_\tau^q}
    \leq
    \|\psi\|_{L^{\frac{pq}{p-q}}}.
\]
For the reverse estimate, Theorem~\ref{thm-wco-q-lt-p}, applied with
\(a=1\) and \(b=0\), gives
\[
    \left(\frac{\sigma}{\tau}\right)^{\frac{1}{p}}
    \|\psi\|_{L^{\frac{pq}{p-q}}}
    \lesssim_{p,q}
    \|M_\psi\|_{B_\sigma^p\to B_\tau^q}.
\]
This proves
\[
    \left(\frac{\sigma}{\tau}\right)^{\frac{1}{p}}
    \|\psi\|_{L^{\frac{pq}{p-q}}}
    \lesssim_{p,q}
    \|M_\psi\|_{B_\sigma^p\to B_\tau^q}
    \leq
    \|\psi\|_{L^{\frac{pq}{p-q}}}.
\]
Finally, let \(p=\infty\). For \(f\in B_\sigma^\infty\),
\[
    \|M_\psi f\|_{L^q}
    \leq
    \|\psi\|_{L^q}\|f\|_{L^\infty},
\]
and therefore
\[
    \|M_\psi\|_{B_\sigma^\infty\to B_\tau^q}
    \leq
    \|\psi\|_{L^q}.
\]
Conversely, the constant function \(1\) belongs to
\(B_\sigma^\infty\) and satisfies $\|1\|_{L^\infty}=1.$ Thus
\[
    \|M_\psi\|_{B_\sigma^\infty\to B_\tau^q}
    \geq
    \|M_\psi1\|_{L^q}
    =
    \|\psi\|_{L^q}.
\]
Consequently,
\[
    \|M_\psi\|_{B_\sigma^\infty\to B_\tau^q}
    =
    \|\psi\|_{L^q}.
\]
The proof is complete.
\end{proof}

\begin{remark}
Theorem~\ref{thm-multiplication-operators} may be summarized as
\[
\operatorname{Mult}\bigl(B_\sigma^p,B_\tau^q\bigr)
=
\begin{cases}
B_{\tau-\sigma}^{\infty},
    & 0<p\le q\le\infty,\ \sigma\le\tau,\\[2mm]
B_{\tau-\sigma}^{\frac{pq}{p-q}},
    & 0<q<p\le\infty,\ \sigma<\tau,\\[2mm]
\{0\},
    & \text{otherwise}.
\end{cases}
\]
In particular, when \(\sigma=\tau\), every bounded multiplier in the
range \(p\le q\) is constant, whereas for \(q<p\) there are no nonzero
bounded multipliers.
\end{remark}

\section{Differences of weighted composition operators}
\label{sec:differences}

We now study boundedness and compactness of differences of weighted
composition operators
\[
    W_{\psi,\varphi}-W_{\vartheta,\phi}
    :B_\sigma^p\longrightarrow B_\tau^q,
    \qquad
    \sigma,\tau>0,\quad 0<p,q\le\infty.
\]
Throughout this section, we write
\begin{equation*}
\label{eq:difference-notation}
    L:=W_{\psi,\varphi}-W_{\vartheta,\phi},
    \qquad
    \Delta:=\varphi-\phi.
\end{equation*}
When the two symbols coincide, the difference reduces to a single
weighted composition operator,
\[
    W_{\psi,\varphi}-W_{\vartheta,\varphi}
    =W_{\psi-\vartheta,\varphi},
\]
and the results of the preceding section apply directly. Thus the main
issue is to understand what boundedness or compactness of \(L\) forces
when the two symbols are distinct. As in the preceding section, the
analysis separates naturally into the ranges
\(0<p\le q\le\infty\) and \(0<q<p\le\infty\).

\subsection{Structural reduction}

The quantitative separation of point evaluations obtained in
Lemma~\ref{lem-two-point-separation} will be used repeatedly below.
We begin with the structural consequences of boundedness.

\begin{proposition}\label{prop-difference-structural}
Let \(\sigma,\tau>0\) and \(0<p,q\leq\infty\). Let
\(\psi,\vartheta,\varphi,\phi\) be entire functions such that
\[
    \psi\not\equiv0,\qquad
    \vartheta\not\equiv0,\qquad
    \varphi\not\equiv\phi.
\]
Suppose that the difference
\(
    L:B_\sigma^p\longrightarrow B_\tau^q
\)
is bounded. Then
\[
    \varphi(z)=az+b,\qquad
    \phi(z)=cz+d,
\]
where
\[
    a,c\in\R,\qquad
    \sigma|a|\leq\tau,\qquad
    \sigma|c|\leq\tau,
    \qquad b,d\in\C.
\]
Moreover, 
\begin{equation}\label{eq-difference-weight-separation-psi}
    |\psi(z)|
    \Phi_{\sigma,p}(|\im\varphi(z)|)
    \frac{\sigma|\Delta(z)|}
    {1+\sigma|\im\varphi(z)|+\sigma|\Delta(z)|}
    \lesssim
    \Phi_{\tau,q}(|\im z|)\|L\|,
    \qquad z\in\C,
\end{equation}
and
\begin{equation}\label{eq-difference-weight-separation-vartheta}
    |\vartheta(z)|
    \Phi_{\sigma,p}(|\im\phi(z)|)
    \frac{\sigma|\Delta(z)|}
    {1+\sigma|\im\phi(z)|+\sigma|\Delta(z)|}
    \lesssim
    \Phi_{\tau,q}(|\im z|)\|L\|,
    \qquad z\in\C,
\end{equation}
where the implicit constants depend only on \(p\) and \(q\).
In addition, \(\psi\) and \(\vartheta\) have exponential type at most $\tau-\sigma|a|$ and $\tau-\sigma|c|$, respectively.
\end{proposition}

\begin{proof}
We divide the proof into four steps.

\smallskip
\noindent
\emph{Step 1. The symbols are affine.}
Choose \(\varepsilon>0\) such that $3\varepsilon<\frac{\sigma}{2}$ and $e^{-i\varepsilon\varphi}
    \not\equiv
    e^{-i\varepsilon\phi}$. Such a choice is possible because \(\varphi-\phi\not\equiv0\). Indeed, if
\(\varphi-\phi\) is nonconstant, the two exponentials are distinct for every
\(\varepsilon>0\), while if \(\varphi-\phi\) is a nonzero constant, it is
enough to avoid a discrete set of values of \(\varepsilon\).

\smallskip
\noindent
Choose an even integer \(N>\frac{1}{p}\) and
a point \(\zeta\in\C\) such that the function
\[
    h(z)
    :=
    \left(
        \frac{\sin\bigl(\sigma(z-\zeta)/(2N)\bigr)}
        {\sigma(z-\zeta)/(2N)}
    \right)^N
\]
satisfies $h\circ\varphi\not\equiv0$ and $h\circ\phi\not\equiv0.$ This is always possible: for a nonconstant symbol the corresponding
composition cannot vanish identically, while for a constant symbol one
chooses \(\zeta\) so that \(h\) does not vanish at its constant value.
Moreover, $h\in B_{\sigma/2}^p$.

\smallskip
\noindent
For \(j=0,1,2,3\), set
\[
    h_j(z):=e^{-ij\varepsilon z}h(z).
\]
Then
\[
    \operatorname{type}(h_j)
    \leq \frac{\sigma}{2}+j\varepsilon<\sigma,
\qquad \text{and} \qquad
    |h_j(x)|=|h(x)|,\qquad x\in\R.
\]
Hence \(h_j\in B_\sigma^p\). Since \(L\) is bounded,
\[
    A_j:=Lh_j\in B_\tau^q,
    \qquad j=0,1,2,3.
\]
Put
\[
    U:=\psi(h\circ\varphi),\qquad
    V:=\vartheta(h\circ\phi),
\qquad \text{and} \qquad
    X:=e^{-i\varepsilon\varphi},\qquad
    Y:=e^{-i\varepsilon\phi}.
\]
Then
\[
    A_j=UX^j-VY^j,
    \qquad j=0,1,2,3.
\]
Define
\[
    P_0:=A_0A_2-A_1^2,\qquad
    P_1:=A_1A_3-A_2^2,\qquad
    P_2:=A_0A_3-A_1A_2.
\]
A direct calculation gives
\begin{equation}\label{eq-difference-algebra-P}
    P_0=-UV(X-Y)^2, \qquad P_1=XY\,P_0, \qquad P_2=(X+Y)P_0.
\end{equation}
By the choice of \(h\) and \(\varepsilon\),
\[
    U\not\equiv0,\qquad
    V\not\equiv0,\qquad
    X-Y\not\equiv0,
\]
and hence
\[
    P_0\not\equiv0.
\]
Since each \(A_j\) has exponential type at most \(\tau\), the functions
\(P_0,P_1,P_2\) have finite exponential type.

\smallskip
\noindent
From \eqref{eq-difference-algebra-P} and
Lemma~\ref{lem-product-rigidity}, it follows that
\[
    XY=e^{-i\varepsilon(\varphi+\phi)}
\]
has finite exponential type. Similarly, if \(P_2\not\equiv0\), then
\eqref{eq-difference-algebra-P} and Lemma~\ref{lem-product-rigidity} show
that \(X+Y\) has finite exponential type; if \(P_2\equiv0\), then
\(X+Y\equiv0\), so the same conclusion is immediate. Therefore
\[
    (X-Y)^2=(X+Y)^2-4XY
\]
has finite exponential type, and consequently \(X-Y\) has finite
exponential type. Hence
\[
    X=\frac{(X+Y)+(X-Y)}2,
    \qquad
    Y=\frac{(X+Y)-(X-Y)}2
\]
are both of finite exponential type.

\smallskip
\noindent
Since \(X\) and \(Y\) are zero-free, the Hadamard factorization theorem
gives
\[
    X(z)=e^{P(z)},\qquad
    Y(z)=e^{Q(z)},
\]
where \(P\) and \(Q\) are polynomials of degree at most one. Thus
\[
    -i\varepsilon\varphi(z)-P(z)\in2\pi i\Z,
    \qquad
    -i\varepsilon\phi(z)-Q(z)\in2\pi i\Z.
\]
Both left-hand sides are entire and take values in a discrete set, so they
are constant. Consequently,
\[
    \varphi(z)=az+b,\qquad
    \phi(z)=cz+d
\]
for some \(a,b,c,d\in\C\).

\smallskip
\noindent
\emph{Step 2. Pointwise separation of the weights.}
For every \(z\in\C\),
\[
    \delta_z\circ L
    =
    \psi(z)\delta_{\varphi(z)}
    -
    \vartheta(z)\delta_{\phi(z)}.
\]
Hence, by Theorem~\ref{thm:point-evaluations},
\[
\begin{aligned}
    \bigl\|
        \psi(z)\delta_{\varphi(z)}
        -
        \vartheta(z)\delta_{\phi(z)}
    \bigr\|_{(B_\sigma^p)^*}
    \leq
    \|\delta_z\|_{(B_\tau^q)^*}\|L\|
    \lesssim_q
    \Phi_{\tau,q}(|\im z|)\|L\|.
\end{aligned}
\]
The distance from \(\psi(z)\delta_{\varphi(z)}\) to
\(\C\delta_{\phi(z)}\) is bounded above by the left-hand side. Therefore
Lemma~\ref{lem-two-point-separation} gives
\[
    |\psi(z)|
    \Phi_{\sigma,p}(|\im\varphi(z)|)
    \frac{\sigma|\varphi(z)-\phi(z)|}
    {1+\sigma|\im\varphi(z)|
      +\sigma|\varphi(z)-\phi(z)|}
    \lesssim_{p,q}
    \Phi_{\tau,q}(|\im z|)\|L\|,
\]
which is \eqref{eq-difference-weight-separation-psi}. Interchanging the two
terms gives \eqref{eq-difference-weight-separation-vartheta}.

\smallskip
\noindent
\emph{Step 3. The coefficients \(a\) and \(c\) are real.}
Set
\[
    F:=\psi\Delta,
    \qquad
    G:=\vartheta\Delta.
\]
Since neither weight nor \(\Delta\) is identically zero,
\[
    F\not\equiv0,
    \qquad
    G\not\equiv0.
\]
From \eqref{eq-difference-weight-separation-psi},
\[
    |F(z)|
    \lesssim_{p,q}
    \frac{\Phi_{\tau,q}(|\im z|)}
         {\Phi_{\sigma,p}(|\im\varphi(z)|)}
    \frac{
        1+\sigma|\im\varphi(z)|+\sigma|\Delta(z)|
    }{\sigma}
    \|L\|.
\]
Since \(\varphi\) and \(\Delta\) are affine, the right-hand side is bounded
by
\[
    C(1+|z|)^M e^{\tau|\im z|}
\]
for suitable \(C,M>0\). Thus \(F\) has finite exponential type. The same
argument shows that \(G\) has finite exponential type.

\smallskip
\noindent
Write
\[
    a=a_1+ia_2,\qquad
    b=b_1+ib_2,
\]
with \(a_j,b_j\in\R\). Suppose that \(a_2\neq0\). Taking \(z=x\in\R\) in
the preceding estimate gives
\[
    |F(x)|
    \lesssim
    (1+|x|)^M
    e^{-\sigma|a_2x+b_2|},
    \qquad x\in\R.
\]
Since \(a_2\neq0\), the polynomial factor can be absorbed into the
exponential. Hence there exist \(C_0,\gamma_0>0\) such that
\[
    |F(x)|\leq C_0e^{-\gamma_0|x|},
    \qquad x\in\R.
\]
Lemma~\ref{lem-exponential-decay-rigidity} gives \(F\equiv0\), a
contradiction. Therefore \(a_2=0\), and hence \(a\in\R\).
The same argument applied to \(G\) shows that \(c\in\R\).

\smallskip
\noindent
\emph{Step 4. The bandwidth restrictions and the types of the weights.}
Suppose that $\sigma|a|>\tau.$ Since \(a\in\R\), taking \(z=iy\) in the estimate for \(F\) gives
\[
    |F(iy)|
    \lesssim
    (1+|y|)^M
    e^{\tau|y|-\sigma|ay+\im b|}\lesssim
    (1+|y|)^M
    e^{-(\sigma|a|-\tau)|y|}.
\]
After absorbing the polynomial factor, there exist \(C_1,\gamma_1>0\)
such that
\[
    |F(iy)|\leq C_1e^{-\gamma_1|y|},
    \qquad y\in\R.
\]
The entire function \(H(z):=F(iz)\) has finite exponential type, so
Lemma~\ref{lem-exponential-decay-rigidity} yields \(H\equiv0\). Hence
\(F\equiv0\), a contradiction. Therefore $\sigma|a|\leq\tau.$ Likewise, $\sigma|c|\leq\tau.$

\smallskip
\noindent
It remains to record the exponential types of the weights. Since
\(\Delta\) is a nonzero affine function, there exist \(R>0\) and \(M>0\)
such that
\[
    \frac{
        1+\sigma|\im\varphi(z)|+\sigma|\Delta(z)|
    }
    {\sigma|\Delta(z)|}
    \lesssim
    (1+|z|)^M,
    \qquad |z|\geq R.
\]
Hence \eqref{eq-difference-weight-separation-psi} gives
\[
    |\psi(z)|
    \lesssim
    (1+|z|)^M
    \frac{\Phi_{\tau,q}(|\im z|)}
         {\Phi_{\sigma,p}(|\im\varphi(z)|)},
    \qquad |z|\geq R.
\]
Writing \(z=x+iy\) and using $\im\varphi(z)=ay+\im b,$ we obtain
\[
    |\psi(z)|
    \lesssim
    (1+|z|)^{M+1/p}
    e^{\tau|y|-\sigma|ay+\im b|}
    \lesssim
    (1+|z|)^{M+1/p}
    e^{(\tau-\sigma|a|)|y|}.
\]
It follows that \(\psi\) has exponential type at most $\tau-\sigma|a|.$ Applying the same argument to
\eqref{eq-difference-weight-separation-vartheta} shows that \(\vartheta\) has exponential type at most $\tau-\sigma|c|.$ The proof is complete.
\end{proof}

\subsection{The case \(0<p\leq q\leq\infty\)}

In this range, boundedness of a difference with distinct symbols cannot
result from cancellation between two unbounded weighted composition
operators.

\begin{theorem}\label{thm-difference-bounded-p-le-q}
Let \(\sigma,\tau>0\) and \(0<p\leq q\leq\infty\). Let
\(\psi,\vartheta,\varphi,\phi\) be entire functions. Then the difference
\[
    W_{\psi,\varphi}-W_{\vartheta,\phi}
    :B_\sigma^p\longrightarrow B_\tau^q
\]
is bounded if and only if one of the following conditions holds:
\begin{enumerate}
\item[\textnormal{(a)}] \(\varphi\not\equiv\phi\) and both
\(
    W_{\psi,\varphi},\,
    W_{\vartheta,\phi}
    :B_\sigma^p\longrightarrow B_\tau^q
\)
are bounded;

\item[\textnormal{(b)}] \(\varphi=\phi\) and
\(
    W_{\psi-\vartheta,\varphi}
    :B_\sigma^p\longrightarrow B_\tau^q
\)
is bounded.
\end{enumerate}
\end{theorem}

\begin{proof}
The sufficiency is immediate. Indeed, in the first case the difference of
two bounded operators is bounded, while in the second case
\[
    W_{\psi,\varphi}-W_{\vartheta,\phi}
    =
    W_{\psi-\vartheta,\varphi}.
\]
Conversely, 
suppose that the difference
\(
    L:B_\sigma^p\longrightarrow B_\tau^q
\)
is bounded. If \(\varphi=\phi\), then $L=W_{\psi-\vartheta,\varphi},$ so the second condition holds. We may therefore assume that
\[
    \varphi\not\equiv\phi.
\]
If \(\psi\equiv0\) or \(\vartheta\equiv0\), the conclusion is immediate.
Thus suppose that $\psi\not\equiv0$ and $\vartheta\not\equiv0.$

\smallskip
\noindent
By Proposition~\ref{prop-difference-structural},
\[
    \varphi(z)=az+b,\qquad
    \phi(z)=cz+d,
\]
where
\[
    a,c\in\R,\qquad
    \sigma|a|\leq\tau,\qquad
    \sigma|c|\leq\tau.
\]
Moreover,
\[
    \operatorname{type}(\psi)\leq\tau-\sigma|a|,
    \qquad
    \operatorname{type}(\vartheta)\leq\tau-\sigma|c|.
\]
We first show that both weights are bounded on \(\R\). Since \(a,c\in\R\), for \(x\in\R\),
\[
    \im\varphi(x)=\im b,
    \qquad
    \im\phi(x)=\im d.
\]
If \(a\neq c\), then
\[
    |\Delta(x)|\longrightarrow\infty
    \qquad (|x|\to\infty).
\]
Hence there exist \(R>0\) and \(\kappa>0\) such that
\[
    \frac{\sigma|\Delta(x)|}
    {1+\sigma|\im b|+\sigma|\Delta(x)|}
    \geq\kappa,
    \qquad |x|\geq R,
\]
and likewise
\[
    \frac{\sigma|\Delta(x)|}
    {1+\sigma|\im d|+\sigma|\Delta(x)|}
    \geq\kappa,
    \qquad |x|\geq R.
\]
If \(a=c\), then \(\Delta=b-d\neq0\), and the same estimates hold for all
\(x\in\R\), after possibly decreasing \(\kappa\).

\smallskip
\noindent
Taking \(z=x\in\R\) in
\eqref{eq-difference-weight-separation-psi} and
\eqref{eq-difference-weight-separation-vartheta}, we therefore obtain
\[
    \sup_{|x|\geq R}|\psi(x)|<\infty,
    \qquad
    \sup_{|x|\geq R}|\vartheta(x)|<\infty.
\]
Since the weights are continuous, $\psi,\vartheta\in L^\infty(\R).$ Consequently,
\begin{equation}\label{eq-difference-bounded-weight-classes}
    \psi\in B^\infty_{\tau-\sigma|a|},
    \qquad
    \vartheta\in B^\infty_{\tau-\sigma|c|}.
\end{equation}
We now distinguish three cases.

\smallskip
\noindent
Suppose first that \(a\neq0\) and \(c\neq0\). Then
\[
    0<\sigma|a|\leq\tau,
    \qquad
    0<\sigma|c|\leq\tau,
\]
and \eqref{eq-difference-bounded-weight-classes}, together with
Theorem~\ref{thm-wco-p-le-q-bounded}, shows that both
\[
    W_{\psi,\varphi},\,
    W_{\vartheta,\phi}
    :B_\sigma^p\longrightarrow B_\tau^q
\]
are bounded.

\smallskip
\noindent
Next suppose that exactly one of \(a,c\) is zero. By symmetry, assume that
\[
    a=0,\qquad c\neq0.
\]
Then \(\varphi\equiv b\), whereas \(\phi\) is nonconstant. By
\eqref{eq-difference-bounded-weight-classes} and
Theorem~\ref{thm-wco-p-le-q-bounded},
\(
    W_{\vartheta,\phi}
    :B_\sigma^p\longrightarrow B_\tau^q
\)
is bounded. Hence
\[
    W_{\psi,\varphi}
    =
    L+W_{\vartheta,\phi}
\]
is bounded as well.

\smallskip
\noindent
It remains to consider \(a=c=0\). Then
\[
    \varphi\equiv b,\qquad
    \phi\equiv d,
\]
and \(b\neq d\). By Lemma~\ref{lem-two-point-separation}, the evaluation
functionals \(\delta_b\) and \(\delta_d\) are linearly independent.
Therefore there exist \(f,g\in B_\sigma^p\) such that
\[
    f(d)=0,\qquad f(b)\neq0,
\qquad \text{and} \qquad
    g(b)=0,\qquad g(d)\neq0.
\]
It follows that
\[
    Lf=f(b)\psi\in B_\tau^q,
    \qquad
    Lg=-g(d)\vartheta\in B_\tau^q.
\]
Hence $\psi,\vartheta\in B_\tau^q.$ Proposition~\ref{prop-wco-constant-symbol} now shows that both
\[
    W_{\psi,\varphi},\,
    W_{\vartheta,\phi}
    :B_\sigma^p\longrightarrow B_\tau^q
\]
are bounded.

\smallskip
\noindent
Thus the first condition holds whenever
\(\varphi\not\equiv\phi\), and the proof is complete.
\end{proof}

\begin{theorem}\label{thm-difference-compact-p-le-q}
Let \(\sigma,\tau>0\) and \(0<p\leq q\leq\infty\). Let
\(\psi,\vartheta,\varphi,\phi\) be entire functions. Then the difference
\[
    W_{\psi,\varphi}-W_{\vartheta,\phi}
    :B_\sigma^p\longrightarrow B_\tau^q
\]
is compact if and only if one of the following conditions holds:
\begin{enumerate}
\item[\textnormal{(a)}] \(\varphi\not\equiv\phi\) and both
\(
    W_{\psi,\varphi},\,
    W_{\vartheta,\phi}
    :B_\sigma^p\longrightarrow B_\tau^q
\)
are compact;

\item[\textnormal{(b)}] \(\varphi=\phi\) and
\(
    W_{\psi-\vartheta,\varphi}
    :B_\sigma^p\longrightarrow B_\tau^q
\)
is compact.
\end{enumerate}
\end{theorem}

\begin{proof}
The sufficiency is immediate. Indeed, in the first case the difference of
two compact operators is compact, while in the second case
\[
    W_{\psi,\varphi}-W_{\vartheta,\phi}
    =
    W_{\psi-\vartheta,\varphi}.
\]
Conversely, suppose that the difference
\(
    L:B_\sigma^p\longrightarrow B_\tau^q
\)
is compact. If \(\varphi=\phi\), then $L=W_{\psi-\vartheta,\varphi},$ so the second condition holds. We may therefore assume that $\varphi\not\equiv\phi.$

\smallskip
\noindent
Since \(L\) is bounded, Theorem~\ref{thm-difference-bounded-p-le-q}
shows that both
\[
    W_{\psi,\varphi},\,
    W_{\vartheta,\phi}
    :B_\sigma^p\longrightarrow B_\tau^q
\]
are bounded. If one of them is compact, then so is the other, because
\[
    W_{\vartheta,\phi}=W_{\psi,\varphi}-L
\qquad \text{and} \qquad
    W_{\psi,\varphi}=L+W_{\vartheta,\phi}.
\]
It therefore remains to show that the two operators cannot both be
noncompact.

\smallskip
\noindent
Suppose, to the contrary, that neither operator is compact. A bounded
weighted composition operator with constant symbol is compact by
Proposition~\ref{prop-wco-constant-symbol}. Hence both symbols are
nonconstant. The uniqueness hypothesis in Lemma~\ref{lem-compactness-criterion}
is automatic for \(L\). Indeed, \(f_j\to0\) locally uniformly implies
\(Lf_j\to0\) locally uniformly, while convergence in \(B_\tau^q\)
implies locally uniform convergence by
Theorem~\ref{thm:point-evaluations}. By Theorem~\ref{thm-difference-bounded-p-le-q},
\[
    \varphi(z)=az+b,\qquad
    \phi(z)=cz+d,
\]
where
\[
    a,c\in\R,\qquad
    0<\sigma|a|\leq\tau,\qquad
    0<\sigma|c|\leq\tau,
\qquad \text{and} \qquad
    \psi\in B^\infty_{\tau-\sigma|a|},
    \qquad
    \vartheta\in B^\infty_{\tau-\sigma|c|}.
\]
In particular, $\psi,\vartheta\in L^\infty(\R).$

\smallskip
\noindent
Since \(W_{\psi,\varphi}\) is not compact, the compactness
characterization in Theorem~\ref{thm-wco-p-le-q-compact} gives
\[
    \limsup_{|x|\to\infty}|\psi(x)|>0.
\]
Indeed, this follows directly if \(\sigma|a|<\tau\), while if
\(\sigma|a|=\tau\), then $\psi\in B^\infty_0=\C,$ and \(\psi\) is a nonzero constant. Hence there exist
\(\varepsilon>0\) and a sequence \((x_n)\subset\R\) such that
\[
    |x_n|\longrightarrow\infty,
    \qquad
    |\psi(x_n)|\geq\varepsilon,
    \qquad n\in\N.
\]
Put
\[
    w_n:=\varphi(x_n)=ax_n+b,
    \qquad
    v_n:=\phi(x_n)=cx_n+d.
\]
We first show that \(a=c\). Suppose that \(a\neq c\). Choose
\(m\in\N\) such that $m>\frac1p$ and define
\[
    f_n
    :=
    \frac{G_{\sigma,m,w_n}}
         {\|G_{\sigma,m,w_n}\|_{L^p}}.
\]
Then
\[
    \|f_n\|_{L^p}=1.
\]
Since \(a\neq0\),
\[
    |\re w_n|\longrightarrow\infty,
    \qquad
    \im w_n=\im b,
\]
and Lemma~\ref{lem:packet-decay} gives
\[
    f_n\longrightarrow0
\]
uniformly on compact subsets of \(\C\).  Since \(L\) is compact, the
compactness criterion in Lemma~\ref{lem-compactness-criterion} yields
\begin{equation}\label{eq-compact-difference-Lfn-zero}
    \|Lf_n\|_{L^q}\longrightarrow0.
\end{equation}
By the point-evaluation estimate on the real axis,
\[
    |Lf_n(x_n)|
    \lesssim_q
    \tau^{1/q}\|Lf_n\|_{L^q}
    \longrightarrow0.
\]
On the other hand,
\[
    Lf_n(x_n)
    =
    \psi(x_n)f_n(w_n)
    -
    \vartheta(x_n)f_n(v_n).
\]
By Lemma~\ref{lem:packet-estimates},
\[
    |f_n(w_n)|
    \simeq_p
    \Phi_{\sigma,p}(|\im b|),
\]
and hence these quantities are bounded below by a positive constant
independent of \(n\). Moreover,
\(
    v_n-w_n
    =
    (c-a)x_n+d-b,
\)
so
\[
    |v_n-w_n|\longrightarrow\infty.
\]
Since \(\im w_n=\im b\), \(\im v_n=\im d\), and
\(|\re v_n-\re w_n|\longrightarrow\infty\), the estimate used in the
proof of Lemma~\ref{lem:packet-decay} gives
\[
    f_n(v_n)\longrightarrow0.
\]
As \(\vartheta\in L^\infty(\R)\),
\[
    \vartheta(x_n)f_n(v_n)\longrightarrow0.
\]
Consequently,
\[
    |Lf_n(x_n)|
    \geq
    \varepsilon |f_n(w_n)|-o(1),
\]
which contradicts \eqref{eq-compact-difference-Lfn-zero}. Therefore $a=c.$

\smallskip
\noindent
Since \(\varphi\not\equiv\phi\), we must have \(b\neq d\). We now derive
a second contradiction. Choose \(m\in\N\) such that $m>1+\frac1p,$ and set
\[
    \zeta_b:=i\,\im b,
    \qquad
    \zeta_d:=\re d-\re b+i\,\im d.
\]
Then \(\zeta_b\neq\zeta_d\). Define
\[
    F(z)
    :=
    (z-\zeta_d)G_{\sigma,m,\zeta_b}(z).
\]
As in the proof of Lemma~\ref{lem-two-point-separation}, $F\in B_\sigma^p,$ and
\[
    F(\zeta_d)=0,
    \qquad
    F(\zeta_b)
    =
    (\zeta_b-\zeta_d)G_{\sigma,m,\zeta_b}(\zeta_b)
    \neq0.
\]
Put
\[
    t_n:=ax_n+\re b
\]
and define
\[
    g_n(z)
    :=
    \frac{F(z-t_n)}{\|F\|_{L^p}}.
\]
Since \(t_n\in\R\), the norm $\|g_n\|_{L^p}=1.$ Moreover, \(|t_n|\to\infty\), and the decay estimate used in the proof of
Lemma~\ref{lem-two-point-separation} shows that
\[
    g_n\longrightarrow0
\]
uniformly on compact subsets of \(\C\). Compactness of \(L\) therefore
gives
\[
    \|Lg_n\|_{L^q}\longrightarrow0,
\]
and hence
\begin{equation}\label{eq-compact-difference-Lgn-point}
    |Lg_n(x_n)|\longrightarrow0.
\end{equation}
Since \(a=c\),
\[
    \varphi(x_n)-t_n
    =
    i\,\im b
    =
    \zeta_b,
\]
while
\[
    \phi(x_n)-t_n
    =
    \re d-\re b+i\,\im d
    =
    \zeta_d.
\]
Thus
\[
\begin{aligned}
    Lg_n(x_n)
    =
    \psi(x_n)g_n(\varphi(x_n))
    -
    \vartheta(x_n)g_n(\phi(x_n))
    =
    \psi(x_n)
    \frac{F(\zeta_b)}{\|F\|_{L^p}},
\end{aligned}
\]
because \(F(\zeta_d)=0\). Therefore
\[
    |Lg_n(x_n)|
    \geq
    \varepsilon
    \frac{|F(\zeta_b)|}{\|F\|_{L^p}}
    >0,
\]
contradicting \eqref{eq-compact-difference-Lgn-point}.

\smallskip
\noindent
Hence the assumption that both individual operators are noncompact is
impossible. Therefore, when \(\varphi\not\equiv\phi\), both $W_{\psi,\varphi}$ and $W_{\vartheta,\phi}$ are compact. This proves the first condition and completes the proof.
\end{proof}

\subsection{The case \(0<q<p\leq\infty\)}

In this range, boundedness of a difference must force the weights themselves
to provide the missing integrability. The main difficulty is to separate the
two summands locally when their symbols are distinct. We begin by constructing
packets that are uniformly large along one affine symbol and uniformly small
along the other.

\begin{lemma}\label{lem-separating-packets}
Let \(\sigma>0\) and \(0<p\leq\infty\). Let
\[
    \varphi(z)=az+b,\qquad
    \phi(z)=cz+d,
\qquad \text{where} \qquad
    a,c\in\R,\qquad a\neq0,\qquad
    \varphi\not\equiv\phi.
\]
Then there exist a partition
\[
    \R=\bigcup_{n\in\Z}E_n
\]
into pairwise disjoint half-open intervals of equal length, a sequence
\((h_n)_{n\in\Z}\subset B_\sigma^p\), a bounded interval \(K\subset\R\),
and a constant \(0<\kappa<1\) such that the following properties hold:
\begin{enumerate}
\item[\textnormal{(a)}] For every finite set \(A\subset\Z\) and every complex sequence
\((\lambda_n)_{n\in A}\),
\[
    \left\|
        \sum_{n\in A}\lambda_n h_n
    \right\|_{L^p}
    \lesssim
    \|(\lambda_n)_{n\in A}\|_{\ell^p}.
\]

\item[\textnormal{(b)}] For every \(n\in\Z\),
\[
    |h_n(\varphi(x))|\geq1,
    \qquad x\in E_n.
\]

\item[\textnormal{(c)}] For every \(n\in\Z\),
\[
    |h_n(\phi(x))|\leq\kappa,
    \qquad x\in E_n\setminus K.
\]
\end{enumerate}
The implicit constant in \emph{(a)}, as well as \(K\) and \(\kappa\),
may depend on \(\sigma,p,\varphi\), and \(\phi\), but not on \(n\), \(A\),
or the coefficients \((\lambda_n)_{n \in A}\).
\end{lemma}

\begin{proof}
We distinguish the cases \(a\neq c\) and \(a=c\).

\smallskip
\noindent
\emph{Case 1. \(a\neq c\).}
Put $\rho:=\frac{\sigma}{2},$ and let \(h_{\rho,n}\), \(n\in\Z\), be the localized packets from
Lemma~\ref{lem:localized-packets}, corresponding to the exponent \(p\).
Thus
\[
    h_{\rho,n}(z)
    =
    H\left(\rho z-n-\frac12\right),
\qquad \text{where} \qquad
    H(z)
    =
    \left(
        \frac{\sin(z/N_p)}{z/N_p}
    \right)^{N_p},
\]
and there exists \(\eta>0\), depending only on \(p\), such that
\[
    h_{\rho,n}(x)\geq\eta,
    \qquad x\in I_{\rho,n}.
\]
Moreover,
\begin{equation}\label{eq-separating-packets-synthesis-base}
    \left\|
        \sum_{n\in A}\lambda_n h_{\rho,n}
    \right\|_{L^p}
    \lesssim_p
    \rho^{-1/p}
    \|(\lambda_n)_{n\in A}\|_{\ell^p}
\end{equation}
for every finite \(A\subset\Z\).

\smallskip
\noindent
Define
\[
    E_n
    :=
    \{x\in\R:a x+\re b\in I_{\rho,n}\},
    \qquad n\in\Z.
\]
Since \(a\neq0\), the sets \(E_n\) form a partition of \(\R\) into
pairwise disjoint half-open intervals of length $\frac{1}{\rho|a|}.$ For \(n\in\Z\), set
\[
    h_n(z)
    :=
    \frac{1}{\eta}
    h_{\rho,n}(z-i\im b).
\]
Each \(h_n\) belongs to \(B_\rho^p\subset B_\sigma^p\). By the vertical
translation estimate and \eqref{eq-separating-packets-synthesis-base},
\[
\begin{aligned}
    \left\|
        \sum_{n\in A}\lambda_n h_n
    \right\|_{L^p}
    =
    \frac1\eta
    \left\|
        T_{-\im b}
        \left(
            \sum_{n\in A}\lambda_nh_{\rho,n}
        \right)
    \right\|_{L^p}
    \lesssim_p
    e^{\rho|\im b|}
    \rho^{-1/p}
    \|(\lambda_n)_{n\in A}\|_{\ell^p}.
\end{aligned}
\]
This proves (a).

\smallskip
\noindent
If \(x\in E_n\), then $ a x+\re b\in I_{\rho,n}.$ Since
\(
    \varphi(x)-i\im b=a x+\re b,
\)
we obtain
\[
    |h_n(\varphi(x))|
    =
    \frac1\eta
    h_{\rho,n}(a x+\re b)
    \geq1.
\]
Thus (b) holds.

\smallskip
\noindent
It remains to estimate \(h_n(\phi(x))\). From the explicit form of \(H\)
and the elementary estimate for the sinc function, for every \(M>0\),
\begin{equation}\label{eq-separating-profile-decay}
    |H(u+iv)|
    \lesssim_{p,M}
    (1+|u|)^{-N_p},
    \qquad u\in\R,\quad |v|\leq M.
\end{equation}
Let \(x\in E_n\), and write
\[
    u
    :=
    \rho(a x+\re b)-n-\frac12.
\]
Then $|u|\leq\frac12.$ Moreover,
\[
\begin{aligned}
    \rho\bigl(\phi(x)-i\im b\bigr)-n-\frac12
    &=
    \rho(c-a)x
    +\rho(\re d-\re b)
    +u
    +i\rho(\im d-\im b).
\end{aligned}
\]
Hence
\[
\begin{aligned}
    \left|
        \re\left(
            \rho\bigl(\phi(x)-i\im b\bigr)
            -n-\frac12
        \right)
    \right|
    \geq
    \rho|a-c|\,|x|
    -\rho|\re d-\re b|
    -\frac12.
\end{aligned}
\]
Since \(a\neq c\), the right-hand side tends to infinity as
\(|x|\to\infty\), uniformly in \(n\) and \(x\in E_n\). The imaginary
part $\rho(\im d-\im b)$ is fixed. Therefore \eqref{eq-separating-profile-decay} gives
\[
    \sup_{\substack{n\in\Z,\ x\in E_n\\ |x|\geq R}}
    |h_n(\phi(x))|
    \longrightarrow0
    \qquad (R\to\infty).
\]
Choose \(R>0\) sufficiently large that
\[
    |h_n(\phi(x))|\leq\frac12,
    \qquad
    x\in E_n,\quad |x|\geq R.
\]
Thus (c) holds with
\[
    K=[-R,R],
    \qquad
    \kappa=\frac12.
\]

\smallskip
\noindent
\emph{Case 2. \(a=c\).}
Since \(\varphi\not\equiv\phi\), we have \(b\neq d\). Put
\[
    \zeta_b:=i\im b,
    \qquad
    \zeta_d:=\re d-\re b+i\im d.
\]
Then $\zeta_b\neq\zeta_d.$ Choose \(0<\rho<\sigma\) and an integer \(m\) such that
\[
    m>\max\left\{2,1+\frac1p\right\}.
\]
Define
\[
    F(z)
    :=
    (z-\zeta_d)G_{\rho,m,\zeta_b}(z).
\]
Since multiplication by a polynomial does not increase exponential type,
\(F\) has exponential type at most \(\rho\). Moreover, the choice of \(m\)
and the estimate used in the proof of Lemma~\ref{lem:packet-decay} give
\begin{equation}\label{eq-separating-F-decay}
    |F(x+iy)|
    \lesssim_M
    (1+|x|)^{1-m},
    \qquad x\in\R,\quad |y|\leq M,
\end{equation}
for every \(M>0\). In particular,
\[
    F\in B_\rho^p\subset B_\sigma^p.
\]
Also, $F(\zeta_d)=0,$ whereas
\[
    F(\zeta_b)
    =
    (\zeta_b-\zeta_d)
    G_{\rho,m,\zeta_b}(\zeta_b)
    \neq0.
\]
Put
\[
    M_0:=|F(\zeta_b)|>0.
\]
By continuity, there exists \(\delta>0\) such that
\[
    |F(u+\zeta_b)|\geq\frac{M_0}{2},
    \qquad |u|\leq\delta,
\qquad \text{and} \qquad
    |F(u+\zeta_d)|\leq\frac{M_0}{4},
    \qquad |u|\leq\delta.
\]
Set
\[
    P(z):=\frac{2}{M_0}F(z).
\]
Then
\begin{equation}\label{eq-separating-parallel-values}
    |P(u+\zeta_b)|\geq1,
    \qquad
    |P(u+\zeta_d)|\leq\frac12,
    \qquad |u|\leq\delta.
\end{equation}
For \(n\in\Z\), put
\[
    s_n:=2\delta n,
    \qquad
    J_n:=[s_n-\delta,s_n+\delta),
\]
and define
\[
    E_n
    :=
    \{x\in\R:a x+\re b\in J_n\}.
\]
Again, \((E_n)_{n\in\Z}\) is a partition of \(\R\) into pairwise
disjoint half-open intervals, now of length $\frac{2\delta}{|a|}.$ Set
\[
    h_n(z):=P(z-s_n),
    \qquad n\in\Z.
\]
If \(x\in E_n\), put $u:=a x+\re b-s_n.$ Then \(|u|\leq\delta\), and
\[
    \varphi(x)-s_n=u+\zeta_b,
\]
while, since \(a=c\),
\[
    \phi(x)-s_n=u+\zeta_d.
\]
It follows from \eqref{eq-separating-parallel-values} that
\[
    |h_n(\varphi(x))|\geq1
\qquad \text{and} \qquad
    |h_n(\phi(x))|\leq\frac12,
    \qquad x\in E_n.
\]
Thus (b) and (c) hold, with \(\kappa=\frac{1}{2}\) and any fixed bounded interval
\(K\).

\smallskip
\noindent
It remains to verify (a). Suppose first that \(0<p<1\). By the
\(p\)-triangle inequality and real translation invariance,
\[
\begin{aligned}
    \left\|
        \sum_{n\in A}\lambda_nh_n
    \right\|_{L^p}^p
    \leq
    \sum_{n\in A}
    |\lambda_n|^p\|h_n\|_{L^p}^p
    =
    \|P\|_{L^p}^p
    \sum_{n\in A}|\lambda_n|^p.
\end{aligned}
\]
Hence
\[
    \left\|
        \sum_{n\in A}\lambda_nh_n
    \right\|_{L^p}
    \leq
    \|P\|_{L^p}
    \|(\lambda_n)_{n\in A}\|_{\ell^p}.
\]
Now suppose that \(1\leq p\leq\infty\). For \(j\in\Z\), put
\[
    \alpha_j
    :=
    \sup_{|u|\leq\delta}
    |P(2\delta j+u)|.
\]
By \eqref{eq-separating-F-decay} and \(m>2\),
\[
    \alpha: = (\alpha_j)_{j\in\Z}\in\ell^1.
\]
If \(x\in J_k\), then $x=s_k+u$ for some \(|u|\leq\delta\), and therefore
\[
    |h_n(x)|
    =
    |P(2\delta(k-n)+u)|
    \leq \alpha_{k-n}.
\]
Extending \((\lambda_n)_{n\in A}\) by zero outside \(A\), we obtain
\[
    \left|
        \sum_{n\in A}\lambda_nh_n(x)
    \right|
    \leq
    \sum_{n\in\Z}|\lambda_n|\alpha_{k-n},
    \qquad x\in J_k.
\]
If \(1\leq p<\infty\), then
\[
\begin{aligned}
    \left\|
        \sum_{n\in A}\lambda_nh_n
    \right\|_{L^p}^p
    \leq
    2\delta
    \sum_{k\in\Z}
    \left(
        \sum_{n\in\Z}|\lambda_n|\alpha_{k-n}
    \right)^p
    =
    2\delta
    \|\alpha*|\lambda|\|_{\ell^p}^p.
\end{aligned}
\]
Young's convolution inequality gives
\[
    \left\|
        \sum_{n\in A}\lambda_nh_n
    \right\|_{L^p}
    \leq
    (2\delta)^{1/p}
    \|\alpha\|_{\ell^1}
    \|(\lambda_n)_{n\in A}\|_{\ell^p}.
\]
For \(p=\infty\), the same pointwise estimate yields
\[
    \left\|
        \sum_{n\in A}\lambda_nh_n
    \right\|_{L^\infty}
    \leq
    \|\alpha\|_{\ell^1}
    \|(\lambda_n)_{n\in A}\|_{\ell^\infty}.
\]
This proves (a) and completes the proof.
\end{proof}

\begin{lemma}\label{lem-difference-local-mass}
Let \(\sigma,\tau>0\) and \(0<q<p\leq\infty\). Let
\[
    \varphi(z)=az+b,\qquad
    \phi(z)=cz+d,
\qquad \text{where} \qquad
    a,c\in\R,\qquad a\neq0,\qquad
    \varphi\not\equiv\phi.
\]
Suppose that the difference
\(
    L    :B_\sigma^p\longrightarrow B_\tau^q
\)
is bounded. Let \((E_n)_{n\in\Z}\), \((h_n)_{n\in\Z}\), \(K\), and
\(0<\kappa<1\) be given by Lemma~\ref{lem-separating-packets}, and set
\[
    A_n
    :=
    \int_{E_n\setminus K}
    \bigl(|\psi(x)|-\kappa|\vartheta(x)|\bigr)_+^q\,dx,
    \qquad n\in\Z,
\qquad \text{where} \qquad
    t_+:=\max\{t,0\}.
\]
Then, for every finite set \(A\subset\Z\) and every complex sequence
\((\lambda_n)_{n\in A}\),
\begin{equation}\label{eq-difference-local-mass}
    \sum_{n\in A}|\lambda_n|^q A_n
    \lesssim
    \|L\|^q
    \|(\lambda_n)_{n\in A}\|_{\ell^p}^q.
\end{equation}
Consequently,
\[
    (A_n)_{n\in\Z}
    \in
    \ell^{\left(\frac pq\right)'}
\qquad \text{and} \qquad
    \|(A_n)_{n\in\Z}\|_{\ell^{\bigl(\frac pq\bigr)'}}
    \lesssim
    \|L\|^q.
\]
The implicit constants depend only on \(p,q\) and on the constants in
Lemma~\ref{lem-separating-packets}.
\end{lemma}

\begin{proof}
Let \((r_n)_{n\in\Z}\) be the Rademacher functions on \([0,1]\). For a
finite set \(A\subset\Z\), a complex sequence
\((\lambda_n)_{n\in A}\), and \(t\in[0,1]\), define
\[
    f_t
    :=
    \sum_{n\in A}r_n(t)\lambda_n h_n.
\]
Since \(|r_n(t)|=1\), Lemma~\ref{lem-separating-packets}(a) gives
\[
    \|f_t\|_{L^p}
    \lesssim
    \|(\lambda_n)_{n\in A}\|_{\ell^p},
    \qquad 0\leq t\leq1.
\]
Hence, by the boundedness of \(L\),
\[
    \|Lf_t\|_{L^q}
    \lesssim
    \|L\|\cdot
    \|(\lambda_n)_{n\in A}\|_{\ell^p}.
\]
Since \(q<\infty\), raising to the power \(q\) and integrating over
\([0,1]\), we obtain
\begin{equation}\label{eq-difference-rademacher-upper}
    \int_0^1\|Lf_t\|_{L^q}^q\,dt
    \lesssim
    \|L\|^q
    \|(\lambda_n)_{n\in A}\|_{\ell^p}^q.
\end{equation}
Since the sum defining \(f_t\) is finite and \(L\) is linear,
\[
    Lf_t
    =
    \sum_{n\in A}r_n(t)\lambda_nLh_n.
\]
Therefore, by Tonelli's theorem,
\[
\begin{aligned}
    \int_0^1\|Lf_t\|_{L^q}^q\,dt
    &=
    \int_{\R}\int_0^1
    \left|
        \sum_{n\in A}
        r_n(t)\lambda_nLh_n(x)
    \right|^q
    dt\,dx.
\end{aligned}
\]
Khinchine's inequality, valid for every \(0<q<\infty\), gives
\[
    \int_0^1
    \left|
        \sum_{n\in A}
        r_n(t)\lambda_nLh_n(x)
    \right|^q
    dt
    \simeq_q
    \left(
        \sum_{n\in A}
        |\lambda_n|^2|Lh_n(x)|^2
    \right)^{q/2}.
\]
Combining this with \eqref{eq-difference-rademacher-upper}, we get
\begin{equation}\label{eq-difference-square-function}
    \int_{\R}
    \left(
        \sum_{n\in A}
        |\lambda_n|^2|Lh_n(x)|^2
    \right)^{q/2}
    dx
    \lesssim
    \|L\|^q
    \|(\lambda_n)_{n\in A}\|_{\ell^p}^q.
\end{equation}
For \(n\in A\) and \(x\in E_n\setminus K\),
Lemma~\ref{lem-separating-packets}(b)--(c) gives
\[
    |h_n(\varphi(x))|\geq1,
    \qquad
    |h_n(\phi(x))|\leq\kappa.
\]
Hence
\[
\begin{aligned}
    |Lh_n(x)|
    &=
    \bigl|
        \psi(x)h_n(\varphi(x))
        -
        \vartheta(x)h_n(\phi(x))
    \bigr| \\
    &\geq
    |\psi(x)|\,|h_n(\varphi(x))|
    -
    |\vartheta(x)|\,|h_n(\phi(x))| \\
    &\geq
    |\psi(x)|-\kappa|\vartheta(x)|.
\end{aligned}
\]
Since the left-hand side is nonnegative,
\[
    |Lh_n(x)|
    \geq
    \bigl(|\psi(x)|-\kappa|\vartheta(x)|\bigr)_+.
\]
Moreover,
\[
    \left(
        \sum_{j\in A}
        |\lambda_j|^2|Lh_j(x)|^2
    \right)^{q/2}
    \geq
    |\lambda_n|^q|Lh_n(x)|^q.
\]
Since the intervals \(E_n\) are pairwise disjoint, it follows that
\[
    \sum_{n\in A}|\lambda_n|^qA_n
    \leq
    \sum_{n\in A}
    |\lambda_n|^q
    \int_{E_n\setminus K}|Lh_n(x)|^q\,dx \leq
    \int_{\R}
    \left(
        \sum_{j\in A}
        |\lambda_j|^2|Lh_j(x)|^2
    \right)^{q/2}
    dx.
\]
Together with \eqref{eq-difference-square-function}, this proves
\eqref{eq-difference-local-mass}.

\smallskip
\noindent
It remains to derive the asserted sequence condition. Suppose first that
\(p<\infty\). Then
\[
    \left(\frac pq\right)'
    =
    \frac{p}{p-q}.
\]
Let \(A\subset\Z\) be finite and set
\[
    \lambda_n:=A_n^{\frac{1}{p-q}},
    \qquad n\in A.
\]
Then \eqref{eq-difference-local-mass} gives
\[
    \sum_{n\in A}A_n^{\frac{p}{p-q}}
    \lesssim
    \|L\|^q
    \left(
        \sum_{n\in A}A_n^{\frac{p}{p-q}}
    \right)^{\frac{q}{p}}.
\]
If the sum on the left is nonzero, division yields
\[
    \left(
        \sum_{n\in A}A_n^{\frac{p}{p-q}}
    \right)^{\frac{p-q}{p}}
    \lesssim
    \|L\|^q.
\]
Since \(A\) is arbitrary,
\[
    (A_n)_{n\in\Z}
    \in
    \ell^{\frac{p}{p-q}}
    =
    \ell^{\bigl(\frac pq\bigr)'}
\qquad \text{and} \qquad
    \|(A_n)_{n\in\Z}\|_{\ell^{\bigl(\frac pq\bigr)'}}
    \lesssim
    \|L\|^q.
\]
Finally, let \(p=\infty\). For a finite set \(A\subset\Z\), take
\[
    \lambda_n=1,
    \qquad n\in A.
\]
Then \eqref{eq-difference-local-mass} gives
\[
    \sum_{n\in A}A_n
    \lesssim
    \|L\|^q.
\]
Since \(A\) is arbitrary,
\[
    (A_n)_{n\in\Z}\in\ell^1
    =
    \ell^{\bigl(\frac{\infty}{q}\bigr)'}
\qquad \text{and} \qquad
    \|(A_n)_{n\in\Z}\|_{\ell^1}
    \lesssim
    \|L\|^q.
\]
This completes the proof.
\end{proof}

\begin{lemma}\label{lem-difference-decoupling}
Let \(0<q<p\leq\infty\). Let $(E_n)_{n\in\Z}$ and $(\widetilde E_m)_{m\in\Z}$ be partitions of \(\R\) into half-open intervals of fixed positive lengths.
Suppose that there exist bounded intervals \(K,\widetilde K\subset\R\)
and constants \(0<\kappa_1,\kappa_2<1\) such that
\[
    A_n
    :=
    \int_{E_n\setminus K}
    \bigl(|\psi(x)|-\kappa_1|\vartheta(x)|\bigr)_+^q\,dx,
    \qquad \text{and} \qquad
    B_m
    :=
    \int_{\widetilde E_m\setminus\widetilde K}
    \bigl(|\vartheta(x)|-\kappa_2|\psi(x)|\bigr)_+^q\,dx,
\]
satisfy
\[
    (A_n)_{n\in\Z},
    \ (B_m)_{m\in\Z}
    \in
    \ell^{(\frac pq)'}.
\]
Then
\[
    \left(
        \int_{E_n}|\psi(x)|^q\,dx
    \right)_{n\in\Z}
    \in
    \ell^{(\frac pq)'}
\qquad \text{and} \qquad
    \left(
        \int_{E_n}|\vartheta(x)|^q\,dx
    \right)_{n\in\Z}
    \in
    \ell^{(\frac pq)'}.
\]
\end{lemma}

\begin{proof}
Choose a bounded interval \(K_0\) containing $K\cup\widetilde K,$ and put $\kappa:=\max\{\kappa_1,\kappa_2\}.$ Then \(0<\kappa<1\), and
\[
    \bigl(|\psi|-\kappa|\vartheta|\bigr)_+
    \leq
    \bigl(|\psi|-\kappa_1|\vartheta|\bigr)_+,
\qquad \text{while} \qquad
    \bigl(|\vartheta|-\kappa|\psi|\bigr)_+
    \leq
    \bigl(|\vartheta|-\kappa_2|\psi|\bigr)_+.
\]
Consequently, if
\[
    \alpha_n
    :=
    \int_{E_n\setminus K_0}
    \bigl(|\psi(x)|-\kappa|\vartheta(x)|\bigr)_+^q\,dx
\qquad \text{and} \qquad
    \beta_m
    :=
    \int_{\widetilde E_m\setminus K_0}
    \bigl(|\vartheta(x)|-\kappa|\psi(x)|\bigr)_+^q\,dx,
\]
then
\[
    0\leq\alpha_n\leq A_n,
    \qquad
    0\leq\beta_m\leq B_m.
\]
Hence
\begin{equation}\label{eq-decoupling-alpha-beta}
    (\alpha_n)_{n\in\Z},
    \ (\beta_m)_{m\in\Z}
    \in
    \ell^{(\frac pq)'}.
\end{equation}
We first transfer the second sequence to the partition \((E_n)\). Put
\[
    \gamma_n
    :=
    \int_{E_n\setminus K_0}
    \bigl(|\vartheta(x)|-\kappa|\psi(x)|\bigr)_+^q\,dx,
    \qquad n\in\Z.
\]
For \(n\in\Z\), let
\[
    \mathcal N(n)
    :=
    \{m\in\Z:E_n\cap\widetilde E_m\neq\varnothing\}.
\]
Since both partitions consist of intervals of fixed positive lengths,
there exist integers \(N_1,N_2\geq1\) such that
\[
    \#\mathcal N(n)\leq N_1,
    \qquad n\in\Z,
\]
and each interval \(\widetilde E_m\) intersects at most \(N_2\) of the
intervals \(E_n\). Thus
\[
    \gamma_n
    \leq
    \sum_{m\in\mathcal N(n)}\beta_m.
\]
Since
\(
    \left(\frac pq\right)'\geq1,
\)
the elementary inequality for finite sums gives
\[
    \gamma_n^{\bigl(\frac pq\bigr)'}
    \lesssim
    \sum_{m\in\mathcal N(n)}
    \beta_m^{\bigl(\frac pq\bigr)'}.
\]
Summing over \(n\) and using the bounded overlap in the reverse direction,
we obtain
\[
    \sum_{n\in\Z}
    \gamma_n^{\bigl(\frac pq\bigr)'}
    \lesssim
    \sum_{m\in\Z}
    \beta_m^{\bigl(\frac pq\bigr)'}.
\]
Hence
\begin{equation}\label{eq-decoupling-gamma}
    (\gamma_n)_{n\in\Z}
    \in
    \ell^{(\frac pq)'}.
\end{equation}
We now remove the mixed terms. For \(r,t\geq0\),
\begin{equation}\label{eq-decoupling-elementary}
    r^q+t^q
    \lesssim_{\kappa,q}
    (r-\kappa t)_+^q
    +(t-\kappa r)_+^q.
\end{equation}
Indeed, if \(r\geq t\), then $r-\kappa t
    \geq
    (1-\kappa)r,$ and therefore
\[
    r^q+t^q
    \leq
    2r^q
    \leq
    \frac{2}{(1-\kappa)^q}
    (r-\kappa t)_+^q.
\]
The case \(t\geq r\) is symmetric.

\smallskip
\noindent
Applying \eqref{eq-decoupling-elementary} with $r=|\psi(x)|$, $t=|\vartheta(x)|$ and integrating over \(E_n\setminus K_0\), we obtain
\[
\begin{aligned}
    \int_{E_n\setminus K_0}|\psi(x)|^q\,dx
    +
    \int_{E_n\setminus K_0}|\vartheta(x)|^q\,dx
    \lesssim_{\kappa,q}
    \alpha_n+\gamma_n.
\end{aligned}
\]
In particular,
\[
    \int_{E_n\setminus K_0}|\psi(x)|^q\,dx
    \lesssim_{\kappa,q}
    \alpha_n+\gamma_n,\quad \int_{E_n\setminus K_0}|\vartheta(x)|^q\,dx
    \lesssim_{\kappa,q}
    \alpha_n+\gamma_n.
\]
By \eqref{eq-decoupling-alpha-beta} and
\eqref{eq-decoupling-gamma},
\[
    \left(
        \int_{E_n\setminus K_0}|\psi(x)|^q\,dx
    \right)_{n\in\Z}
    \in
    \ell^{\bigl(\frac pq\bigr)'},\quad\left(
        \int_{E_n\setminus K_0}|\vartheta(x)|^q\,dx
    \right)_{n\in\Z}
    \in
    \ell^{\bigl(\frac pq\bigr)'}.
\]
Finally, only finitely many intervals \(E_n\) intersect \(K_0\).
Since \(\psi\) and \(\vartheta\) are entire, they belong locally to
\(L^q(\R)\). Hence the sequences
\[
    \left(
        \int_{E_n\cap K_0}|\psi(x)|^q\,dx
    \right)_{n\in\Z},
    \qquad
    \left(
        \int_{E_n\cap K_0}|\vartheta(x)|^q\,dx
    \right)_{n\in\Z}
\]
are finitely supported. Adding these terms proves
\[
    \left(
        \int_{E_n}|\psi(x)|^q\,dx
    \right)_{n\in\Z}
    \in
    \ell^{\bigl(\frac pq\bigr)'},\quad\left(
        \int_{E_n}|\vartheta(x)|^q\,dx
    \right)_{n\in\Z}
    \in
    \ell^{\bigl(\frac pq\bigr)'}.
\]
\end{proof}

\begin{theorem}\label{thm-difference-q-less-p}
Let \(\sigma,\tau>0\) and \(0<q<p\leq\infty\). Let
\(\psi,\vartheta,\varphi,\phi\) be entire functions. Then the following
conditions are equivalent:
\begin{enumerate}
\item[\textnormal{(a)}]
The difference \(
    W_{\psi,\varphi}-W_{\vartheta,\phi}
    :B_\sigma^p\longrightarrow B_\tau^q
\)
is bounded.

\item[\textnormal{(b)}]
The difference \(
    W_{\psi,\varphi}-W_{\vartheta,\phi}
    :B_\sigma^p\longrightarrow B_\tau^q
\)
is compact.

\item[\textnormal{(c)}] At least one of the following alternatives holds:
\begin{enumerate}
\item[\textnormal{(i)}] \(\varphi\not\equiv\phi\) and both
\(
    W_{\psi,\varphi},\,
    W_{\vartheta,\phi}
    :B_\sigma^p\longrightarrow B_\tau^q
\)
are bounded, and hence compact;

\item[\textnormal{(ii)}] \(\varphi=\phi\) and
\(
    W_{\psi-\vartheta,\varphi}
    :B_\sigma^p\longrightarrow B_\tau^q
\)
is bounded, and hence compact.
\end{enumerate}
\end{enumerate}
\end{theorem}

\begin{proof}
The implication (b)\(\Rightarrow\)(a) is immediate.

\smallskip
\noindent
We prove (a)\(\Rightarrow\)(c). Suppose that the difference
\(
    L:B_\sigma^p\longrightarrow B_\tau^q
\)
is bounded. If $\varphi=\phi,$ then $L=W_{\psi-\vartheta,\varphi},$ and (c-ii) holds.

\smallskip
\noindent
Suppose now that $\varphi\not\equiv\phi.$ If one of the weights is identically zero, then the corresponding weighted
composition operator is the zero operator, while the other one coincides,
up to sign, with \(L\). Thus (c-i) holds. We may therefore assume
that $\psi\not\equiv0$ and $\vartheta\not\equiv0$.

\smallskip
\noindent
By Proposition~\ref{prop-difference-structural},
\[
    \varphi(z)=az+b,\qquad
    \phi(z)=cz+d,
\qquad \text{where} \qquad
    a,c\in\R,\qquad
    \sigma|a|\leq\tau,\qquad
    \sigma|c|\leq\tau.
\]
Moreover,
\begin{equation}\label{eq-difference-q-less-p-type}
    \operatorname{type}(\psi)\leq\tau-\sigma|a|,
    \qquad
    \operatorname{type}(\vartheta)\leq\tau-\sigma|c|.
\end{equation}
We distinguish three cases.

\smallskip
\noindent
\emph{Case 1. Both symbols are nonconstant.}
Suppose that $a\neq0,c\neq0.$ 
Apply Lemma~\ref{lem-separating-packets} to the ordered pair
\((\varphi,\phi)\). Then Lemma~\ref{lem-difference-local-mass} gives a
partition \((E_n)_{n\in\Z}\) of \(\R\) into intervals of fixed positive
length, a bounded interval \(K\), and a number \(0<\kappa_1<1\) such that
\[
    \left(
        \int_{E_n\setminus K}
        \bigl(|\psi(x)|-\kappa_1|\vartheta(x)|\bigr)_+^q\,dx
    \right)_{n\in\Z}
    \in
    \ell^{\bigl(\frac pq\bigr)'}.
\]
Interchanging the two ordered pairs and observing that
\[
    W_{\vartheta,\phi}-W_{\psi,\varphi}=-L,
\]
we obtain another partition
\((\widetilde E_m)_{m\in\Z}\), a bounded interval \(\widetilde K\), and
a number \(0<\kappa_2<1\) such that
\[
    \left(
        \int_{\widetilde E_m\setminus\widetilde K}
        \bigl(|\vartheta(x)|-\kappa_2|\psi(x)|\bigr)_+^q\,dx
    \right)_{m\in\Z}
    \in
    \ell^{\bigl(\frac pq\bigr)'}.
\]
Lemma~\ref{lem-difference-decoupling} therefore yields
\begin{equation}\label{eq-difference-q-less-p-local-psi}
    \left(
        \int_{E_n}|\psi(x)|^q\,dx
    \right)_{n\in\Z}
    \in
    \ell^{\bigl(\frac pq\bigr)'}
\end{equation}
and
\begin{equation}\label{eq-difference-q-less-p-local-vartheta}
    \left(
        \int_{E_n}|\vartheta(x)|^q\,dx
    \right)_{n\in\Z}
    \in
    \ell^{\bigl(\frac pq\bigr)'}.
\end{equation}
Let \(\ell>0\) be the common length of the intervals \(E_n\), and put $\rho:=\frac1\ell.$ Since the two partitions $(E_n)_{n\in\Z}$ and $(I_{\rho,k})_{k\in\Z}$ consist of intervals of the same fixed length, they have uniformly bounded
mutual overlap. Since
\(
    \left(\frac pq\right)'\geq1,
\)
the bounded-overlap argument applied to
\eqref{eq-difference-q-less-p-local-psi} and
\eqref{eq-difference-q-less-p-local-vartheta} gives
\begin{equation}\label{eq-difference-q-less-p-natural-cells}
    \left(
        \int_{I_{\rho,k}}|\psi(x)|^q\,dx
    \right)_{k\in\Z},
    \quad
    \left(
        \int_{I_{\rho,k}}|\vartheta(x)|^q\,dx
    \right)_{k\in\Z}
    \in
    \ell^{\bigl(\frac pq\bigr)'}.
\end{equation}
Suppose first that \(p<\infty\). Then
\(
    \left(\frac pq\right)'
    =
    \frac{p}{p-q}
    =
    \frac1q\frac{pq}{p-q}.
\)
Since both weights have finite exponential type, Lemma~\ref{lem:local-q-to-p},
applied with the exponent pair
\[
    q<\frac{pq}{p-q},
\]
and with the partition \((I_{\rho,k})\), gives
\[
    \psi,\vartheta
    \in
    L^{\frac{pq}{p-q}}(\R).
\]
If \(p=\infty\), then
\(
    \left(\frac{\infty}{q}\right)'=1,
\)
and \eqref{eq-difference-q-less-p-natural-cells} gives directly
\[
    \int_{\R}|\psi(x)|^q\,dx<\infty,
    \qquad
    \int_{\R}|\vartheta(x)|^q\,dx<\infty.
\]
Thus,
we have in all cases
\begin{equation}\label{eq-difference-q-less-p-integrability}
    \psi,\vartheta
    \in
    L^{\frac{pq}{p-q}}(\R).
\end{equation}
We next show that the bandwidth inequalities are strict. Suppose, for
example, that $\sigma|a|=\tau.$ Then \eqref{eq-difference-q-less-p-type} shows that \(\psi\) has
exponential type zero. Hence, for every \(\varepsilon>0\),
\[
    \psi
    \in
    B_\varepsilon^{\frac{pq}{p-q}}.
\]
By Nikol'skii's inequality,
\[
    \|\psi\|_{L^\infty}
    \lesssim
    \varepsilon^{\,\frac1q-\frac1p}
    \|\psi\|_{L^{\frac{pq}{p-q}}}.
\]
Letting \(\varepsilon\to0^+\) gives $\psi\equiv0,$ a contradiction. Therefore, $\sigma|a|<\tau.$ Similarly, $\sigma|c|<\tau.$ Combining this with
\eqref{eq-difference-q-less-p-type} and
\eqref{eq-difference-q-less-p-integrability}, we obtain
\[
    \psi
    \in
    B_{\tau-\sigma|a|}^{\frac{pq}{p-q}},
    \qquad
    \vartheta
    \in
    B_{\tau-\sigma|c|}^{\frac{pq}{p-q}}.
\]
Theorem~\ref{thm-wco-q-lt-p} now shows that both
\[
    W_{\psi,\varphi},
    \,\,
    W_{\vartheta,\phi}
    :B_\sigma^p\longrightarrow B_\tau^q
\]
are bounded, and hence compact.

\smallskip
\noindent
\emph{Case 2. Exactly one symbol is constant.}
By symmetry, suppose that $a\neq0$, $c=0,$ so that $\phi\equiv d.$ 
We first isolate the weight \(\psi\). Put $\rho:=\frac{\sigma}{2},$ and let \(h_{\rho,n}\), \(n\in\Z\), be the localized packets from
Lemma~\ref{lem:localized-packets}, corresponding to the exponent \(p\).
Let \(\eta>0\) be the constant in
Lemma~\ref{lem:localized-packets}(b), and set
\[
    E_n
    :=
    \{x\in\R:a x+\re b\in I_{\rho,n}\}.
\]
Define
\[
    g_n(z)
    :=
    \frac1\eta h_{\rho,n}(z-i\im b).
\]
Then \((E_n)_{n\in\Z}\) is a partition of \(\R\) into intervals of fixed
positive length,
\[
    |g_n(\varphi(x))|\geq1,
    \qquad x\in E_n,
\]
and Lemma~\ref{lem:localized-packets}(d), together with the vertical
translation estimate, gives
\begin{equation}\label{eq-difference-constant-symbol-gn-synthesis}
    \left\|
        \sum_{n\in A}\lambda_ng_n
    \right\|_{L^p}
    \lesssim
    \|(\lambda_n)_{n\in A}\|_{\ell^p}
\end{equation}
for every finite \(A\subset\Z\).

\smallskip
\noindent
Choose \(u\in B_\sigma^p\) such that $u(d)=1;$ for instance, one may take
\[
    u(z)
    :=
    \frac{G_{\sigma,m,d}(z)}
         {G_{\sigma,m,d}(d)}
\]
with \(m>\frac{1}{p}\). Put $\widetilde h_n(z)
    :=
    g_n(z)-g_n(d)u(z).$ Then $\widetilde h_n(d)=0.$ The estimate used in the proof of
Lemma~\ref{lem:localized-packets}(c) gives
\[
    |g_n(d)|
    \lesssim
    (1+|n|)^{-N_p}.
\]
Since \(N_p>1\), this estimate, together with the \(p\)-triangle inequality
when \(0<p<1\), Hölder's inequality when \(1\leq p<\infty\), and the
\(\ell^1\)-estimate when \(p=\infty\), shows that
\[
    \left|
        \sum_{n\in A}\lambda_ng_n(d)
    \right|
    \lesssim
    \|(\lambda_n)_{n\in A}\|_{\ell^p}.
\]
Together with
\eqref{eq-difference-constant-symbol-gn-synthesis}, this gives
\[
    \left\|
        \sum_{n\in A}\lambda_n\widetilde h_n
    \right\|_{L^p}
    \lesssim
    \|(\lambda_n)_{n\in A}\|_{\ell^p}.
\]
Since \(a\neq0\) and \(a\in\R\), the set \(\varphi(\R)\) is a horizontal
line. Lemmas~\ref{lem:vertical-translation} and
\ref{lem:nikolskii} therefore give
\[
    M:=\sup_{x\in\R}|u(\varphi(x))|<\infty.
\]
Moreover,
\[
    g_n(d)\longrightarrow0
    \qquad (|n|\to\infty).
\]
Hence there exists \(n_0\) such that
\[
    |\widetilde h_n(\varphi(x))|
    \geq
    1-M|g_n(d)|
    \geq\frac12,
    \qquad
    x\in E_n,\quad |n|\geq n_0.
\]
Let
\[
    K:=\bigcup_{|n|<n_0}E_n
\]
and put $h_n:=2\widetilde h_n.$ Then \(K\) is bounded and
\[
    h_n(d)=0,
    \qquad
    |h_n(\varphi(x))|\geq1,
    \quad x\in E_n\setminus K,
\]
while
\[
    \left\|
        \sum_{n\in A}\lambda_nh_n
    \right\|_{L^p}
    \lesssim
    \|(\lambda_n)_{n\in A}\|_{\ell^p}.
\]
Repeating the Rademacher argument in the proof of
Lemma~\ref{lem-difference-local-mass}, we obtain
\[
    \left(
        \int_{E_n\setminus K}|\psi(x)|^q\,dx
    \right)_{n\in\Z}
    \in
    \ell^{\bigl(\frac pq\bigr)'}.
\]
Since \(K\) is bounded and only finitely many intervals \(E_n\) meet \(K\),
\[
    \left(
        \int_{E_n}|\psi(x)|^q\,dx
    \right)_{n\in\Z}
    \in
    \ell^{(\frac pq)'}.
\]
Suppose first that \(p<\infty\). As in Case 1, a bounded-overlap
argument followed by Lemma~\ref{lem:local-q-to-p}, applied with the
exponent pair
\(
    q<\frac{pq}{p-q},
\)
gives
\[
    \psi\in L^{\frac{pq}{p-q}}(\R).
\]
If \(p=\infty\), then
\(
    \left(\frac{\infty}{q}\right)'=1,
\)
and hence
\[
    \sum_{n\in\Z}\int_{E_n}|\psi(x)|^q\,dx<\infty.
\]
Since \((E_n)_{n\in\Z}\) is a partition of \(\R\), it follows directly that
\[
    \psi\in L^q(\R).
\]
Thus,
we have in all cases
\[
    \psi\in L^{\frac{pq}{p-q}}(\R).
\]
The argument used in Case 1
also shows that $\sigma|a|<\tau.$ Together with \eqref{eq-difference-q-less-p-type},
\[
    \psi
    \in
    B_{\tau-\sigma|a|}^{\frac{pq}{p-q}}.
\]
By Theorem~\ref{thm-wco-q-lt-p},
\(
    W_{\psi,\varphi}
    :B_\sigma^p\longrightarrow B_\tau^q
\)
is bounded. Since
\(
    W_{\vartheta,d}
    =
    W_{\psi,\varphi}-L,
\)
the constant-symbol operator
\[
    W_{\vartheta,d}
    :B_\sigma^p\longrightarrow B_\tau^q
\]
is bounded as well. The case \(a=0\) and \(c\neq0\) is symmetric.

\smallskip
\noindent
\emph{Case 3. Both symbols are constant.}
Suppose that $\varphi\equiv b$ and $\phi\equiv d.$ Since \(\varphi\not\equiv\phi\), $b\neq d.$ Choose \(m\in\N\) such that $m>1+\frac1p$ and define
\[
    f_b(z):=(z-d)G_{\sigma,m,b}(z),
    \qquad
    f_d(z):=(z-b)G_{\sigma,m,d}(z).
\]
As in the proof of Lemma~\ref{lem-two-point-separation}, $f_b,f_d\in B_\sigma^p.$ Moreover,
\[
    f_b(d)=0,
    \qquad
    f_b(b)\neq0,
\qquad \text{and} \qquad
    f_d(b)=0,
    \qquad
    f_d(d)\neq0.
\]
Consequently,
\[
    Lf_b=f_b(b)\psi,
    \qquad
    Lf_d=-f_d(d)\vartheta.
\]
Since
\(
    Lf_b,Lf_d\in B_\tau^q,
\)
we obtain
\[
    \psi,\vartheta\in B_\tau^q.
\]
Proposition~\ref{prop-wco-constant-symbol} therefore shows that both
constant-symbol weighted composition operators are bounded, and hence
compact.

\smallskip
\noindent
Thus, whenever $\varphi\not\equiv\phi,$ both \(W_{\psi,\varphi}\) and \(W_{\vartheta,\phi}\) are bounded. Hence
(c-i) holds, and the implication
(a)\(\Rightarrow\)(c) is proved.

\smallskip
\noindent
It remains to prove (c)\(\Rightarrow\)(b). In the range
\(0<q<p\leq\infty\), every bounded weighted composition operator
\(
    B_\sigma^p\longrightarrow B_\tau^q
\)
is compact: this follows from Theorem~\ref{thm-wco-q-lt-p} for a
nonconstant symbol and from Proposition~\ref{prop-wco-constant-symbol}
for a constant symbol, while the zero operator is trivial. Hence, under
(c-i), both summands are compact, and so is their difference.
Under (c-ii), $W_{\psi,\varphi}-W_{\vartheta,\varphi}
    =
    W_{\psi-\vartheta,\varphi},$ which is compact by the same observation. Thus (b) holds.

\smallskip
\noindent
This proves (c)\(\Rightarrow\)(b) and completes the proof.
\end{proof}

\begin{remark}\label{rem-difference-rigidity}
Theorems~\ref{thm-difference-bounded-p-le-q},
\ref{thm-difference-compact-p-le-q}, and
\ref{thm-difference-q-less-p} show that, when the composition symbols are
distinct, cancellation between two weighted composition operators cannot
produce boundedness or compactness that is absent from the individual
summands. The only genuine cancellation occurs when the symbols coincide,
in which case
\[
    W_{\psi,\varphi}-W_{\vartheta,\varphi}
    =
    W_{\psi-\vartheta,\varphi}.
\]
In particular, in the range \(0<q<p\leq\infty\), boundedness and compactness
of the difference are equivalent. This rigidity will be used in the next
section to study the topological structure of weighted composition operators.
\end{remark}

\section{Topological structure}\label{sec-topological-structure}

We conclude by studying the topological structure of bounded composition
and weighted composition operators between Paley--Wiener spaces. The
difference results of Section~\ref{sec:differences} will play a central role
in separating operators induced by distinct affine symbols.

\subsection{The space of composition operators}\label{subsec-topology-composition}

For \(\sigma,\tau>0\) and \(0<p,q\leq\infty\), let
\[
    \operatorname{CO}_{p,q}(\sigma,\tau)
    :=
    \left\{
        C_\varphi:B_\sigma^p\longrightarrow B_\tau^q:
        C_\varphi\text{ is bounded}
    \right\},
\]
endowed with the operator norm topology, or with the operator quasi-norm
topology when \(q<1\).

\smallskip
\noindent
By Corollary~\ref{cor-composition-classification}, if
\(0<p\leq q\leq\infty\), every bounded nonconstant composition operator
\(C_\varphi:B_\sigma^p\to B_\tau^q\) is induced by
\[
    \varphi(z)=az+b,
    \qquad
    a\in\R,\quad 0<\sigma|a|\leq\tau,\quad b\in\C.
\]
Accordingly, for each \(a\in\R\) with \(0<\sigma|a|\leq\tau\), set
\[
    \operatorname{CO}_a
    :=
    \left\{
        C_{az+b}:b\in\C
    \right\}.
\]
When \(q=\infty\), Remark~\ref{rem-constant-composition} shows that
constant symbols also induce bounded composition operators, and we set
\[
    \operatorname{CO}_0
    :=
    \left\{
        C_b:b\in\C
    \right\},
    \qquad
    C_bf=f(b)1.
\]

\smallskip
\noindent
We first record the norm continuity of translations on a fixed
Paley--Wiener space. Vertical translations \(T_y\) were treated in
Lemma~\ref{lem:vertical-translation}. For \(h\in\mathbb C\), set
\[
S_hf(z):=f(z+h).
\]

\begin{lemma}\label{lem-translation-continuity}
Let \(\sigma>0\) and \(0<p\leq\infty\).
Then
\[
    \|S_h-I\|_{B_\sigma^p\to B_\sigma^p}
    \leq
    \left(
        \sum_{n=1}^{\infty}
        \left(
            \frac{(\sigma|h|)^n}{n!}
        \right)^{\min\{1,p\}}
    \right)^{\frac{1}{\min\{1,p\}}}.
\]
In particular,
\[
    \|S_h-I\|_{B_\sigma^p\to B_\sigma^p}
    \longrightarrow0
    \quad\text{as} \quad h\to0.
\]
\end{lemma}

\begin{proof}
Let \(f\in B_\sigma^p\). By Bernstein's inequality for entire functions of exponential type
\cite{Bernstein1923, RahmanSchmeisser1990},
\[
    \|f^{(n)}\|_{L^p}
    \leq
    \sigma^n\|f\|_{L^p},
    \qquad n\geq1.
\]
Hence the Taylor series
\[
    \sum_{n=1}^{\infty}\frac{h^n}{n!}f^{(n)}
\]
converges in \(B_\sigma^p\). By Theorem~\ref{thm:point-evaluations},
convergence in \(B_\sigma^p\) implies locally uniform convergence, so its
sum is \(S_hf-f\). Therefore, using the
\(\min\{1,p\}\)-triangle inequality,
\[
    \|S_hf-f\|_{L^p}^{\min\{1,p\}}
    \leq
    \sum_{n=1}^{\infty}
    \left(
        \frac{(\sigma|h|)^n}{n!}
    \right)^{\min\{1,p\}}
    \|f\|_{L^p}^{\min\{1,p\}}.
\]
Taking the supremum over the unit ball of \(B_\sigma^p\) gives the
asserted estimate, and the final conclusion follows immediately.
\end{proof}

\begin{proposition}\label{prop-fixed-slope-composition}
Let \(\sigma,\tau>0\) and \(0<p\leq q\leq\infty\).

\begin{enumerate}
\item[\textnormal{(a)}]
For every \(a\in\R\) with \(0<\sigma|a|\leq\tau\), the set
\(\operatorname{CO}_a\) is path connected. If \(q=\infty\), then
\(\operatorname{CO}_0\) is also path connected.

\item[\textnormal{(b)}]
Let $\varphi(z)=az+b$, $\phi(z)=cz+d$, where \(a,c\in\R\), \(a\neq c\), and both \(C_\varphi\) and \(C_\phi\)
belong to \(\operatorname{CO}_{p,q}(\sigma,\tau)\). Then
\[
    \|C_\varphi-C_\phi\|_{B_\sigma^p\to B_\tau^q}
    \gtrsim
    \tau^{-1/q}\Phi_{\sigma,p}(|\im b|),
\]
where the implicit constant depends only on \(p\) and \(q\).
\end{enumerate}
\end{proposition}

\begin{proof}
For part~(a), fix \(b_0,b_1\in\C\) and put
\[
    b_t:=(1-t)b_0+tb_1,
    \qquad 0\leq t\leq1.
\]
Suppose first that \(a\neq0\), and define $\Gamma(t):=C_{az+b_t}.$ By Corollary~\ref{cor-composition-classification}, every \(\Gamma(t)\) belongs
to \(\operatorname{CO}_a\).

\smallskip
\noindent
Fix \(t_0\in[0,1]\). For \(f\in B_\sigma^p\), put $g:=
    \bigl(S_{b_t-b_{t_0}}-I\bigr)f.$ Then
\[
    \bigl(\Gamma(t)-\Gamma(t_0)\bigr)f(x)
    =
    g(ax+b_{t_0}).
\]
Hence, by a change of variables and \eqref{eq:vertical-translation},
\[
    \|\bigl(\Gamma(t)-\Gamma(t_0)\bigr)f\|_{L^q}
    \lesssim_{p,q}
    |a|^{-1/q}
    \sigma^{1/p-1/q}
    e^{\sigma|\im b_{t_0}|}
    \|g\|_{L^p}.
\]
Lemma~\ref{lem-translation-continuity} therefore gives
\[
    \|\Gamma(t)-\Gamma(t_0)\|_{B_\sigma^p\to B_\tau^q}
    \longrightarrow0
    \qquad (t\to t_0).
\]
Thus \(\Gamma\) is a continuous path joining
\(C_{az+b_0}\) and \(C_{az+b_1}\).

\smallskip
\noindent
If \(q=\infty\) and \(a=0\), then
\[
    \|C_{b_t}-C_{b_{t_0}}\|_{B_\sigma^p\to B_\tau^\infty}
    =
    \|\delta_{b_t}-\delta_{b_{t_0}}\|_{(B_\sigma^p)^*}.
\]
Since
\[
    f(b_t)-f(b_{t_0})
    =
    \bigl((S_{b_t-b_{t_0}}-I)f\bigr)(b_{t_0}),
\]
Theorem~\ref{thm:point-evaluations} and
Lemma~\ref{lem-translation-continuity} show that the last norm tends to
zero as \(t\to t_0\). Hence \(\operatorname{CO}_0\) is path connected.

\smallskip
\noindent
For part~(b), for every \(x\in\R\),
\[
    \delta_x\circ(C_\varphi-C_\phi)
    =
    \delta_{\varphi(x)}-\delta_{\phi(x)}.
\]
By Theorem~\ref{thm:point-evaluations},
\[
    \|\delta_{\varphi(x)}-\delta_{\phi(x)}\|_{(B_\sigma^p)^*}
    \lesssim_q
    \tau^{1/q}
    \|C_\varphi-C_\phi\|_{B_\sigma^p\to B_\tau^q}.
\]
On the other hand, Lemma~\ref{lem-two-point-separation} gives
\[
    \|\delta_{\varphi(x)}-\delta_{\phi(x)}\|_{(B_\sigma^p)^*}
    \geq
    d_{\sigma,p}(\varphi(x),\phi(x))
\]
and hence
\[
    d_{\sigma,p}(\varphi(x),\phi(x))
    \gtrsim_p
    \Phi_{\sigma,p}(|\im b|)
    \frac{\sigma|\varphi(x)-\phi(x)|}
         {1+\sigma|\im b|+\sigma|\varphi(x)-\phi(x)|}.
\]
Since \(a\neq c\),
\[
    |\varphi(x)-\phi(x)|
    =
    |(a-c)x+b-d|
    \longrightarrow\infty
    \qquad (|x|\to\infty).
\]
Letting \(|x|\to\infty\) gives
\[
    \|C_\varphi-C_\phi\|_{B_\sigma^p\to B_\tau^q}
    \gtrsim
    \tau^{-1/q}\Phi_{\sigma,p}(|\im b|),
\]
as required.
\end{proof}

\begin{theorem}\label{thm-components-composition}
Let \(\sigma,\tau>0\) and \(0<p,q\leq\infty\).

\begin{enumerate}
\item[(a)]
If \(0<q<p\leq\infty\), then
\(
    \operatorname{CO}_{p,q}(\sigma,\tau)=\varnothing.
\)

\item[(b)]
If \(0<p\leq q<\infty\), then
\[
    \operatorname{CO}_{p,q}(\sigma,\tau)
    =
\bigsqcup_{\substack{a\in\R\\0<\sigma|a|\leq\tau}}
    \operatorname{CO}_a.
\]

\item[(c)]
If \(0<p\leq q=\infty\), then
\[
    \operatorname{CO}_{p,\infty}(\sigma,\tau)
    =
\bigsqcup_{\substack{a\in\R\\\sigma|a|\leq\tau}}
    \operatorname{CO}_a.
\]
\end{enumerate}
In cases {\rm (b)} and {\rm (c)}, each \(\operatorname{CO}_a\) is both a connected component and a
path-connected component of
\(\operatorname{CO}_{p,q}(\sigma,\tau)\).
\end{theorem}

\begin{proof}
Part~(a) follows from Corollary~\ref{cor-composition-classification}
and Remark~\ref{rem-constant-composition}.

\smallskip
\noindent
Suppose now that \(0<p\leq q\leq\infty\). The decompositions in
(b) and (c) follow from the same two results and from the
definition of the classes \(\operatorname{CO}_a\). By
Proposition~\ref{prop-fixed-slope-composition}(a), each
\(\operatorname{CO}_a\) is path connected.

\smallskip
\noindent
It remains to show that different fixed-slope classes lie in different
connected components. Fix $C_{az+b}\in\operatorname{CO}_a.$ If \(c\neq a\) and
\(C_{cz+d}\in\operatorname{CO}_c\), then
Proposition~\ref{prop-fixed-slope-composition}(b) gives
\[
    \|C_{az+b}-C_{cz+d}\|_{B_\sigma^p\to B_\tau^q}
    \gtrsim
    \tau^{-1/q}\Phi_{\sigma,p}(|\im b|)>0,
\]
with a lower bound independent of \(c\) and \(d\). Hence
\(\operatorname{CO}_a\) is relatively open in
\(\operatorname{CO}_{p,q}(\sigma,\tau)\). Since the fixed-slope classes
form a partition, each of them is also relatively closed.

\smallskip
\noindent
Therefore no connected subset can meet two distinct classes
\(\operatorname{CO}_a\). Since each class is path connected, it is
simultaneously a connected component and a path-connected component.
\end{proof}

\subsection{The space of nonzero weighted composition operators}
\label{subsec-topology-wco}

For \(\sigma,\tau>0\) and \(0<p,q\leq\infty\), let
\[
    \operatorname{WCO}_{p,q}(\sigma,\tau)
    :=
    \left\{
        W_{\psi,\varphi}:
        W_{\psi,\varphi}:B_\sigma^p\longrightarrow B_\tau^q
        \text{ is bounded and nonzero}
    \right\},
\]
endowed with the operator norm topology, or with the operator quasi-norm
topology when \(q<1\). The exclusion of the zero operator is essential,
since otherwise every bounded weighted composition operator would be
joined to zero by the path
\[
    t\longmapsto W_{t\psi,\varphi},
    \qquad 0\leq t\leq1.
\]
Unlike composition operators, weighted composition operators have an
additional degree of freedom coming from the weight. We first examine
operators with a fixed symbol.
For an entire function \(\varphi\), set
\[
    \operatorname{WCO}_{p,q}(\sigma,\tau;\varphi)
    :=
    \left\{
        W_{\psi,\varphi}\in
        \operatorname{WCO}_{p,q}(\sigma,\tau)
    \right\}.
\]

\begin{lemma}\label{lem-wco-uniqueness-fibers}
Let \(\sigma,\tau>0\) and \(0<p,q\leq\infty\).

\begin{enumerate}
\item[\textnormal{(a)}]
If
\(
    W_{\psi,\varphi},W_{\vartheta,\phi}
    \in\operatorname{WCO}_{p,q}(\sigma,\tau),
\)
then
\(
    W_{\psi,\varphi}=W_{\vartheta,\phi}
\)
if and only if
\(\psi=\vartheta\) and \(\varphi=\phi\).

\item[\textnormal{(b)}]
For every entire function \(\varphi\), the nonempty fiber
\(\operatorname{WCO}_{p,q}(\sigma,\tau;\varphi)\)
is path connected.
\end{enumerate}
\end{lemma}
\begin{proof}
The sufficiency in (a) is immediate. Conversely, suppose that $W_{\psi,\varphi}=W_{\vartheta,\phi}.$ 

\smallskip
\noindent
Since the common operator is nonzero, neither weight is identically zero.
For every \(z\in\C\),
\[
    \psi(z)\delta_{\varphi(z)}
    =
    \vartheta(z)\delta_{\phi(z)}
\]
as functionals on \(B_\sigma^p\).

\smallskip
\noindent
Let $\Omega:=\{z\in\C:\psi(z)\neq0\}.$ Then \(\Omega\) is nonempty and open. For \(z\in\Omega\), the preceding
identity shows that $\delta_{\varphi(z)}\in\C\delta_{\phi(z)}.$ Hence
\[
    d_{\sigma,p}(\varphi(z),\phi(z))=0,
\]
and Lemma~\ref{lem-two-point-separation} gives
\[
    \varphi(z)=\phi(z),\qquad z\in\Omega.
\]
Returning to the functional identity, we obtain
\[
    \psi(z)=\vartheta(z),\qquad z\in\Omega.
\]
The identity theorem therefore yields $\varphi=\phi$ and $\psi=\vartheta.$

\smallskip
\noindent
For (b), let
\[
    W_{\psi_0,\varphi},
    W_{\psi_1,\varphi}
    \in
    \operatorname{WCO}_{p,q}(\sigma,\tau;\varphi).
\]
If \(\psi_0\) and \(\psi_1\) are linearly independent, set
\[
    \psi_t:=(1-t)\psi_0+t\psi_1,
    \qquad 0\leq t\leq1.
\]
Then \(\psi_t\not\equiv0\) for every \(t\), and
\[
    W_{\psi_t,\varphi}
    =
    (1-t)W_{\psi_0,\varphi}
    +
    tW_{\psi_1,\varphi}.
\]
Thus \(t\mapsto W_{\psi_t,\varphi}\) is a continuous path in
\(\operatorname{WCO}_{p,q}(\sigma,\tau;\varphi)\).

\smallskip
\noindent
If \(\psi_1=\lambda\psi_0\) for some
\(\lambda\in\C\setminus\{0\}\), choose a continuous path $\gamma:[0,1]\longrightarrow\C\setminus\{0\}$ with \(\gamma(0)=1\) and \(\gamma(1)=\lambda\). Then
\[
    t\longmapsto W_{\gamma(t)\psi_0,\varphi}
\]
is a continuous path in the same fiber joining
\(W_{\psi_0,\varphi}\) to \(W_{\psi_1,\varphi}\).
Hence every nonempty fixed-symbol fiber is path connected.
\end{proof}

\smallskip
\noindent
To connect operators whose affine symbols vary within the strict bandwidth
range, we need a single nonzero weight that remains admissible along a
continuous path of symbols. The next lemma provides such a bridge.

\begin{lemma}\label{lem-bridge-weight}
Let \(\sigma,\tau>0\) and \(0<p,q\leq\infty\). Suppose that
\[
    a:[0,1]\to\R,
    \qquad
    b:[0,1]\to\C
\]
are continuous and
\[
    A:=\max_{0\leq t\leq1}|a(t)|,
    \qquad
    \sigma A<\tau.
\]
For \(0\leq t\leq1\), put
\[
    \varphi_t(z):=a(t)z+b(t).
\]
Choose
\(
    0<\delta<\tau-\sigma A
\)
and an integer \(m>\frac{1}{q}\), and set
\[
    h:=G_{\delta,m,0}.
\]
Then
\[
    W_{h,\varphi_t}\in
    \operatorname{WCO}_{p,q}(\sigma,\tau),
    \qquad 0\leq t\leq1,
\]
and the map
\(
    [0,1]\ni t\longmapsto W_{h,\varphi_t}
\)
is continuous.
\end{lemma}

\begin{proof}
By the packet estimates, \(h\) has exponential type at most \(\delta\),
belongs to \(L^q(\R)\) when \(q<\infty\), and satisfies
\[
    h(x)\longrightarrow0
    \qquad (|x|\to\infty).
\]
Moreover, \(h\in L^\infty(\R)\).

\smallskip
\noindent
We first verify boundedness. Suppose that \(0<p\leq q\leq\infty\).
If \(a(t)\neq0\), then
\[
    \delta<\tau-\sigma A
    \leq\tau-\sigma|a(t)|,
\]
so
\[
    h\in B^\infty_{\tau-\sigma|a(t)|}.
\]
Hence \(W_{h,\varphi_t}\) is bounded by
Theorem~\ref{thm-wco-p-le-q-bounded}. If \(a(t)=0\), then
\(\varphi_t\) is constant, and $h\in B_\tau^q$ because \(\delta<\tau\); Proposition~\ref{prop-wco-constant-symbol}
therefore applies.

\smallskip
\noindent
Now let \(0<q<p\leq\infty\). If \(p<\infty\), then
\(
    \frac{p-q}{pq}<\frac1q<m,
\)
and the packet estimates give
\[
    h\in B_\delta^{\frac{pq}{p-q}}
    \subset
    B_{\tau-\sigma|a(t)|}^{\frac{pq}{p-q}}.
\]
For \(p=\infty\), the same conclusion reads
\[
    h\in B_\delta^q
    \subset B_{\tau-\sigma|a(t)|}^q.
\]
Thus Theorem~\ref{thm-wco-q-lt-p} gives boundedness whenever
\(a(t)\neq0\), while the constant-symbol case again follows from
Proposition~\ref{prop-wco-constant-symbol}. Since \(h\not\equiv0\),
each \(W_{h,\varphi_t}\) is nonzero.

\smallskip
\noindent
It remains to prove continuity. Put
\[
    B:=\max_{0\leq t\leq1}|\im b(t)|.
\]
By Theorem~\ref{thm:point-evaluations},
\[
    |f(\varphi_t(x))|
    \lesssim_{p,\sigma,B}
    \|f\|_{L^p},
    \qquad
    x\in\R,\quad 0\leq t\leq1.
\]
Fix \(R>0\). The convex hull of $\{\varphi_t(x):|x|\leq R,\ 0\leq t\leq1\}$ is compact. Hence Theorem~\ref{thm:point-evaluations} and Cauchy's
estimate give
\[
    \sup_{|x|\leq R}
    |f(\varphi_t(x))-f(\varphi_s(x))|
    \leq C_R
    \bigl(R|a(t)-a(s)|+|b(t)-b(s)|\bigr)
    \|f\|_{L^p},
\]
where \(C_R\) is independent of \(s,t\) and \(f\).

\smallskip
\noindent
Suppose first that \(q<\infty\). By the preceding two estimates,
\[
\begin{aligned}
    \|(W_{h,\varphi_t}-W_{h,\varphi_s})f\|_{L^q}^q
    &\lesssim C_R
    \bigl(R|a(t)-a(s)|+|b(t)-b(s)|\bigr)^q
    \|h\|_{L^q([-R,R])}^q\|f\|_{L^p}^q
    \\
    &\quad+
    \|h\|_{L^q(\{|x|>R\})}^q\|f\|_{L^p}^q .
\end{aligned}
\]
The second term can be made arbitrarily small by choosing \(R\) large,
uniformly in \(s,t\), while for fixed \(R\) the first tends to zero as
\(t\to s\). Hence
\[
    \|W_{h,\varphi_t}-W_{h,\varphi_s}\|_{B_\sigma^p\to B_\tau^q}
    \longrightarrow0
    \qquad (t\to s).
\]
If \(q=\infty\), the same argument gives
\[
\begin{aligned}
    \|(W_{h,\varphi_t}-W_{h,\varphi_s})f\|_{L^\infty}
    &\lesssim_R
    \|h\|_{L^\infty}
    \bigl(R|a(t)-a(s)|+|b(t)-b(s)|\bigr)\|f\|_{L^p}
    \\
    &\quad+
    \sup_{|x|>R}|h(x)|\,\|f\|_{L^p}.
\end{aligned}
\]
Since \(h(x)\to0\) as \(|x|\to\infty\), the same two-step argument proves
continuity. 
\end{proof}

\smallskip
\noindent
By Proposition~\ref{prop-wco-constant-symbol} and
Theorems~\ref{thm-wco-p-le-q-bounded} and~\ref{thm-wco-q-lt-p},
the symbol of every operator in
\(\operatorname{WCO}_{p,q}(\sigma,\tau)\) can be written as
\[
    \varphi(z)=az+b,
    \qquad a\in\R,\quad \sigma|a|\leq\tau,\quad b\in\C,
\]
where \(a=0\) corresponds to a constant symbol. Moreover, if \(q<p\),
then necessarily \(\sigma|a|<\tau\).

\smallskip
\noindent
Accordingly, set
\[
    \operatorname{WCO}_0
    :=
    \left\{
        W_{\psi,\varphi}\in\operatorname{WCO}_{p,q}(\sigma,\tau):
        \varphi(z)=az+b,\quad
        a\in\R,\quad \sigma|a|<\tau
    \right\}.
\]
When \(0<p\leq q\leq\infty\), also set
\[
    \operatorname{WCO}_\pm
    :=
    \left\{
        W_{c,\varphi}:
        \varphi(z)=\pm\frac{\tau}{\sigma}z+b,\quad
        c\in\C\setminus\{0\},\quad b\in\C
    \right\}.
\]
Indeed, Theorem~\ref{thm-wco-p-le-q-bounded} shows that at the endpoint
slopes the weight belongs to \(B_0^\infty=\C\).

\begin{proposition}\label{prop-wco-component-families}
Let \(\sigma,\tau>0\) and \(0<p,q\leq\infty\).

\begin{enumerate}
\item[\textnormal{(a)}]
The set \(\operatorname{WCO}_0\) is path connected.

\item[\textnormal{(b)}]
If \(0<p\leq q\leq\infty\), then
\(\operatorname{WCO}_+\) and \(\operatorname{WCO}_-\) are path connected.
For either fixed choice of sign, let
\[
    \varphi_\pm(z)=\pm\frac{\tau}{\sigma}z+b,
    \qquad
    W_{c,\varphi_\pm}\in\operatorname{WCO}_\pm,
\]
and suppose that
\[
    W_{\psi,\varphi}
    \in
    \operatorname{WCO}_{p,q}(\sigma,\tau)
    \setminus\operatorname{WCO}_\pm.
\]
Then
\[
    \|W_{c,\varphi_\pm}-W_{\psi,\varphi}\|_{B_\sigma^p\to B_\tau^q}
    \gtrsim
    |c|\,\tau^{-\frac{1}{q}}\Phi_{\sigma,p}(|\im b|),
\]
where the implicit constant depends only on \(p\) and \(q\).
Consequently,
\(\operatorname{WCO}_0\), \(\operatorname{WCO}_+\), and
\(\operatorname{WCO}_-\) are relatively open and relatively closed in
\(\operatorname{WCO}_{p,q}(\sigma,\tau)\).
\end{enumerate}
\end{proposition}

\begin{proof}
For (a), let $W_{\psi_j,\varphi_j}\in\operatorname{WCO}_0,
    \,j=0,1,$ where
\[
    \varphi_j(z)=a_jz+b_j,
    \qquad
    a_j\in\R,\quad \sigma|a_j|<\tau.
\]
Put
\[
    a(t):=(1-t)a_0+ta_1,
    \qquad
    b(t):=(1-t)b_0+tb_1,
    \qquad 0\leq t\leq1.
\]
Then
\[
    \sigma\max_{0\leq t\leq1}|a(t)|
    \leq
    \sigma\max\{|a_0|,|a_1|\}<\tau.
\]
By Lemma~\ref{lem-bridge-weight}, there exists a nonzero entire function
\(h\) such that $t\longmapsto W_{h,a(t)z+b(t)}$ is a continuous path in \(\operatorname{WCO}_0\).
By Lemma~\ref{lem-wco-uniqueness-fibers}(b),
\(W_{\psi_j,\varphi_j}\) and \(W_{h,\varphi_j}\) can be joined within
their fixed-symbol fibers. Concatenating these paths proves that
\(\operatorname{WCO}_0\) is path connected.

\smallskip
\noindent
Now assume \(0<p\leq q\leq\infty\), and fix one of the two endpoint
families \(\operatorname{WCO}_\pm\). Let
\[
    W_{c_j,\pm\frac{\tau}{\sigma}z+b_j}
    \in\operatorname{WCO}_\pm,
    \qquad j=0,1,
\]
where the same sign is fixed throughout. By
Lemma~\ref{lem-wco-uniqueness-fibers}(b), we may first join the
two weights while keeping the symbol $\pm\frac{\tau}{\sigma}z+b_0$ fixed. Since \(W_{c_1,\varphi}=c_1C_\varphi\), Proposition~\ref{prop-fixed-slope-composition}(a) joins $W_{c_1,\pm\frac{\tau}{\sigma}z+b_0}$ and $W_{c_1,\pm\frac{\tau}{\sigma}z+b_1}$ by a continuous path in \(\operatorname{WCO}_\pm\). Hence \(\operatorname{WCO}_\pm\) is path connected.

\smallskip
\noindent
We next prove the separation estimate. Let
\[
    \varphi_\pm(z)=\pm\frac{\tau}{\sigma}z+b,
    \qquad
    W_{c,\varphi_\pm}\in\operatorname{WCO}_\pm,
\]
and let
\[
    W_{\psi,\varphi}
    \in
    \operatorname{WCO}_{p,q}(\sigma,\tau)
    \setminus\operatorname{WCO}_\pm.
\]
By the boundedness classification, $\varphi(z)=az+d$ for some \(a\in\R\) and \(d\in\C\), and $a\neq\pm\frac{\tau}{\sigma},$ where the sign is the fixed one above. For every \(x\in\R\),
\[
    \delta_x\circ
    \bigl(W_{c,\varphi_\pm}-W_{\psi,\varphi}\bigr)
    =
    c\,\delta_{\varphi_\pm(x)}
    -
    \psi(x)\delta_{\varphi(x)}.
\]
Hence Theorem~\ref{thm:point-evaluations} gives
\[
    \left\|
        c\,\delta_{\varphi_\pm(x)}
        -
        \psi(x)\delta_{\varphi(x)}
    \right\|_{(B_\sigma^p)^*}
    \lesssim_q
    \tau^{1/q}
    \|W_{c,\varphi_\pm}-W_{\psi,\varphi}\|_{B_\sigma^p \to B_\tau^q}.
\]
On the other hand,
\[
    \left\|
        c\,\delta_{\varphi_\pm(x)}
        -
        \psi(x)\delta_{\varphi(x)}
    \right\|_{(B_\sigma^p)^*}
    \geq
    |c|\,d_{\sigma,p}(\varphi_\pm(x),\varphi(x)).
\]
By Lemma~\ref{lem-two-point-separation},
\[
    d_{\sigma,p}(\varphi_\pm(x),\varphi(x))
    \gtrsim_p
    \Phi_{\sigma,p}(|\im b|)
    \frac{\sigma|\varphi_\pm(x)-\varphi(x)|}
         {1+\sigma|\im b|+\sigma|\varphi_\pm(x)-\varphi(x)|}.
\]
Since the two slopes are distinct,
\[
    |\varphi_\pm(x)-\varphi(x)|
    \longrightarrow\infty
    \qquad (|x|\to\infty).
\]
Letting \(|x|\to\infty\), we obtain
\[
    \|W_{c,\varphi_\pm}-W_{\psi,\varphi}\|_{B_\sigma^p\to B_\tau^q}
    \gtrsim
    |c|\,\tau^{-1/q}\Phi_{\sigma,p}(|\im b|).
\]
Thus \(\operatorname{WCO}_\pm\) is relatively open.

\smallskip
\noindent
We now prove relative closedness. Suppose that $W_{c_n,\pm\frac{\tau}{\sigma}z+b_n}
    \in\operatorname{WCO}_\pm$ and
\[
    W_{c_n,\pm\frac{\tau}{\sigma}z+b_n}
    \longrightarrow
    W_{\psi,\varphi}
    \in\operatorname{WCO}_{p,q}(\sigma,\tau).
\]
Choose \(x_0\in\R\) such that \(\psi(x_0)\neq0\). Then
\[
    c_n\delta_{\pm\frac{\tau}{\sigma}x_0+b_n}
    \longrightarrow
    \psi(x_0)\delta_{\varphi(x_0)}
\]
in \((B_\sigma^p)^*\). Therefore
\[
\begin{aligned}
    |\psi(x_0)|
    d_{\sigma,p}
    \left(
        \varphi(x_0),
        \pm\frac{\tau}{\sigma}x_0+b_n
    \right)
    &=
    \operatorname{dist}
    \left(
        \psi(x_0)\delta_{\varphi(x_0)},
        \C\delta_{\pm\frac{\tau}{\sigma}x_0+b_n}
    \right)  \\
    &\leq
    \left\|
        \psi(x_0)\delta_{\varphi(x_0)}
        -
        c_n\delta_{\pm\frac{\tau}{\sigma}x_0+b_n}
    \right\|
    \longrightarrow0.
\end{aligned}
\]
Lemma~\ref{lem-two-point-separation} now gives
\[
    \pm\frac{\tau}{\sigma}x_0+b_n
    \longrightarrow
    \varphi(x_0),
\]
and hence
\[
    b_n\longrightarrow
    b:=
    \varphi(x_0)\mp\frac{\tau}{\sigma}x_0.
\]
Choose \(f\in B_\sigma^p\) such that
\(f(\varphi(x_0))\neq0\). Applying the functional convergence above to
\(f\), and using \(b_n\to b\), gives
\[
    c_n\longrightarrow c:=\psi(x_0)\neq0.
\]
For every \(x\in\R\) and \(f\in B_\sigma^p\),
convergence in the operator (quasi-)norm gives
\[
    c_n
    f\left(
        \pm\frac{\tau}{\sigma}x+b_n
    \right)
    \longrightarrow
    \psi(x)f(\varphi(x)).
\]
Since \(b_n\to b\) and \(c_n\to c\), it follows that
\[
    c\,\delta_{\pm\frac{\tau}{\sigma}x+b}
    =
    \psi(x)\delta_{\varphi(x)},
    \qquad x\in\R.
\]
As \(c\neq0\), Lemma~\ref{lem-two-point-separation} yields
\[
    \varphi(x)=\pm\frac{\tau}{\sigma}x+b,
    \qquad
    \psi(x)=c,
    \qquad x\in\R.
\]
By the identity theorem,
\[
    \varphi(z)=\pm\frac{\tau}{\sigma}z+b,
    \qquad
    \psi(z)=c,
    \qquad z\in\C.
\]
Hence $W_{\psi,\varphi}\in\operatorname{WCO}_\pm,$ so \(\operatorname{WCO}_\pm\) is relatively closed.

\smallskip
\noindent
Finally, when \(p\leq q\), the boundedness classification gives the
disjoint decomposition
\[
    \operatorname{WCO}_{p,q}(\sigma,\tau)
    =
    \operatorname{WCO}_0
    \bigsqcup
    \operatorname{WCO}_+
    \bigsqcup
    \operatorname{WCO}_-.
\]
Since \(\operatorname{WCO}_+\) and \(\operatorname{WCO}_-\) are both
relatively open and relatively closed, the same holds for
\(\operatorname{WCO}_0\). This completes the proof.
\end{proof}

\begin{theorem}\label{thm-wco-components}
Let \(\sigma,\tau>0\) and \(0<p,q\leq\infty\).

\begin{enumerate}
\item[(a)]
If \(0<q<p\leq\infty\), then \(\operatorname{WCO}_{p,q}(\sigma,\tau)
    =
    \operatorname{WCO}_0\). In particular,
\(\operatorname{WCO}_{p,q}(\sigma,\tau)\) is path connected.

\item[(b)]
If \(0<p\leq q\leq\infty\), then
\[
    \operatorname{WCO}_{p,q}(\sigma,\tau)
    =
    \operatorname{WCO}_0
    \bigsqcup
    \operatorname{WCO}_+
    \bigsqcup
    \operatorname{WCO}_-.
\]
Moreover,
\(    \operatorname{WCO}_0\), \(\operatorname{WCO}_+\), \(\operatorname{WCO}_-\)
are precisely the connected components, and also the path-connected
components, of
\(\operatorname{WCO}_{p,q}(\sigma,\tau)\).
\end{enumerate}
In all cases, \(\operatorname{WCO}_{p,q}(\sigma,\tau)\) has no isolated
points.
\end{theorem}

\begin{proof}
Suppose first that \(0<q<p\leq\infty\). By
Theorem~\ref{thm-wco-q-lt-p}, every nonconstant symbol of a bounded
nonzero weighted composition operator satisfies
\[
    \varphi(z)=az+b,
    \qquad
    a\in\R,\quad \sigma|a|<\tau.
\]
Constant symbols correspond to \(a=0\). Hence
\[
    \operatorname{WCO}_{p,q}(\sigma,\tau)
    =
    \operatorname{WCO}_0,
\]
which is path connected by
Proposition~\ref{prop-wco-component-families}(a).

\smallskip
\noindent
Now let \(0<p\leq q\leq\infty\). By
Proposition~\ref{prop-wco-constant-symbol},
Theorem~\ref{thm-wco-p-le-q-bounded}, and Lemma~\ref{lem-wco-uniqueness-fibers}(a), every operator in
\(\operatorname{WCO}_{p,q}(\sigma,\tau)\) belongs to exactly one of
\(    \operatorname{WCO}_0\), \(\operatorname{WCO}_+\), and  \(\operatorname{WCO}_-\).
Thus
\[
    \operatorname{WCO}_{p,q}(\sigma,\tau)
    =
    \operatorname{WCO}_0
    \bigsqcup
    \operatorname{WCO}_+
    \bigsqcup
    \operatorname{WCO}_-.
\]
By Proposition~\ref{prop-wco-component-families}, each of these sets is
path connected and relatively open and relatively closed. Hence no
connected subset of
\(\operatorname{WCO}_{p,q}(\sigma,\tau)\) can meet two of them.
Therefore they are precisely the connected components. Since each is
path connected, they are also precisely the path-connected components.

\smallskip
\noindent
Finally, if
\(W_{\psi,\varphi}\in\operatorname{WCO}_{p,q}(\sigma,\tau)\), then
\[
    W_{(1+1/n)\psi,\varphi}
    \longrightarrow
    W_{\psi,\varphi}
\]
in operator (quasi-)norm, while
\(W_{(1+1/n)\psi,\varphi}\neq W_{\psi,\varphi}\) for every \(n\).
Hence \(\operatorname{WCO}_{p,q}(\sigma,\tau)\) has no isolated points.
\end{proof}

\begin{remark}\label{rem-topology-comparison}
The preceding results show a sharp contrast between the topological
structures of composition and weighted composition operators.

\smallskip
\noindent
For composition operators, the slope is a complete component invariant:
each admissible fixed-slope class \(\operatorname{CO}_a\) forms a
connected and path-connected component. By contrast, all nonzero weighted composition operators with constant
symbols or with affine symbols satisfying the strict bandwidth condition $\sigma|a|<\tau$ belong to the single component \(\operatorname{WCO}_0\). When
\(0<p\leq q\leq\infty\), only the two endpoint slopes $a=\pm\frac{\tau}{\sigma}$ give rise to the additional components
\(\operatorname{WCO}_+\) and \(\operatorname{WCO}_-\). If
\(0<q<p\leq\infty\), the strict bandwidth condition excludes these
endpoint families, and the whole space
\(\operatorname{WCO}_{p,q}(\sigma,\tau)\) is path connected.
\end{remark}

\section*{Acknowledgements} 
The research is funded by the Vietnam National Foundation for Science and Technology Development (NAFOSTED) under grant number 101.02-2025.08. Part of this work was carried out during the stays of Pham Viet Hai and
Pham Trong Tien at the Vietnam Institute for Advanced Study in Mathematics.
They thank the Institute for its hospitality and support.

\bigskip
\noindent
\textbf{AI Statement.}
The authors used artificial-intelligence tools for language editing,
\LaTeX{} formatting, and limited assistance in checking local mathematical
arguments. The proof strategy and mathematical development are entirely the
authors' own. All mathematical statements, proofs, and conclusions were
independently verified by the authors, who take full responsibility for
the content of the paper.

\end{document}